\documentclass[11pt,amsfonts]{amsart}

\usepackage{fullpage}

\newtheorem{specialtheorem}{Theorem}
\def\endproof{\qed\smallskip}
\def\blacksquare{\hbox to .60em{\vrule width .60em height .60em}}

\newtheorem{theorem}{Theorem}[section]

\newtheorem{corollary}[theorem]{Corollary}

\newtheorem{definition}[theorem]{Definition}

\newtheorem{lemma}[theorem]{Lemma}

\newtheorem{proposition}[theorem]{Proposition}

\newtheorem{remark}[theorem]{Remark}

\newtheorem{example}[theorem]{Example}

\begin{document}

\title[]{Einstein Metrics with Prescribed Conformal Infinity on 
4-Manifolds}

\author[]{Michael T. Anderson}

\thanks{Partially supported by NSF Grant DMS 0604735}

\maketitle

\abstract
This paper considers the existence of conformally compact Einstein 
metrics on 4-manifolds. A reasonably complete understanding is obtained 
for the existence of such metrics with prescribed conformal infinity, 
when the conformal infinity is of positive scalar curvature. We find in 
particular that general solvability depends on the topology of the 
filling manifold. The obstruction to extending these results to 
arbitrary boundary values is also identified. While most of the paper 
concerns dimension 4, some general results on the structure of the space 
of such metrics hold in all dimensions. 
\endabstract

%\tableofcontents

\setcounter{section}{0}

\section{Introduction.}
\setcounter{equation}{0}

 This paper is concerned with the existence of conformally compact 
Einstein metrics on a given manifold $M$ with boundary $\partial M$. 
The main results are restricted to dimension 4, although some of the 
results hold in all dimensions.
  
  This existence problem was raised by Fefferman and Graham in [19] in 
connection with a study of conformal invariants of Riemannian manifolds. 
More recently, the study of such metrics has become of strong interest 
through the AdS/CFT correspondence, relating gravitational theories on 
$M$ with conformal field theories on $\partial M$, cf.~[18], [38] and 
references therein.

 Let $M$ be a compact, oriented manifold with non-empty boundary 
$\partial M$; $M$ is assumed to be connected, but apriori $\partial M$ 
may be connected or disconnected. A defining function $\rho$ for 
$\partial M$ in $M$ is a non-negative $C^{\infty}$ function on the closure 
$\bar M = M\cup\partial M$ such that $\rho^{-1}(0) = \partial M$ and 
$d\rho  \neq 0$ on $\partial M$. A complete Riemannian 
metric $g$ on $M$ is {\it  conformally compact} if there is a defining 
function $\rho$ such that the conformally equivalent metric
\begin{equation} \label{e1.1}
\widetilde g = \rho^{2}\cdot  g 
\end{equation}
extends at least continuously to a Riemannian metric $\widetilde g$ on 
$\bar M$. The metric $\gamma  = \widetilde g|_{\partial M}$ induced on 
$\partial M$ is the boundary metric induced by $g$ and the compactification 
$\rho$. Any compact manifold with boundary carries many conformally compact 
metrics; for instance, one may let $\widetilde g$ be any smooth metric on 
$\bar M$ and define $g$ by $g = \rho^{-2}\cdot \widetilde g$, for any choice 
of defining function $\rho$. 

 Defining functions for $\partial M$ are unique only up to 
multiplication by positive smooth functions on $\bar M$ and hence 
the compactification (1.1) is not uniquely determined by $(M, g)$. On 
the other hand, the conformal class $[\gamma]$ of a boundary metric is 
uniquely determined by $g$; the class $[\gamma]$ is called the {\it  
conformal infinity}  of $(M, g)$.

 A conformal compactification is called $C^{m,\alpha}$ or $L^{k,p}$ if 
the metric $\widetilde g$ extends to a $C^{m,\alpha}$ or $L^{k,p}$ metric 
on the closure $\bar M;$ here $C^{m,\alpha}$ and $L^{k,p}$ are the usual 
H\"older and Sobolev function spaces.

\medskip

 In this paper, we consider complete conformally compact Einstein 
metrics $g$ on $(n+1)$-dimensional manifolds $M$, normalized so that
\begin{equation} \label{e1.2}
Ric_{g} = -ng. 
\end{equation}
It is easy to see, cf. the Appendix, that any $C^{2}$ conformally compact 
Einstein metric $g$ satisfies $|K_{g} + 1| = O(\rho^{2})$, where $K_{g}$ denotes 
the sectional curvature of $g$. Hence, such metrics are asymptotically 
hyperbolic (AH), in that the local geometry tends to that of hyperbolic 
space at infinity.

\medskip

 Let $E_{AH} = E_{AH}^{m,\alpha}$ be the space of AH Einstein metrics 
on $M$ which admit a $C^{m,\alpha}$ compactification $\widetilde g$ as 
in (1.1). We require that $m \geq 3$, $\alpha\in (0,1)$ but otherwise 
allow any value of $m$, including $m = \infty$ or $m = \omega$, 
(for real-analytic). The space $E_{AH}^{m,\alpha}$ is given the 
$C^{m,\alpha'}$ topology on metrics on $\bar M$, for any fixed 
$\alpha' < \alpha$, via a fixed compactification as in (1.1). Let 
${\mathcal E}_{AH} = E_{AH}/{\mathcal D}_{1}(\bar M),$ where ${\mathcal D}_{1}(\bar M)$ 
is the group of orientation preserving $C^{m+1,\alpha}$ diffeomorphisms of 
$\bar M$ inducing the identity on $\partial M,$ acting on $E_{AH}$ in the 
usual way by pullback. 

  As boundary data, let $Met(\partial M) = Met^{m,\alpha}(\partial M)$ be the 
space of $C^{m,\alpha}$ metrics on $\partial M$ and ${\mathcal C}  = {\mathcal C}^{m,\alpha}
(\partial M)$ the corresponding space of pointwise conformal classes, as 
above endowed with the $C^{m,\alpha'}$ topology. There is a natural boundary map, 
(for any fixed $(m, \alpha)$),
\begin{eqnarray} \label{e1.3}
\Pi : {\mathcal E}_{AH} \rightarrow  {\mathcal C} \\
\Pi [g] = [\gamma], \nonumber
\end{eqnarray}
which takes an AH Einstein metric $g$ on $M$ to its conformal infinity 
on $\partial M$. 

  It is proved in [7] that when $dim M = 4$, and $\pi_{1}(M, \partial M) = 0$, the 
space ${\mathcal E}_{AH}^{m,\alpha}$ is a $C^{\infty}$ infinite dimensional Banach 
manifold, and the boundary map $\Pi$ is $C^{\infty}$ of Fredholm index 0, provided 
${\mathcal E}_{AH}^{m,\alpha}$ is non-empty. Moreover, the spaces 
${\mathcal E}_{AH}^{m,\alpha}$ are essentially independent of $(m, \alpha)$, 
in that these spaces are all diffeomorphic, with ${\mathcal E}_{AH}^{m',\alpha'}$ 
dense in ${\mathcal E}_{AH}^{m,\alpha}$ whenever $m'+\alpha' > m+\alpha$. 
Essentially the same results also hold in in all dimensions $n+1 \geq 4$, 
see \S 2 for a more detailed discussion. 

 The issues of existence and uniqueness of AH Einstein metrics with a 
given conformal infinity are thus equivalent to the surjectivity and 
injectivity of the boundary map $\Pi$. 

\medskip

 The main result of the paper which is used to approach the existence 
question is the following. Let ${\mathcal C}^{o}$ be the space of 
non-negative conformal classes $[\gamma]$ on $\partial M,$ in the 
sense that $[\gamma]$ has a {\it  non-flat}  representative $\gamma $ 
of non-negative scalar curvature. Let ${\mathcal E}_{AH}^{o} = 
\Pi^{-1}({\mathcal C}^{o})$ be the space of AH Einstein metrics on $M$ with 
conformal infinity in ${\mathcal C}^{o}.$ Thus, one has the restricted 
boundary map $\Pi^{o} = \Pi|_{{\mathcal E}_{AH}^{o}}: {\mathcal E}_{AH}^{o} 
\rightarrow  {\mathcal C}^{o}.$
\begin{specialtheorem} \label{t A.}
  Let $M$ be a $4$-manifold satisfying $\pi_{1}(M, \partial M) = 0$ and for 
which the inclusion $\iota: \partial M \rightarrow \bar M$ induces a surjection
\begin{equation} \label{e1.4}
H_{2}(\partial M, {\mathbb F} ) \rightarrow  H_{2}(\bar M, {\mathbb F} ) 
\rightarrow  0, 
\end{equation}
for all fields ${\mathbb F}$. Then for any $(m, \alpha)$, $m \geq 4$, 
the boundary map 
\begin{equation} \label{e1.5}
\Pi^{o}: {\mathcal E}_{AH}^{o} \rightarrow  {\mathcal C}^{o} 
\end{equation}
is proper.
\end{specialtheorem}

 Unfortunately (or fortunately), the boundary map $\Pi$ is not proper 
in general. In fact, a sequence of distinct AH Einstein metrics $\{g_{i}\}$ 
on $M = {\mathbb R}^{2}\times T^{2}$ was constructed in [6], whose conformal 
infinity is an arbitrary fixed flat metric on $\partial M = T^{3}$, but 
which has no convergent subsequence to an AH Einstein metric on $M$. 
These metrics are twisted versions of the AdS $T^{2}$ black hole metrics, 
cf.~Remark 2.6. The sequence $\{g_{i}\}$ converges, (in a natural sense), 
to a complete hyperbolic cusp metric
\begin{equation} \label{e1.6}
g_{C} = dr^{2} + e^{2r}g_{T^{3}} 
\end{equation}
on the manifold $N = {\mathbb R} \times T^{3}.$ Since the boundary metric 
on $T^{3}$ is flat and $M$ satisfies (1.4), Theorem A is sharp, at least 
without further restrictions.

 We will show in Theorem 5.4 that this behavior is the only way that 
$\Pi$ is non-proper, in that divergent sequences $\{g_{i}\}$ of AH 
Einstein metrics on a fixed 4-manifold $M$ and with a fixed conformal 
infinity not in ${\mathcal C}^{o}$, (or a convergent sequence of conformal 
infinities outside ${\mathcal C}^{o}$), necessarily converge to AH Einstein 
manifolds $(N, g)$ with cusps, although the cusps are not necessarily of 
the simple form (1.6). Thus, the completion $\bar {\mathcal E}_{AH}$ of 
${\mathcal E}_{AH}$ in a natural topology consists of AH Einstein metrics on 
$M$, together with AH Einstein metrics with cusps, cf.~Theorem 5.4 and 
Corollary 5.5 for further details. Using this, we show in Corollary 5.6 
that the boundary map $\Pi$ is proper onto the subset $\hat {\mathcal C}$ of 
${\mathcal C}$ consisting of conformal classes which are not the boundary 
metrics of AH Einstein metrics with cusps.

 The infinite sequence of AH Einstein metrics $g_{i}$ on ${\mathbb R}^{2} \times T^{2}$ 
above all lie in distinct components of the space ${\mathcal E}_{AH}$, so that this 
space has infinitely many components. On the other hand, the fact that $\Pi^{o}$ is 
proper of course implies that ${\mathcal E}_{AH}^{o}$ has only finitely many components 
mapping into any compact region of ${\mathcal C}^{o}$. 

 The topological condition (1.4) is used solely to rule out orbifold 
degenerations of AH Einstein sequences $\{g_{i}\}.$ This condition may 
not be necessary, but it has been included in order to not overly 
complicate the paper; some situations where it is not needed are 
described in \S 7, cf.~Remark 7.4. Similarly it may well be possible 
to weaken or drop the condition $\pi_{1}(M, \partial M) = 0$ on the 
fundamental group, but this will not be pursued further here.

\medskip

   Theorem A implies that at least on ${\mathcal E}_{AH}^{o}$, the map $\Pi^{o}$ is a 
proper Fredholm map of index 0. It follows from a result of Smale [34] that $\Pi^{o}$ 
has a well-defined mod 2 degree, $deg_{2}\Pi^{o} \in  {\mathbb Z}_{2}$, on each component 
of ${\mathcal E}_{AH}^{o}$. In fact, building on the work of Tromba [35] and White 
[36], [37] we show in Theorem 6.1 that $\Pi^{o}$ has a ${\mathbb Z}$-valued degree 
\begin{equation} \label{e1.7}
deg\Pi^{o} \in  {\mathbb Z} ,
\end{equation}
again on components of ${\mathcal E}_{AH}^{o}$. Tautologically, 
if $deg\Pi^{o} \neq 0$, then $\Pi^{o}$ is surjective.

 We compute $deg\Pi^{o}$ in certain cases via symmetry arguments. 
Namely, it is proved in [8], (cf.~also Theorem 2.5 below) that any connected 
group of isometries of the conformal infinity $[\gamma]$ on $\partial M$ 
extends to a group of isometries of any AH Einstein filling metric 
$(M, g)$, with $\Pi [g] = [\gamma]$, in any dimension, again provided 
$\pi_{1}(M, \partial M) = 0$. Using this, and a straightforward classification 
of Einstein metrics with large symmetry groups leads to the following result.

\begin{specialtheorem} \label{t B.}
  Let $M = B^{4}$ be the $4$-ball, with $\partial M = S^{3}$, and, 
(abusing notation slightly), let ${\mathcal C}^{o}$ be the component of the 
non-negative $C^{m,\alpha}$ conformal classes containing the round metric 
on $S^{3}$. Also let ${\mathcal E}_{AH}^{o}$ be the component of 
$\Pi^{-1}({\mathcal C}^{o})$ containing the Poincar\'e metric on $B^{4}$. 
Then
\begin{equation} \label{e1.8}
deg_{B^{4}}\Pi^{o} = 1. 
\end{equation}
In particular, for any $(m, \alpha)$, $m \geq 4$, any conformal class 
$[\gamma] \in {\mathcal C}^{o}$ on $S^{3}$ is the conformal infinity of 
an AH Einstein metric on $B^{4}$. 
\end{specialtheorem}

 Explicit examples of AH Einstein metrics on $B^{4}$, whose the 
conformal infinity is an arbitrary Berger sphere, were constructed by 
Pedersen [32]. More generally, Hitchin [26] has constructed such 
metrics with conformal infinity an arbitrary left-invariant metric on 
$S^{3}$. 

\smallskip

 It remains an open question whether the boundary map $\Pi$ is 
surjective onto all of ${\mathcal C}$ when $M = B^{4}$. Based on reasoning 
from the AdS/CFT correspondence, Witten [38] has remarked that this may 
not be the case. It appears that physically, the natural class of 
boundary metrics are those of positive scalar curvature; in this 
regard, see also the work of Witten-Yau [39] proving that $\partial M$ 
is necessarily connected when the conformal infinity has a component of 
positive scalar curvature. In any case, as remarked above in connection 
with Theorem 5.4, the only obstruction to surjectivity is the possible 
existence of AH Einstein cusp metrics associated to $B^{4}$, cf.~also 
Remark 6.8. 

\medskip

 Theorem B also holds when $M$ is a disc bundle over $S^{2}$, 
so that $\partial M = S^{3}/{\mathbb Z}_{k}$, provided $k \geq 10$, with again 
${\mathcal C}^{o}$ the component containing the standard round metric on 
$S^{3}/{\mathbb Z}_{k}$, cf.~Proposition 7.6. We conjecture this result holds 
for any $k \geq 3$. However, for $M = S^{2}\times {\mathbb R}^{2}$, (i.e.~$k = 0$), 
it is shown in Proposition 7.2, that
\begin{equation} \label{e1.9}
deg_{{\mathbb R}^{2}\times S^{2}}\Pi^{o} = 0, 
\end{equation}
and further that $\Pi^{o}$ is not surjective onto ${\mathcal C}^{o}$. In 
fact, the conformal class of $S^{2}(1)\times S^{1}(L),$ for any $L >  
2\pi / \sqrt{3} $ is not the boundary metric of any AH Einstein metric 
on $S^{2}\times {\mathbb R}^{2}$. This result is based on the work of 
Hawking-Page [23] on the AdS-Schwarzschild metric. Similarly, for $M$ 
the degree 1 disc bundle over $S^{2},$ i.e.~$M = {\mathbb C}{\mathbb P}^{2}
\setminus B^{4}$, as well as the degree 2 disc bundle, $\Pi^{o}$ is not 
surjective. These results for $\Pi^{o}$ on disc bundles use properties of 
the AdS-Taub Bolt metrics analysed in [24],[31]. These results and others 
on deg $\Pi^{o}$ are given in \S 7.

  Theorem B and the related results above represent the first {\it global} results 
on an existence theory for complete Einstein metrics, outside the context of 
K\"ahler-Einstein metrics.

  I would like to thank Piotr Chru\'sciel, Robin Graham, Rafe Mazzeo, 
Michael Singer for their interest and comments on this work. Thanks also to 
the referees for a careful examination of the paper leading to improvements in the 
exposition.

\section{Background Results}
\setcounter{equation}{0}

  In this section, we discuss a number of background results from [6]-[8], [11] 
[19], [21] needed for the work to follow.

  We begin with the structure of the moduli space ${\mathcal E}_{AH}^{m,\alpha}$ of 
$C^{m,\alpha}$ conformally compact Einstein metrics on a given $(n+1)$-dimensional 
manifold $M$. The results discussed below are from [7]. 
  
  Let $E_{AH}^{m,\alpha}$ be the space of AH Einstein metrics on $M$ which have 
a $C^{m,\alpha}$ conformal compactification with respect to a smooth (or $C^{\omega}$) 
defining function $\rho$, as in (1.1). The space $E_{AH}^{m,\alpha}$ is given the 
$C^{m,\alpha'}$ topology on the space of $C^{m,\alpha}$ metrics $Met^{m,\alpha}(\bar M)$ 
on $\bar M$, for some $\alpha' < \alpha$. Thus, a neighborhood of $g \in 
{E}_{AH}^{m,\alpha}$ is defined to be the set of metrics $g'$ whose compactification 
$\widetilde g'$ by $\rho$ lie in a neighborhood of $\widetilde g$ in $Met^{m,\alpha}(\bar M)$, 
where $Met^{m,\alpha}(\bar M)$ is given the $C^{m,\alpha'}$ topology. Let $Met^{m,\alpha}$ 
be the corresponding space of metrics on $\partial M$, again with the $C^{m,\alpha'}$ 
topology. One thus has a natural boundary map
\begin{equation} \label{e2.1}
\Pi = \Pi^{m,\alpha}: E_{AH}^{m,\alpha} \rightarrow Met^{m,\alpha}(\partial M),
\end{equation}
sending $g$ to its boundary metric $\gamma$. 

  An element $[g]$ in the moduli space ${\mathcal E}_{AH}^{m,\alpha}$ is an equivalence 
class of a $C^{m,\alpha}$ conformally compact Einstein metric $g$ on $M$, where 
$g \sim g'$ if $g' = \phi^{*}g$, with $\phi$ a $C^{m+1,\alpha}$ diffeomorphism of $\bar M$, 
equal to the identity on $\partial M$, i.e. $\phi \in {\mathcal D}_{1}$. Changing $g$ by 
a diffeomorphism in ${\mathcal D}_{1}$, or changing the defining function $\rho$ in (1.1) 
changes the boundary metric conformally. Hence if ${\mathcal C}^{m,\alpha}$ denotes 
the space of conformal classes of $C^{m,\alpha}$ metrics on $\partial M$, then 
the boundary map
\begin{equation} \label{e2.2}
\Pi = \Pi^{m,\alpha}: {\mathcal E}_{AH}^{m,\alpha} \rightarrow {\mathcal C}^{m,\alpha}
\end{equation}
is well defined. When there is no danger of confusion, we will usually work with 
any given representative $g \in [g]$ and $\gamma \in [\gamma]$.

  Suppose $dim M = n+1 = 4$, and the 4-manifold $M$ satisfies $\pi_{1}(M, \partial M) = 0$. 
Then if ${\mathcal E}_{AH}^{m,\alpha}$ is non-empty, the spaces $E_{AH}^{m,\alpha}$ and 
${\mathcal E}_{AH}^{m,\alpha}$ are $C^{\infty}$ infinite dimensional separable Banach 
manifolds, and the boundary maps $\Pi$ in (2.1) and (2.2) are $C^{\infty}$ smooth maps 
of Banach manifolds, of Fredholm index $0$. This result holds for any $m \geq 3$, 
and $\alpha \in (0,1)$, and, replacing Banach by Fr\'echet, also for $m = \infty$ and 
$m = \omega$, (i.e. the real-analytic case). Implicit in this statement is the following 
boundary regularity result: an AH Einstein metric with a $C^{2}$ conformal compactification 
which has a $C^{m,\alpha}$ boundary metric, $m \geq 2$, has a $C^{m,\alpha}$ conformal 
compactification. 

  In addition, the spaces ${\mathcal E}_{AH}^{m,\alpha}$ are all diffeomorphic and for any 
$(m', \alpha')$ with $m'+ \alpha ' > m+\alpha$, the space ${\mathcal E}_{AH}^{m',\alpha'}$ 
is dense in ${\mathcal E}_{AH}^{m,\alpha}$. In particular,
\begin{equation} \label{e2.3} 
\overline{{\mathcal E}_{AH}^{\omega}} = {\mathcal E}_{AH}^{m, \alpha},
\end{equation}
where the completion is taken in the $C^{m,\alpha'}$ topology.

  If $dim M = n+1 > 4$, then these results hold when $3 \leq m \leq n-1$, using 
the boundary regularity results in [30]. In addition, the smooth manifold result 
holds when $m = \infty$, if, when $n$ is even, ${\mathcal E}_{AH}^{\infty}$ is understood 
to be the space of AH Einstein metrics which are $C^{\infty}$ polyhomogeneous; this 
uses the boundary regularity result of [16], cf.~also the discussion following 
(2.10) below. 

  Let ${\mathbb S}^{m,\alpha}(\bar M)$ denote the space of $C^{m,\alpha}$ symmetric 
bilinear forms on $\bar M$, and similarly for $M$ and $\partial M$. A 
tangent vector $h$ to the Banach manifold ${\mathcal E}_{AH}^{m,\alpha}$ at a 
representative $g \in [g]$, i.e. an infinitesimal Einstein deformation of $g$, 
is a form $h \in {\mathbb S}^{m,\alpha}(\bar M)$ satisfying the equation
\begin{equation} \label{e2.4}
L(h) = \tfrac{1}{2}D^{*}Dh - R(h) = 0, \ \ \beta_{g}(h) = 0, 
\end{equation}
where $R$ is the action of the curvature tensor of $g$ on symmetric bilinear forms 
and $\beta_{g} = \delta_{g} + \frac{1}{2}dtr_{g}$ is the Bianchi operator with 
respect to $g$. Such $h$ are transverse to the orbits of the diffeomorphism group 
${\mathcal D}_{1}$ acting on $\bar M$. Let $T_{g}{\mathcal E}_{AH}^{m,\alpha}$ denote the 
space of such tangent vectors.

  An Einstein metric $(M, g)$ is {\it non-degenerate} if the operator $L$ has no 
kernel in $L^{2}(M, g)$, 
\begin{equation} \label{e2.5}
K = L^{2}-Ker L = 0.
\end{equation}
For $g \in {\mathcal E}_{AH}^{m,\alpha}$, the $L^{2}$ kernel $K$ equals the kernel of the 
linear map $D\Pi: T_{g}{\mathcal E}_{AH}^{m,\alpha} \rightarrow T_{\Pi(g)}{\mathcal C}^{m,\alpha}$. 
Hence, $g$ is non-degenerate if and only if $g$ is a regular point of the boundary map 
$\Pi$ and thus $\Pi$ is a local diffeomorphism near $g$. Any element $\kappa \in K$ is 
transverse-traceless, i.e. 
\begin{equation} \label{e2.6}
\delta_{g}\kappa = 0, \ \ tr_{g}\kappa = 0,
\end{equation}
and so satisfies the Bianchi gauge condition (2.4).  

\medskip

 In the next sections, we will frequently consider compactifications 
$\bar g$ of $g$ by a geodesic defining function $t$, for which
\begin{equation} \label{e2.7}
t(x) = dist_{\bar g}(x, \partial M). 
\end{equation}
A compactification $\bar g = t^{2}g$ satisfying (2.7) is called a {\it geodesic 
compactification} and $t$ is a {\it  geodesic} defining function. Such compactifications 
are natural from a number of viewpoints. In particular the curvature of $\bar g$ has a 
particularly simple form; this and related issues are discussed in the the Appendix. 
It is not difficult to see that given a $C^{2}$ conformally compact Einstein metric 
$g$ with boundary metric $\gamma$ in some compactification, there is a unique geodesic 
compactification $\bar g$ of $g$ with boundary metric $\gamma$. Further if $\widetilde g$ 
is some $C^{m,\alpha}$ compactification, then then the geodesic compactification 
$\bar g$ is at least $C^{m-1,\alpha}$ off the cutlocus $\bar C$ of $\partial M$ in 
$(M, \bar g)$.

 Since the integral curves of $\bar \nabla t$ are geodesics, the metric 
$\bar g$ splits as
\begin{equation} \label{e2.8}
\bar g = dt^{2} + g_{t}, 
\end{equation}
within a collar neighborhood $U$ (inside the cutlocus) of $\partial M$; 
here $g_{t}$ is a curve of metrics on the boundary $\partial M$ with 
$g_{0} = \gamma$. Setting $r = -\log t$, a simple calculation shows that
the integral curves of $\nabla r$ are also geodesics in $(M, g)$ and so the metric 
$g$ also splits as $g = dr^{2} + g_{r}$. 

\medskip

  The Fefferman-Graham expansion [19] is the formal expansion of the curve $g_{t}$ 
in a Taylor-like series: the specific form of the expansion depends on $n$. If 
$dim M = n+1$ is even, then the expansion reads 
\begin{equation} \label{e2.9}
g_{t} \sim g_{(0)} + t^{2}g_{(2)} + \cdots + t^{n-1}g_{(n-1)} + t^{n}g_{(n)} + \cdots ,
\end{equation}
where the coefficients are symmetric bilinear forms on $\partial M$. The expansion is 
in even powers of $t$ up to order $n-1$ and the terms $g_{(2k)}$ are intrinsically 
determined by the boundary metric $\gamma = g_{(0)}$ and its tangential derivatives up 
to order $2k$, for $2k \leq n-1$. The term $g_{(n)}$ is transverse-traceless; however, 
beyond this, $g_{(n)}$ is not locally determined by $\gamma$; it depends on global 
properties of the Einstein metric $(M, g)$. For example, when $n = 3$, one has 
$$g_{(3)} = \tfrac{1}{6} \bar \nabla_{N} \bar Ric,$$
where $\bar Ric$ is the Ricci curvature and $N = \bar \nabla t$ is the unit normal of 
$\partial M$ in $(\bar M, \bar g)$. This term is not computable in terms of $\gamma$; 
it depends on the global properties of $(M, g)$. On the other hand, the higher terms 
$g_{(k)}$, $k > n$ depend on two tangential derivatives of $g_{(k-2)}$. 

  If $dim M = n+1$ is even, then the expansion is 
\begin{equation} \label{e2.10}
g_{t} \sim g_{(0)} + t^{2}g_{(2)} + ... +  t^{n}g_{(n)} + t^{n}(\log t){\mathcal H} + \cdots ,
\end{equation}
where again this is an even expansion up to order $n$, with coefficients up to order 
$(n-2)$ locally determined by $\gamma$. The term ${\mathcal H}$, called the obstruction 
tensor in [19], is naturally identified as the stress-energy tensor of the integral $L$ 
of the conformal anomaly, (see (2.16) below), via the AdS/CFT correspondence, cf.~[18], 
[33]; it is transverse-traceless and also locally determined by $\gamma$. The trace 
and divergence of $g_{(n)}$ are determined by $\gamma$; in fact there is a symmetric 
bilinear form $r_{(n)}$ and scalar function $a_{(n)}$, both explicitly computable 
from $\gamma$ and its derivatives, such that 
\begin{equation} \label{e2.11}
\delta(g_{(n)} + r_{(n)}) = 0, \ \ tr(g_{(n)} + r_{(n)}) = a_{(n)}.
\end{equation}
The term $a_{(n)}$ is the conformal anomaly, (so that $\int_{\partial M}a_{(n)} = L$). 
However, as above, $g_{(n)}$ is otherwise only globally determined by $(M, g)$. 
The higher order terms in (2.10) contain terms of the form $t^{m}(\log t)^{p}$. The 
expansion is even in powers of $t$, and each coefficient $g_{(k)}$, $k \neq n$, 
depends on two derivatives of $g_{(k-2)}$. Note that (2.11) also holds when $n$ is 
odd, in which case $r_{(n)} = a_{(n)} = 0$. In both cases, all terms in the expansion 
are thus determined by the two terms $g_{(0)}$ and $g_{(n)}$ or $\tau_{(n)}$. 

  Let
\begin{equation} \label{e2.12}
\tau_{(n)} = g_{(n)} + r_{(n)},
\end{equation}
so that $\tau_{(n)}$ is divergence-free. This term will play an important role 
thoughout much of the paper, in particular in \S 4. In the AdS/CFT correspondence 
$\tau_{(n)}$ corresponds, (up to a constant), to the stress-energy tensor of the 
corresponding CFT on $\partial M$, cf.~(2.19) below. 

  The boundary regularity results of [6], [7], [16] discussed above give the existence 
of the formal series (2.9) and (2.10); in the $C^{\infty}$ case, these are well-defined 
asympotic series for the curve $g_{t}$. If the data $(g_{(0)}, g_{(n)})$ are 
real-analytic, a result of Kichenassamy [28] shows that the series (2.9) and (2.10) 
converge to $g_{t}$, so that such data determine a solution to the Einstein equations 
defined in a neighborhood of $\partial M$. 

  The next local unique continuation result from [8] shows that the terms $g_{(0)}$ 
and $g_{(n)}$ uniquely determine an AH Einstein metric up to isometry, near the boundary. 

\begin{proposition} \label{p2.1}
Let $g^{1}$ and $g^{2}$ be two Einstein metrics defined on a half-ball $U \simeq 
({\mathbb R}^{n+1})^{+}$, which have a $C^{3,\alpha}$ conformal compactification 
to the closed half-ball $\bar U$. If 
\begin{equation} \label{e2.13}
|g^{1} - g^{2}| = o(t^{n}),
\end{equation}
where $t$ is the distance to $\partial U \simeq {\mathbb R}^{n}$, (in either 
compactified metric), then $g^{1}$ is isometric to $g^{2}$ in $U$.  
\end{proposition}

  Of course (2.13) holds if $g^{1}$ and $g^{2}$ are $C^{n}$ polyhomogeneous conformally 
compact with identical $g_{(0)}$ and $g_{(n)}$ terms. Since Einstein metrics are 
real-analytic in local harmonic coordinates, local unique continuation in the interior, 
away from the boundary, is well-known. 

\medskip

   Observe that the terms $g_{(k)}$ depend on the choice of the boundary metric 
$\gamma \in [\gamma]$. If $\widetilde \gamma = \lambda^{2}\gamma$, then the 
coefficients $\widetilde g_{(k)}$ are determined explicitly by $\lambda$ and 
$\gamma$, cf.~[18]. For $n$ odd, the tranformation rule for the undetermined 
coefficient $g_{(n)}$ has the simple form
\begin{equation} \label{e2.14}
\widetilde g_{(n)} = \lambda^{-(n-2)}g_{(n)}. 
\end{equation}
In particular $|\widetilde g_{(n)}|_{\widetilde g} = \lambda^{-n}|g_{(n)}|_{g}$. 
For $n$ even, the transformation rule is more complicated; it has the form (2.14) 
at leading order, but has lower order terms depending explicitly on the derivatives 
of $\gamma$ and $\lambda$. Exact formulas in low dimensions are given in [18], [33]. 

\medskip

 Next we consider certain global issues associated with AH Einstein metrics which 
arise via the AdS/CFT correspondence. Let $g$ be an AH Einstein metric on an 
$(n+1)$-manifold $M$ with boundary metric $\gamma$, and such that the corresponding 
geodesic compactification is at least $C^{n}$ in the polyhomogeneous sense.  
The expansions (2.9)-(2.10) easily lead to an expansion for the volume of the region 
$B(t) = \{x\in M: t(x) \leq t\}$. As before, the form of the expansion depends 
on the parity of $n$. For $n$ odd,
\begin{equation} \label{e2.15}
vol B(t) = v_{(n)}t^{-n} + v_{(n-2)}t^{-(n-2)} + \cdots + V + o(1) , 
\end{equation}
while for $n$ even,
\begin{equation} \label{e2.16}
vol B(t) = v_{(n)}t^{-n} + v_{(n-2)}t^{-(n-2)} + \cdots + L\log t + V + o(1) . 
\end{equation}
The constant term $V$ is called the renormalized volume. More importantly for 
our purposes, let ${\mathcal I}_{EH}$ be the Einstein-Hilbert action, (with 
Gibbons-Hawking-York boundary term), given by 
\begin{equation} \label{e2.17}
{\mathcal I}_{EH} = \int_{M}(s - 2\Lambda)dvol + {\tfrac{1}{2}}\int_{\partial M}H dv,
\end{equation}
where $s$ is the scalar curvature, $\Lambda$ is the cosmological constant and $H$ 
is the mean curvature of the boundary. In the normalization (1.2), $\Lambda = -n(n-1)/2$. 
Einstein metrics, with scalar curvature $-n(n+1)$, are critical points of ${\mathcal I}_{EH}$, 
among variations fixing the boundary metric. (This is of course not the case for the volume 
functional). The action ${\mathcal I}_{EH}$ is infinite on $M$, but since $s$ is constant on 
Einstein metrics, (i.e.~on shell), the cut-off action has an expansion of the same form as 
(2.15) or (2.16). Subtracting the divergent contributions as in (2.15)-(2.16) gives the 
renormalized action ${\mathcal I}_{EH}^{ren}$. On the space $E_{AH}$, i.e.~on shell, when $n$ 
is odd one has
$${\mathcal I}_{EH}^{~ren} = -2nV,$$
as functionals on $E_{AH}$; the boundary term in (2.17) renormalizes to 0. Similarly in this 
case, the renormalized action or volume is in fact independent of the defining function $t$, 
so these functionals descend to functionals on the moduli space ${\mathcal E}_{AH}$. Neither of 
these facts hold however for $n$ even. 

  When $dim M = 4$, it is proved in [5] that on $E_{AH}$, ${\mathcal I}_{EH}^{ren}$ and 
$V$ are related to the square of the $L^{2}$ norm of the Weyl curvature ${\mathcal W}$ 
of $(M, g)$ in the following simple way:
\begin{equation} \label{e2.18}
\frac{1}{8\pi^{2}}\int_{M}|W|^{2} = \chi(M) + \frac{1}{8\pi^{2}}{\mathcal I}_{EH}^{~ren} = 
\chi (M) -  \frac{3}{4\pi^{2}}V .
\end{equation}
The relation (2.18) will play an important role in understanding the structure of the 
space $E_{AH}$, in particular in \S 3 and \S 5. 

  It is proved in [5], [18], cf.~also [33], that as a 1-form on $E_{AH}^{m,\alpha}$ the 
differential of the renormalized action is given by
\begin{equation} \label{e2.19}
d_{g}{\mathcal I}_{EH}^{ren}(h) = c_{n}\int_{\partial M} \langle \tau_{(n)}, h_{(0)} \rangle dv,
\end{equation}
where $\tau_{(n)}$ is as in (2.12), the inner product and volume form are with respect to 
$\gamma$ and $c_{n}$ is a constant depending only on $n$. The term $h_{(0)}$ is the variation 
of the boundary metric induced by $h$, i.e.~$h_{(0)} = \Pi_{*}(h)$. 

  For later purposes, we need a slight but important extension of (2.19). Thus, given 
$g \in E_{AH}$, define an enlarged ``tangent'' space $\widetilde T_{g}E_{AH}$ roughly as 
follows: $h \in \widetilde T_{g}E_{AH}$ if $h \in {\mathbb S}^{m,\alpha}(M)$ and $h$ is an 
infinitesimal Einstein deformation to order at least $n$ at $t = 0$, in the sense that 
the expansion (2.9) or (2.10) holds to order $n$ in $t$, to $1^{\rm st}$ order in the 
variation $h$. More precisely, consider the curve $g_{u} = g + uh$ and assume without 
loss of generality that the geodesic defining function associated to $g_{u}$ is the 
fixed function $t$; (this can always be achieved by modifying $g_{u}$ by a curve of 
diffeomorphisms in ${\mathcal D}_{1}$ if necessary). Then $\bar g_{u}$ has the 
expansion
\begin{equation} \label{e2.20}
\bar g_{u} = dt^{2} + (g_{(0)}+uh_{(0)}) + t(g_{(1)}+uh_{(1)}) + t^{2}(g_{(2)}+uh_{(2)}) + 
\cdots + t^{n}\log t({\mathcal H}_{g+uh})
\end{equation}
$$+ t^{n}(g_{(n)}+uh_{(n)}) + O(t^{n+\alpha}) .$$
Then $h \in \widetilde T_{g}E_{AH}$ means that the coefficients $h_{(i)} = 
\frac{d}{du}(g_{(i)}+uh_{(i)})|_{u=0}$ are the linearizations of the coefficients $g_{(i)}$ 
and ${\mathcal H}$ in (2.9)-(2.10) in the direction $h_{(0)}$, for $i \leq n$. Similarly, the 
linearization of the trace condition in (2.11) should hold in the direction $h_{(0)}$. Here 
we recall that the coefficients $g_{(i)}$, ${\mathcal H}$ and $tr g_{(n)}$ are explicitly 
determined from the Einstein equations by the boundary metric $\gamma = g_{(0)}$. 
\begin{lemma} \label{l 2.2}
For any boundary variation $h_{(0)} \in {\mathbb S}^{2}(\partial M)$, there exists 
$h \in \widetilde T_{g}E_{AH}$ such that $h$ induces $h_{(0)}$ at $\partial M$, 
i.e.~for $\Pi$ as in \eqref{e2.1}, $\Pi_{*}:\widetilde T_{g}E_{AH} \rightarrow  
{\mathbb S}^{2}(\partial M)$ is surjective. 

  Moreover, at any $g \in E_{AH}$, the equation \eqref{e2.19} holds for any $h \in 
\widetilde T_{g}E_{AH}$. 
\end{lemma}

\noindent
{\bf Proof:} The first statement follows immediately from the fact that the equations 
for the terms $h_{(i)}$ in the expansion of $h$ can easily be solved algebraically. 
For example, the $g_{(2)}$ term is given explicitly by
$$g_{(2)} = \frac{1}{n-2}(Ric_{\gamma} - \frac{s_{\gamma}}{2(n-1)}\gamma).$$
The linearization of the right-hand side of this equation in the direction $h_{(0)}$ 
then defines the term $h_{(2)}$. The same reasoning applies to all the other coefficients, 
including the $\log$ coefficient ${\mathcal H}$, as well as the trace constraint (2.11). 
Of course $h_{(i)} = 0$ if $i$ is odd, $i < n$. 

   The second statement follows directly from the proof of (2.19) in [5], to which 
we refer for some further details. First, it is easy to see that ${\mathcal I}_{EH}^{ren}$ 
is well-defined, and so finite, for any metric $g_{u} = t^{-2}\bar g_{u}$ as in (2.20) 
which is Einstein to $n^{\rm th}$ order at $\partial M$; see also [5, Remark 1.2]. Further, 
as noted above, Einstein metrics are critical points of ${\mathcal I}_{EH}^{ren}$, and so its 
variation in any direction depends only on the variation of the metric at the boundary. 
Thus, $d{\mathcal I}_{EH}^{ren}$ depends only on the boundary terms $g_{(j)}$ in the 
Fefferman-Graham expansion (2.9)-(2.10). Since it is exactly the form of these 
terms which leads to the formula (2.19), and since this form is preserved for the metrics 
$g_{u}$, (to 1st order in $u$), it is essentially clear that (2.19) holds for 
$h \in \widetilde T_{g}{\mathcal E}_{AH}$. Alternately, this may verified directly 
by an examination of the proof of (2.19) in [5, Lemma 2.1, Theorem 2.2], with 
${\mathcal I}_{EH}$ in place of volume. 

{\endproof}

  The divergence constraint in (2.11) is not used in the computation of 
${\mathcal I}_{EH}^{ren}$ or its first variation. (It does arise however in the 
second variation). Let ${\mathcal D}_{g}E_{AH} \subset {\mathbb S}^{m,\alpha}(M)$ 
be the subspace consisting of forms satisfying the linearized divergence constraint 
(2.11) at $\partial M$, so that 
\begin{equation} \label{e2.21}
\delta'(\tau_{(n)}) + \delta \tau_{(n)}' = 0,
\end{equation}
where $\delta' = \frac{d}{du}\delta_{(g_{(0)}+uh_{(0)})}$ and $\tau_{(n)}' = 
g_{(n)}' + r_{(n)}' = h_{(n)} + s_{(n)}$, $s_{(n)}' = 
\frac{d}{du}(r_{(n)})_{(g_{(0)}+uh_{(0)})}$. Of course $s_{(n)} = 0$ for $n$ odd. 

  The space ${\mathcal F}_{g}E_{AH} = {\mathcal D}_{g}E_{AH} \cap \widetilde T_{g}E_{AH}$ 
represents the space of formal solutions of the linearized Einstein equations near 
$\partial M$, in that any $h \in {\mathcal F}_{g}E_{AH}$ defines uniquely a formal 
series solution as in (2.9)-(2.10) of the linearized Einstein equations. If $h_{(0)}$ 
and $h_{(n)}$ are real-analytic on $\partial M$, the result of Kichenassamy [28] 
mentioned above implies that the series converges to an actual linearized Einstein 
deformation $h$ defined near $\partial M$. Note that while the divergence and 
trace of $h_{(n)}$ are determined by $h_{(0)}$, the transverse-traceless 
part of $h_{(n)}$ may be freely chosen. 

\begin{proposition} \label{p2.3}
For $g \in E_{AH}$, the map $\Pi_{*}: {\mathcal D}_{g}E_{AH} \rightarrow 
{\mathbb S}^{2}(\partial M)$, $\Pi_{*}(h) = h_{(0)}$, is surjective. 
\end{proposition}

\noindent
{\bf Proof:} Given Lemma 2.2, one needs to show that (2.21) is solvable, for any boundary 
variation $h_{(0)}$, i.e.~$\delta'(\tau_{(n)}) \in Im \delta$, for all variations $h_{(0)}$ 
of $\gamma$. The space of 1-forms on $\partial M$ has the splitting
$$\Omega^{1}(\partial M) = Im \delta \oplus Ker \delta^{*},$$
so that if $Ker \delta^{*} = 0$, i.e. $(\partial M, \gamma)$ has no Killing fields, then 
the result is clear. When $(\partial M, \gamma)$ does have Killing fields, this result is 
far from clear; using Theorem 2.5 below, this result is proved in [8], to which we refer 
for details. 

{\endproof}

  We will see later in \S 4, (in and following \eqref{e4.43}), that $\Pi_{*}$ is surjective 
on the full formal space ${\mathcal F}_{g}E_{AH}$. It is worth noting that the global boundary 
map $\Pi_{*}: T_{g} \rightarrow {\mathbb S}^{2}(\partial M)$ is not surjective in general, 
cf.~\S 7. Thus, the formal solutions of the linearized Einstein equations in 
${\mathcal F}_{g}E_{AH}$, even if they give rise to actual linearized solutions, 
do not extend in general to smooth solutions on the compact manifold $M$. 

  On the other hand, Proposition 2.3 is false for AH Einstein metrics defined only 
in a neighborhood or thickening of $\partial M$. The proof of Proposition 2.3 in 
[8] is global; it requires that $(M, g)$ is conformally {\it compact}, as does 
Theorem 2.5 below. 

\medskip

The following result will play an important role in \S 4, cf.~the proof of Proposition 4.6. 
 
\begin{proposition} \label{p2.4}
For $g \in E_{AH}$ and $h \in {\mathcal F}_{g}E_{AH}$, let $\sigma_{(n)} = \tau_{(n)}'$, 
for $\tau_{(n)}'$ as in \eqref{e2.21}. Similarly, let $a_{(n)}'$ be the variation of the 
conformal anomaly $a_{(n)}$ in the direction $h$. Then for any boundary metric 
$(\partial M, \gamma)$, $\gamma \in [\gamma] = \Pi[g]$, and smooth vector field $X$ on 
$\partial M$, one has 
\begin{equation}\label{e2.22}
\int_{\partial M}\langle {\mathcal L}_{X}\tau_{(n)} + [(1 - {\tfrac{2}{n}})div X]\tau_{(n)}, 
h_{(0)} \rangle dV = 
\int_{\partial M}\langle \sigma_{(n)} + {\tfrac{1}{2}}tr h_{(0)} \tau_{(n)}, 
\hat{\mathcal L}_{X}\gamma \rangle dV + 
\end{equation}
$${\tfrac{1}{n}}\int_{\partial M}div X(a_{(n)}tr h_{(0)} + 2a_{(n)}')dV,$$
where $\hat{\mathcal L}_{X}\gamma$ is the conformal Killing operator on $\partial M$, 
$\hat{\mathcal L}_{X}\gamma = {\mathcal L}_{X}\gamma - \frac{2div X}{n}\gamma$. The 
formula \eqref{e2.22} holds for any variation $h_{(0)}$ of $\gamma$. 
\end{proposition}

\noindent
{\bf Proof:} 
This is a simple consequence of a result proved in [8]. Namely, from [8, Prop.5.4], one has 
\begin{equation}\label{e2.23}
\int_{\partial M}\langle {\mathcal L}_{X}\tau_{(n)}, h_{(0)} \rangle = 
-2\int_{\partial M}\langle \delta'(\tau_{(n)}), X \rangle 
+ \int_{\partial M}\delta X \langle \tau_{(n)}, h_{(0)} \rangle + \langle \tau_{(n)}, 
\delta^{*}X \rangle tr h_{(0)}dV ,
\end{equation}
where $\delta' = \frac{d}{du}\delta_{\gamma+uh_{(0)}}$. Using the definition of 
$\hat{\mathcal L}_{X}\gamma$, this gives 
$$\int_{\partial M}\langle {\mathcal L}_{X}\tau_{(n)} + div X \tau_{(n)}, h_{(0)} \rangle = 
-2\int_{\partial M}\langle \delta'(\tau_{(n)}), X \rangle$$
$$+ {\tfrac{1}{2}}\int_{\partial M} \langle \tau_{(n)}, \hat{\mathcal L}_{X}\gamma \rangle 
tr h_{(0)}dV + {\tfrac{1}{n}}\int_{\partial M} a_{(n)}div X tr h_{(0)}dV.$$
Next, by (2.21) and Proposition 2.3, $\delta'\tau_{(n)} = -\delta \sigma_{(n)}$, 
for any variation $h_{(0)}$ of the boundary metric $\gamma$, so that 
$$-2\int_{\partial M}\langle \delta^{'}\tau_{(n)}, X \rangle = 
2\int_{\partial M}\langle \delta \sigma_{(n)}, X \rangle = 
2\int_{\partial M}\langle \sigma_{(n)}, \delta^{*}X \rangle  = 
\int_{\partial M}\langle \sigma_{(n)}, {\mathcal L}_{X}\gamma \rangle .$$   
Write ${\mathcal L}_{X}\gamma = \hat {\mathcal L}_{X}\gamma + \frac{2div X}{n}\gamma$. Then 
$\langle \sigma_{(n)}, \gamma \rangle = tr \sigma_{(n)} = - tr'(\tau_{(n)}) + a_{(n)}' = 
\langle \tau_{(n)}, h_{(0)} \rangle + a_{(n)}'$, where the middle equality follows 
from the linearization of (2.11) with $h \in \widetilde T_{g}E_{AH}$. Combining these 
computations gives (2.22). 

{\endproof}

  Observe that the first term on the right-hand side of (2.22) vanishes when $X$ is a 
conformal Killing field of $(\partial M, [\gamma])$. Using (2.14), it is easy to see 
that both sides of (2.22) are conformally invariant when $n$ is odd; on the other 
hand, both sides of (2.22) depend on the representative $\gamma \in [\gamma]$ 
when $n$ is even. Note also that (2.22) is invariant under the addition of 
transverse-traceless terms to $\sigma_{(n)}$. 

  The following result, also from [8] and used to prove Proposition 2.3, will be 
frequently used in \S 7. 

\begin{theorem} \label{t 2.5.}
  Let $g$ be an AH Einstein metric on an $(n+1)$-manifold $M$, with 
$C^{\infty}$ boundary metric $\gamma$, (in some compactification). 
Then any connected Lie group of isometries of $(\partial M, \gamma)$ extends 
to an action by isometries on $(M, g)$. 
\end{theorem}

\begin{remark} \label{r 2.6.}
  {\rm Finally, for the work in \S 5, we mention briefly the class of toral 
AdS black hole metrics on $M = {\mathbb R}^{2}\times T^{n-1}$. These are AH Einstein 
metrics given explicitly by
\begin{equation} \label{e2.24}
g_{m} = V^{-1}dr^{2} + Vd\theta^{2} + r^{2}g_{T^{n-1}},
\end{equation}
where $V(r) = r^{2} - \frac{2m}{r^{n-2}}$, $\theta \in [0, 4\pi/nr_{+}]$, $r_{+} = 
(2m)^{1/n} > 0$ and $g_{T^{n-1}}$ is any flat metric on the torus $T^{n-1}$, 
cf.~[6, Prop.4.4] and references therein. The boundary metric of $g_{m}$ is a flat 
(product) metric on the $n$-torus $T^{n}$.

  Twisted versions of these metrics, obtained by taking suitable covering spaces 
and then passing to (different) discrete quotients, give rise to infinitely many 
isometrically distinct AH Einstein metrics $g_{i}$ on ${\mathbb R}^{2}\times T^{n-1}$ 
with a fixed flat boundary metric $T^{n-1}$. These metrics are also the examples 
discussed following the statement of Theorem A. It is worth noting that the metrics 
$g_{m}$ are all locally isometric. } 
\end{remark}

\section{Compactness I: Interior Behavior.}
\setcounter{equation}{0}

 The purpose of the next two sections is to lay the groundwork for the proof of 
Theorem A; this section deals with the interior behavior, while \S 4 is 
concerned with the behavior near the boundary. The property that the 
boundary map $\Pi$ is proper is a compactness issue. Thus, given a 
sequence of boundary metrics $\gamma_{i}$ converging to a limit metric 
$\gamma ,$ one needs to prove that a sequence of AH Einstein metrics 
$(M, g_{i})$ with boundary metrics $\gamma_{i}$ has a subsequence 
converging, modulo diffeomorphisms, to a limit AH Einstein metric $g$ 
on $M$, with boundary metric $\gamma$.

 In \S 3.1, we summarize background material on convergence and 
degeneration of sequences of metrics in general, as well as sequences 
of Einstein metrics. This section may be glanced over and then referred 
to as necessary. The following section \S 3.2 then applies these 
results to the interior behavior of AH Einstein metrics.

\medskip

{\bf \S 3.1.}
 In this section, we discuss $L^{p}$ Cheeger-Gromov theory as well as 
the convergence and degeneration results of Einstein metrics on 
4-manifolds from [1]-[3]. 

 We begin with the $L^{p}$ Cheeger-Gromov theory. The $L^{\infty}$ 
Cheeger-Gromov theory [15], [22], describes the (moduli) space of 
metrics on a manifold, (or sequence of manifolds), with uniformly 
bounded curvature in $L^{\infty}$, i.e.
\begin{equation} \label{e3.1}
|R_{g}|(x) \leq  \Lambda  <  \infty , 
\end{equation}
in that it describes the convergence or possible degenerations of 
sequences of metrics satisfying the bound (3.1). The space $L^{\infty}$ 
is not a good space on which to carry out analysis, and so we replace 
(3.1) by a corresponding $L^{p}$ bound, i.e.
\begin{equation} \label{e3.2}
\int_{M}|R_{g}|^{p}dV \leq  \Lambda  <  \infty . 
\end{equation}
The curvature involves 2 derivatives of the metric, and so (3.2) is 
analogous to an $L^{2,p}$ bound on the metric. The critical exponent 
$p$ with respect to Sobolev embedding $L^{2,p} \subset  C^{o}$ is $(n+1)/2$, 
where $n+1 = dim M$ and hence we will always assume that
\begin{equation} \label{e3.3}
p >  (n+1)/2. 
\end{equation}

 In order to obtain local results, we need the following definitions of 
local invariants of Riemannian metrics, cf. [4].

\begin{definition} \label{d 3.1.}
 {\rm If $(M, g)$ is a Riemannian $(n+1)$-manifold, the $L^{p}$ {\it 
curvature radius}  $\rho (x) \equiv \rho^{p}(x)$ at $x$ is the 
radius of the largest geodesic ball $B_{x}(\rho (x))$ such that, for 
all $B_{y}(s) \subset  B_{x}(\rho (x))$, with $s \leq dist(y, \partial 
B_{x}(\rho (x)))$, one has
\begin{equation} \label{e3.4}
\frac{s^{2p}}{vol B_{y}(s)}\int_{B_{y}(s)}|R|^{p}dV \leq  c_{0}, 
\end{equation}
where $c_{0}$ is a fixed sufficiently small constant. Although $c_{0}$ 
is an essentially free parameter, we will fix $c_{0} = 10^{-2}$ 
throughout the paper. The left-side of (3.4) is a scale-invariant local 
average of the curvature in $L^{p}.$}

 {\rm The {\it  volume radius}  $\nu (x)$ of $(M, g)$ at $x$ is given by
\begin{equation} \label{e3.5}
\nu (x) = \sup\{r: \frac{vol B_{y}(s)}{\omega_{n+1}s^{n+1}} \geq  \mu_{0}, 
\forall B_{y}(s) \subset  B_{x}(r)\}, 
\end{equation}
where $\omega_{n+1}$ is the volume of the Euclidean unit $(n+1)$-ball and 
again $\mu_{0} > 0$ is a free small parameter, which will be fixed 
in any given discussion, e.g. $\mu_{0} = 10^{-2}$.}

 {\rm The $L^{k,p}$ {\it harmonic radius} $r_{h}^{k,p}(x)$ is the largest 
radius such that on the geodesic ball $B_{x}(r)$, $r = r_{h}^{k,p}(x)$, there 
is a harmonic coordinate chart in which the metric components satisfy 
\begin{equation} \label{e3.6}
C^{-1}\delta_{ij} \leq g_{ij} \leq C\delta_{ij}
\end{equation}
$$r^{kp-(n+1)}\int_{B_{x}(r)}|\partial^{k}g_{ij}|^{p}dV_{g} \leq C,$$
where $C$ is a fixed constant; again $C$ may be arbitrary, but will be 
fixed to $C = 2$, (for example). }
\end{definition}

 Observe that $\rho(x)$, $\nu (x)$ and $r_{h}(x)$ scale as distances, 
i.e.~if $g'  = \lambda^{2}\cdot g$, for some constant $\lambda$, then 
$\rho'(x) = \lambda\cdot \rho (x)$, $\nu' (x) = \lambda\cdot \nu (x)$ 
and $r_{h}'(x) = \lambda\cdot r_{h}(x)$. By 
definition, $\rho$, $\nu$ and $r_{h}$ are Lipschitz functions with Lipschitz 
constant 1; in fact for $y\in B_{x}(\rho (x))$, it is immediate from 
the definition that
\begin{equation}\label{e3.7}
\rho (y) \geq  dist(y, \partial B_{x}(\rho (x))),
\end{equation}
and similarly for $\nu$, $r_{h}$.

  One may define $L^{k,p}$, (or $C^{k,\beta}$) curvature radii $\rho^{k,p}$ 
in a manner completely analogous to (3.4). For the $L^{k,p}$ radius, the 
curvature $|R|$ term in (3.4) is replaced by $\sum_{j\leq k}|\nabla^{j}R|$ 
and the power of $s$ is chosen to make the resulting expression scale invariant. 
Clearly 
$$\rho^{k,p} \geq \rho^{k',p'},$$
whenever $k + p \geq k' + p'$. It is proved in [3] that the curvature radius 
and harmonic radius are essentially equivalent, given a lower bound on the 
volume radius. Thus
\begin{equation}\label{e3.8}
r_{h}^{2,p}(x) \geq r_{0}\rho^{p}(x), 
\end{equation}
where $r_{0}$ depends only on a lower bound $\nu_{0}$ for $\nu (x)$. The 
same statement holds for the radii $r_{h}^{k+2,p}$ and $\rho^{k,p}$. 

 A sequence of Riemannian metrics $(\Omega_{i}, g_{i})$ is said to 
converge in the $L^{k,p}$ topology to a limit $L^{k,p}$ metric $g$ on 
$\Omega $ if there is an atlas ${\mathcal A} $ for $\Omega $ and 
diffeomorphisms $F_{i}: \Omega  \rightarrow  \Omega_{i}$ such that 
$F_{i}^{*}(g_{i})$ converges to $g$ in the $L^{k,p}$ topology in local 
coordinates with respect to the atlas ${\mathcal A} .$ Thus, the local components 
$(F_{i}^{*}(g_{i}))_{\alpha\beta} \rightarrow  g_{\alpha\beta}$ in the 
usual $L^{k,p}$ Sobolev topology on functions on ${\mathbb R}^{n+1}.$ Similar 
definitions hold for $C^{m,\alpha}$ convergence, and convergence in the 
weak $L^{k,p}$ topology. Any bounded sequence in $L^{k,p}$ has a weakly 
convergent subsequence, and similarly any bounded sequence in 
$C^{m,\alpha}$ has a convergent subsequence in the $C^{m,\alpha'}$ 
topology, for any $\alpha'  <  \alpha .$ By Sobolev embedding,
$$L^{k,p} \subset  C^{m,\alpha}, $$
for $m+\alpha  <  k-\frac{n+1}{p}.$ In particular, in dimension 4, for 
$p\in (2,4), L^{2,p} \subset  C^{\alpha}, \alpha  = 2-\frac{4}{p},$ and 
for $p > $ 4, $L^{2,p} \subset  C^{1,\alpha},$ for $\alpha  = 
1-\frac{4}{p}.$

 One then has the following result on the convergence and degeneration 
of metrics with bounds on $\rho$, cf.~[4].

\begin{theorem} \label{t 3.2.}
  Let $(\Omega_{i}, g_{i}, x_{i})$ be a pointed sequence of connected 
Riemannian $(n+1)$-manifolds and suppose there are constants $\rho_{o} > $ 
0, $d_{o} > $ 0 and $D <  \infty$ such that, for a fixed $p >  (n+1)/2,$
\begin{equation} \label{e3.9}
\rho^{p}(y) \geq  \rho_{o}, \ diam \Omega_{i} \leq  D, \ dist(x_{i}, 
\partial\Omega_{i}) \geq  d_{o}. 
\end{equation}
Then for any 0 $<  \epsilon  <  d_{o},$ there are domains $U_{i} 
\subset  \Omega_{i},$ with $\epsilon /2 \leq  dist(\partial U_{i}, 
\partial\Omega_{i}) \leq  \epsilon ,$ for which one of the following 
alternatives holds.

(I). (Convergence) If there is constant $\nu_{o} > $ 0 such that,
$$\nu_{i}(x_{i}) \geq  \nu_{o} >  0, $$
then a subsequence of $\{(U_{i}, g_{i}, x_{i})\}$ converges, in the 
weak $L^{2,p}$ topology, to a limit $L^{2,p}$ Riemannian manifold (U, 
g, x), $x = \lim x_{i}$. In particular, $U_{i}$ is diffeomorphic to 
$U$, for $i$ sufficiently large.

(II). (Collapse) If instead,
$$\nu (x_{i}) \rightarrow  0, $$
then $U_{i}$ has an F-structure in the sense of Cheeger-Gromov, cf.~[15]. 
The metrics $g_{i}$ are collapsing everywhere in $U_{i}$, 
i.e.~$\nu_{i}(y_{i}) \rightarrow $ 0 for all $y_{i} \in  U_{i}$, and so 
in particular the injectivity radius $inj_{g_{i}}(y_{i}) \rightarrow 0$. 
Any limit of $\{g_{i}\}$ in the Gromov-Hausdorff topology [22] is a 
lower dimensional length space.
\end{theorem}

\begin{remark} \label{r 3.3.}{\bf (i).}
  {\rm The local hypothesis on $\rho^{p}$ in (3.9) can be replaced by 
the global hypothesis
\begin{equation} \label{e3.10}
\int_{\Omega_{i}}|R|^{p}dV \leq  \Lambda , 
\end{equation}
in case the volume radius satisfies $\nu_{i}(y_{i}) \geq  \nu_{o}$, for 
some $\nu_{o} > 0$ and all $y_{i}\in\Omega_{i}$; in fact under this 
condition (3.9) and (3.10) are then equivalent, with $\Lambda  = \Lambda 
(\rho_{o}, D, \nu_{o})$.

{\bf (ii).}
 It is an easy consequence of the definitions that the $L^{p}$ 
curvature radius $\rho^{p}$ is continuous under convergence in the 
(strong) $L^{2,p}$ topology, i.e. if $g_{i} \rightarrow g$ in the 
$L^{2,p}$ topology, then
\begin{equation} \label{e3.11}
\rho_{i}(x_{i}) \rightarrow  \rho (x), 
\end{equation}
whenever $x_{i} \rightarrow  x$, cf. [4] and references therein. 
However, (3.11) does {\it not}  hold if the convergence is only in the 
weak $L^{2,p}$ topology. 

 If in addition to a bound on $\rho^{p}$ one has $L^{p}$ bounds the 
covariant derivatives of the Ricci curvature up to order $k$, then in 
Case (I) one obtains convergence to the limit in the $L^{k+2,p}$ 
topology. Analogous statements hold for convergence in the 
$C^{m,\alpha'}$ topology.}
\end{remark}

 Next we discuss the convergence and degeneration of Einstein metrics. 
If $M$ is a closed 4-manifold and $g$ an Einstein metric on $M$, then 
the Chern-Gauss-Bonnet theorem gives 
\begin{equation} \label{e3.12}
\frac{1}{8\pi^{2}}\int_{M}|R|^{2}dV = \chi (M), 
\end{equation}
where $\chi (M)$ is the Euler characteristic of $M$. This gives apriori 
control on the $L^{2}$ norm of the curvature of Einstein metrics on 
$M$. However the $L^{2}$ norm is critical in dimension 4 with respect to Sobolev 
embedding, cf. (3.3) and so one may not expect Theorem 3.2 to hold for 
sequences of Einstein metrics on $M$. In fact, there is a further 
possible behavior of such sequences.

\begin{definition} \label{d 3.4.}
 {\rm An {\it Einstein orbifold} $(X, g)$ is a 4-dimensional orbifold with a 
finite number of cone singularities $\{q_{j}\}, j = 1, ..., k$. On 
$X_{o} = X \setminus \cup\{q_{j}\}$, $g$ is a smooth Einstein metric 
while each singular point $q\in\{q_{j}\}$ has a neighborhood $U$ such 
that $U \setminus q$ is diffeomorphic to $C(S^{3}/\Gamma ) \setminus 
\{0\}$, where $\Gamma $ is a finite subgroup of $O(4)$ and $C$ denotes 
the cone with vertex $\{0\}$. Further $\Gamma  \neq  \{e\}$ and when 
lifted to the universal cover $B^{4} \setminus \{0\}$ of 
$C(S^{3}/\Gamma ) \setminus \{0\}$, the metric $g$ extends smoothly 
over $\{0\}$ to a smooth Einstein metric on the 4-ball $B^{4}.$}
\end{definition}

 There are numerous examples of sequences of smooth Einstein metrics $g_{i}$ 
on a compact manifold $M$ which converge, in the Gromov-Hausdorff topology, 
to an Einstein orbifold limit $(X, g)$. Such orbifold metrics are not 
smooth metrics on the manifold $M$, but may be viewed as singular 
metrics on $M$, in that $M$ is a resolution of $X$.

 We then have the following result describing the convergence and 
degeneration of Einstein metrics on 4-manifolds, cf. [3].

\begin{theorem} \label{t 3.5.}
  Let $(\Omega_{i}, g_{i}, x_{i})$ be a pointed sequence of connected 
Einstein $4$-manifolds satisfying
\begin{equation} \label{e3.13}
diam \Omega_{i} \leq  D, \ \ dist(x_{i},\partial\Omega_{i}) \geq  
d_{o}, 
\end{equation}
and
\begin{equation} \label{e3.14}
\int_{\Omega_{i}}|R_{g_{i}}|^{2}dV_{g_{i}} \leq  \Lambda_{o}, 
\end{equation}
for some constants $d_{o} > $ 0, D, $\Lambda_{o} < \infty .$ Then for 
any 0 $<  \epsilon  <  d_{o},$ there are domains $U_{i} \subset  
\Omega_{i},$ with $\epsilon /2 \leq  dist(\partial U_{i}, 
\partial\Omega_{i}) \leq  \epsilon ,$ for which exactly one of the 
following alternatives holds.

(I). (Convergence). A subsequence of $(U_{i}, g_{i}, x_{i})$ converges 
in the $C^{\infty}$ topology to a limit smooth Einstein metric (U, g, 
x), $x = \lim x_{i}.$

(II). (Orbifolds). A subsequence of $(U_{i}, g_{i}, x_{i})$ converges 
to an Einstein orbifold metric $(X, g, x)$ in the Gromov-Hausdorff 
topology. Away from the singular variety, the convergence $g_{i} 
\rightarrow  g$ is $C^{\infty},$ while the curvature of $g_{i}$ blows 
up in $L^{\infty}$ at the singular variety. 

(III). (Collapse). A subsequence of $(U_{i}, g_{i}, x_{i})$ collapses, 
in that $\nu (y_{i}) \rightarrow $ 0 for all $y_{i}\in U_{i}.$ The 
collapse is with ``generalized'' bounded curvature, in that 
$inj_{g_{i}}^{2}(x)|R|_{g_{i}}(x) \rightarrow 0$, metrically away from a 
finite number of singular points; however such singularities might be 
more complicated than orbifold cone singularities.
\end{theorem}

 The cases (I) and (II) occur if and only if
\begin{equation} \label{e3.15}
\nu_{i}(x_{i}) \geq  \nu_{o}, 
\end{equation}
for some $\nu_{o} > $ 0, while (III) occurs if (3.15) fails. One 
obtains $C^{\infty}$ convergence in (I), and (II) away from the 
singularities, since Einstein metrics satisfy an elliptic system of 
PDE, (in harmonic coordinates).

 When $\Omega_{i}$ is a closed 4-manifold $M$, the bound (3.14) follows 
immediately from (3.12). When $(\Omega_{i}, g_{i})$ are complete AHE 
metrics on a fixed 4-manifold $M$, then of course (3.14) does not hold 
with $\Omega_{i} = M$. However, in the normalization (1.2), one has 
$|W|^{2} = |R|^{2}- 6$ where $W$ is the Weyl tensor of $(M, g)$ and so 
by (2.18)
\begin{equation} \label{e3.16}
\frac{1}{8\pi^{2}}\int_{M}(|R_{g}|^{2} -  6)dV_{g} = \chi (M) -  
\frac{3}{4\pi^{2}}V. 
\end{equation}
Thus, a {\it lower} bound on the renormalized volume $V$ (or upper bound 
on the renormalized action) and an upper bound on $\chi (M)$ give a global 
bound on the $L^{2}$ norm of $W$. In particular, under these bounds (3.14) 
holds for $\Omega_{i}$ a geodesic ball $B_{x_{i}}(R)$ of any fixed radius $R$ 
about the base point $x_{i}\in (M, g_{i})$. 

 We point out that Theorem 3.5 is special to dimension 4. A similar 
result holds in higher dimensions only if one has a uniform bound on 
the $L^{n/2}$ norm of curvature in place of (3.14); however there is no 
analogue of (3.16) for (general) Einstein metrics in higher dimensions.

\begin{remark} \label{r 3.6.}
  {\rm Theorems 3.2 and 3.5 are local results. However, we will often 
apply them globally, to sequences of complete manifolds, and with 
complete limits. This is done by a standard procedure as follows, in 
the situation of Theorem 3.2 for example. Suppose $\{g_{i}\}$ is a 
sequence of complete metrics on a manifold $M$, or on a sequence of 
manifolds $M_{i}$, with base points $x_{i}$, and satisfying 
$\rho (y_{i}) \geq  \rho_{o},$ for all $y_{i}\in (M, g_{i})$. One 
may then apply Theorem 3.2 to the domains $(B_{x_{i}}(R), g_{i}, x_{i})$, 
where $B_{x_{i}}(R)$ is the geodesic $R$-ball about $x_{i}$ to obtain a 
limit manifold $(U(R), g_{\infty}, x_{\infty}),$ in the non-collapse case. 
Now take a divergent sequence $R_{j} \rightarrow  \infty$ and carry out 
this process for each $j$. There is then a diagonal subsequence 
$B_{x_{i}}(R_{i})$ of $B_{x_{i}}(R_{j})$ which converges to a complete 
limit $(N, g_{\infty}, x_{\infty}).$ Similar arguments apply in the case 
of collapse and the cases of Theorem 3.5.

 The convergence in these situations is also convergence in the pointed 
Gromov-Hausdorff topology, [22].}
\end{remark}

{\bf \S 3.2.}
 In this section, we study the behavior of sequences of AH Einstein metrics 
in the interior, i.e. away from $\partial M.$ We first discuss the 
hypotheses, and then state and prove the main result, Theorem 3.7.

 Essentially for the rest of the paper, we will assume the topological 
condition on $M = M^{4}$ that
\begin{equation} \label{e3.17}
H_{2}(\partial M, {\mathbb F}) \rightarrow  H_{2}(\bar M, {\mathbb F}) 
\rightarrow  0, 
\end{equation}
where the map is induced by inclusion $\iota: \partial M \rightarrow  
\bar M$ and ${\mathbb F}$ is any field. As will be seen, this serves to 
rule out any orbifold degenerations of AH Einstein sequences. While many 
of the results of this paper carry over to the orbifold setting, we prefer 
for simplicity to exclude this possibility here. 

 Next, let $g_{i}\in E_{AH}$ be a sequence of AH Einstein metrics 
on $M$ with $C^{m,\alpha}$ boundary metrics $\gamma_{i}$, $m \geq  3$. 
As noted in \S 2, the geodesic compactification $\bar g_{i}$ of 
$g_{i}$ with boundary metric $\gamma_{i}$ is at least a $C^{m-1,\alpha}$ 
compactification. We assume the boundary behavior of $\{g_{i}\}$ is 
controlled, in that 
\begin{equation} \label{e3.18}
\gamma_{i} \rightarrow  \gamma , 
\end{equation}
in the $C^{m,\alpha'}$ topology on $\partial M,$ for some $\alpha'  
\leq  \alpha .$ Next assume that the inradius of $(M, \bar g_{i})$ 
has a uniform lower bound, i.e.
\begin{equation} \label{e3.19}
In_{\bar g_{i}}(\partial M) =  dist_{\bar g_{i}}(\bar C_{i}, \partial 
M) \geq  \bar \tau , 
\end{equation}
for some constant $0 < \bar \tau \leq 1$, where $\bar C_{i}$ is the cutlocus 
of $\partial M$ in $(\bar M, \bar g_{i})$, and also assume an upper diameter 
bound
\begin{equation} \label{e3.20}
diam_{\bar g_{i}}S(t_{1}) \leq  T, 
\end{equation}
where $t_{1} = \bar \tau /2$ and $S(t_{1}) = \{x\in (M, \bar g_{i}): 
t_{i}(x) = t_{1}\}$ is the $t_{1}$-level set of the geodesic defining 
function $t_{i}$ for $(g_{i}, \gamma_{i}).$

 Finally, we assume that the Weyl curvature of $g_{i}$ is uniformly 
bounded in $L^{2},$ i.e.
\begin{equation} \label{e3.21}
\int_{M}|W_{g_{i}}|^{2}dV_{g_{i}} \leq  \Lambda_{0} <  \infty . 
\end{equation}

 All of these assumptions will be removed in \S 4 and at the beginning of 
\S 5. The main result of this subsection is the following:

\begin{theorem} \label{t 3.7.}
  Let $\{g_{i}\}$ be a sequence of metrics in $E_{AH}$ on $M$, 
satisfying \eqref{e3.17}-\eqref{e3.21} and let $x_{i}$ be a base point 
of $(M, g_{i})$ satisfying
\begin{equation} \label{e3.22}
d \leq  dist_{\bar g_{i}}(x_{i}, \partial M) \leq  D, 
\end{equation}
for some constants $d > $ 0 and $D <  \infty ,$ where $\bar g_{i}$ is 
the geodesic compactification associated to $\gamma_{i}.$

 Then a subsequence of $(M, g_{i}, x_{i}),$ converges to a complete 
Einstein metric $(N, g, x_{\infty})$, $x_{\infty} = \lim x_{i}.$ The 
convergence is in the $C^{\infty}$ topology, uniformly on compact 
subsets, and the manifold $N$ weakly embeds in $M$,
\begin{equation} \label{e3.23}
N \subset  \subset  M, 
\end{equation}
in the sense that smooth bounded domains in $N$ embed as such in $M$.
\end{theorem}

\noindent
{\bf Proof:}
 By Theorem 3.5 and the discussion following (3.16), together with 
Remark 3.6, a subsequence of $\{(M, g_{i})\}$ based at $x_{i}$ either, 
(i): converges smoothly to a complete limit Einstein manifold $(N, g, 
x_{\infty})$, (ii): converges to an Einstein orbifold, smoothly away from 
the singular variety, or (iii): collapses uniformly on domains of 
arbitrary bounded diameter about $\{x_{i}\}$. 

 The first (and most important) task is to rule out the possibility of 
collapse. To do this, we first prove a useful volume monotonicity 
formula in the Lemma below; this result will also be important in \S 4 
and holds in all dimensions. 

 Let $(M, g)$ be any AH Einstein manifold of dimension $n+1$, with 
$C^{2,\alpha}$ geodesic compactification 
\begin{equation} \label{e3.24}
\bar g = t^{2} g 
\end{equation}
and boundary metric $\gamma$. Let $\bar E$ be the inward normal 
exponential map of $(M, \bar g)$ at $\partial M,$ so that 
$\sigma_{x}(t) = \bar E(x, t)$ is a geodesic in $t$, for each fixed 
$x\in\partial M$. As discussed in \S 2, for $r = -\log t$, the curve 
$\sigma_{x}(r)$ is a geodesic in the Einstein manifold $(M, g)$. Let 
$\bar J(x, t)$ be the Jacobian of $\bar E(x, t)$, so that 
$\bar J(x, t) = d\bar V(x, t)/d\bar V(x, 0)$, where $d\bar V$ is the 
volume form of the 'geodesic sphere' $S(t)$, i.e. the $t$-level set of 
$t$ as following (3.20). Also, let $\tau_{0} = \tau_{0}(x)$ be the 
distance to the cutlocus of $\bar E$ at $x\in\partial M$.

\begin{lemma} \label{l 3.8.}
  In the notation above, the function
\begin{equation} \label{e3.25}
\frac{\bar J(x,t)}{\tau_{0}^{n}(1- (\frac{t}{\tau_{0}})^{2})^{n}}  
\end{equation}
is monotone non-decreasing in $t$, for any fixed $x$.
\end{lemma}

\noindent
{\bf Proof:}
 As above, let $r = \log (\frac{1}{t})$ and, for any fixed $x$, set 
$r_{0} = r_{0}(x) = \log (\frac{1}{\tau_{0}})$. Then the curve 
$\sigma_{x}(r)$ is a geodesic in $(M, g)$  and $S(r)$ is the $r$-rlevel 
set of the distance function $r$. Let $J(x, r)$ be the corresponding 
Jacobian for $S(r)$ along $\sigma_{x}(r)$, so that $J(x,r) = t^{-n}\bar 
J(x,t)$ by (3.24). Since $Ric_{g} = - ng$, the infinitesimal form of 
the Bishop-Gromov volume comparison theorem [22] implies the ratio 
$$\frac{J(x,r)}{\sinh^{n}(r- r_{0})} \downarrow  $$
i.e.~the ratio is monotone non-increasing in $r$, as $r \rightarrow  \infty$. 
Converting this back to $(M, \bar g)$, it follows that, in $t$,
$$\frac{t^{-n}\bar J(x,t)}{t^{n}\sinh^{n}(\log (\frac{\tau_{0}}{t}))} \uparrow , $$
is increasing, since as $r$ increases to $\infty$, $t$ decreases to 0. But 
$\sinh^{n}(\log (\frac{\tau_{0}}{t})) = \frac{1}{2^n}(\frac{\tau_{0}}{t})^{n}(1- 
(\frac{t}{\tau_{0}})^{2})^{n}$, which gives (3.25).
{\endproof}

 Let $q$ be any point in $\partial M$ and consider the metric $s$-ball 
$B_{x}(s) = \{y\in\bar M: dist_{\bar g}(y, q) \leq  s\}$ in $(\bar M, 
\bar g).$ Let $D_{q}(s) = B_{q}(s)\cap\partial M,$ so that $D_{q}(s)$ 
is the metric $s$-ball in $(\partial M, \gamma)$, since $\partial M$ 
is totally geodesic. Observe that there are constants $\mu_{0} >  0$ 
and $\mu_{1} <  \infty$, which depend only on the $C^{0}$ geometry of 
the boundary metric $\gamma$, such that
\begin{equation} \label{e3.26}
\mu_{0} \leq  \frac{vol_{\gamma}(D_{q}(s))}{s^{n}} \leq  \mu_{1}. 
\end{equation}

 Let $D_{q}(s, t) = \{x\in M: x = \bar E(y,t)$, for some $y\in 
D_{q}(s)\}$, so that $D_{q}(s,t) \subset S(t)$. It follows immediately 
from (3.25) by integration over (a domain in) $\partial M$ that, for any 
fixed $s$,
\begin{equation} \label{e3.27}
\frac{vol_{\bar g}D(s,t)}{(s\bar \tau)^{n}(1-(\frac{t}{\bar \tau})^{2})^{n}} 
\uparrow , \ \ {\rm and} \ \ 
\frac{vol_{\bar g}S(t)}{\bar \tau^{n}(1-(\frac{t}{\bar \tau})^{2})^{n}} 
\uparrow , 
\end{equation}
for $\bar \tau$ as in (3.19). Now let $t_{1} = \bar \tau/2$, so that 
(3.27) implies in particular that
\begin{equation} \label{e3.28}
vol_{\bar g}S(t_{1}) \geq  (\frac{3}{4})^{n}vol_{\gamma}\partial M, 
\end{equation}
and hence, with respect to $(M, g)$,
\begin{equation} \label{e3.29}
vol_{g}S(t_{1}) \geq  (\frac{3}{4})^{n}t_{1}^{-n}vol_{\gamma}\partial 
M. 
\end{equation}

 This leads easily to the following lower bound on volumes of balls.

\begin{corollary} \label{c 3.9.}
  Let $(M, g)$ be an AH Einstein metric, satisfying \eqref{e3.19}-\eqref{e3.20}, 
and choose a point $x$ satisfying \eqref{e3.22}, i.e.~$d \leq  dist_{\bar g}(x, \partial M) 
\leq  D$. Then
\begin{equation} \label{e3.30}
vol_{g}B_{x}(1) \geq  \nu_{0} >  0, 
\end{equation}
where $\nu_{0}$ depends only on $d,D,\bar \tau, T$ and the $C^{0}$ 
geometry of the boundary metric $\gamma$.
\end{corollary}

\noindent
{\bf Proof:}
 Again, by the Bishop-Gromov volume comparison theorem on $(M, g)$, one 
has for any $R \geq r$,
\begin{equation}\label{e3.31}
\frac{volB_{x}(r)}{\sinh^{n}(r)} \geq  
\frac{volB_{x}(R)}{\sinh^{n}(R)}.
\end{equation}
Suppose $r(x) = D_{1}$, so that $D_{1}\in 
[-\log D,-\log d]$. Then for any $y\in S(t_{1}),$ the 
triangle inequality and (3.20) imply that 
$$dist_{g}(x,y) \leq  |D_{1}|+ t_{1}^{-1}\cdot  T \equiv  D_{2}. $$
Hence, choose $R$ so that $R = D_{2}+1$, which implies that $S(t_{1}) 
\subset  B_{x}(R-1).$ But $volB_{x}(R) \geq  vol A(R-2, R) \geq  
\frac{1}{2}vol S(t_{1})$, where the second inequality follows from the 
coarea formula, (changing $t_{1}$ slightly if necessary). Combining 
this with the comparison estimate above implies
$$vol B_{x}(1) \geq  c_{2}vol S(t_{1})\sinh^{-n}(D_{2}) \geq  
c_{3}t_{1}^{-n}\sinh^{-n}(D_{2})vol_{\gamma}\partial M,$$ 
where the last inequality follows from (3.29). This gives (3.30).
{\endproof}

 Corollary 3.9 gives a uniform lower bound on the volume radius of each 
$g_{i}$ at $x_{i}$ satisfying (3.22) and so there is no possibility of 
collapse. Theorem 3.5 and Remark 3.6 then imply that $(M, g_{i}, 
x_{i})$ has a subsequence converging in the Gromov-Hausdorff topology 
either to a complete Einstein manifold $(N, g, x_{\infty})$ or to a complete 
Einstein orbifold $(V, g, x_{\infty})$.

\medskip

 Next, we use the hypothesis (3.17) to rule out orbifold limits. Thus, 
suppose the second alternative above holds, so that $(M, g_{i})$ 
converges in the Gromov-Hausdorff topology to a complete Einstein 
orbifold $(X, g)$. With each orbifold singularity $q\in X$ with 
neighborhood of the form $C(S^{3}/\Gamma)$, there is associated a 
(preferred) smooth complete Ricci-flat 4-manifold $(E, g_{0})$, which 
is asymptotically locally Euclidean (ALE), in that $(E, g)$ is 
asymptotic to a flat cone $C(S^{3}/\Gamma_{E})$, $\Gamma_{E}  \neq 
\{e\}$, cf.~[1,\S 5], [3,\S 3]. The complete manifold $(E, g_{o})$ 
is obtained as a limit of blow-ups or rescalings of the metrics $g_{i}$ 
restricted to small balls $B_{y_{i}}(\delta_{i})$ at the {\it maximal} 
curvature scale; thus one does not necessarily have $\Gamma = \Gamma_{E}$, 
since there may be other Ricci-flat ALE orbifolds arising from blow-ups 
at other curvature scales. In any case, the manifold $E$ is embedded in 
the ambient manifold $M$; in fact for any $\delta > 0$ and points $y_{i} 
\rightarrow  q$, $E$ is topologically embedded in $(B_{y_{i}}(\delta), 
g_{i})$. The Einstein metrics $g_{i}$ crush the topology 
of $E$ to a point in that $E \subset  (B_{y_{i}}(\delta), g_{i})$ and 
$B_{y_{i}}(\delta)$ converges to the cone $C(S^{3}/\Gamma)$, in the 
Gromov-Hausdorff topology for any fixed $\delta > 0$ sufficiently 
small.

  Next we point out that any such ALE space $E$ has non-trivial topology. 
This result corrects a minor inaccuracy in [1, Lemma 6.3]. 
\begin{proposition}\label{p3.10}
Let $(E, g_{0})$ be a complete, non-flat, oriented Ricci-flat ALE manifold. 
Then 
\begin{equation} \label{e3.32}
H_{2}(E, {\mathbb F} ) \neq  0, 
\end{equation}
for some field ${\mathbb F}$. Moreover, if $E$ is simply connected, 
then $H_{2}(E, {\mathbb Z})$ is torsion-free, so that $H_{2}(E, {\mathbb R}) 
\neq 0$, while $H_{2}(E, {\mathbb F}) = 0$, for any finite field ${\mathbb F}$. 
\end{proposition}

\noindent
{\bf Proof:} It is well-known that $|\pi_{1}(E)| < \infty$, so that $b_{1}(E) = 0$, 
and also $b_{3}(E) = 0$, (cf.~[1] for example). Hence
$$\chi(E) = 1 + b_{2}(E) \geq 1.$$
If $b_{2}(E) > 0$, then of course \eqref{e3.32} follows, so suppose $b_{2}(E) = 0$, 
so that $\chi (E) = 1$. The Euler characteristic may be computed with respect to 
homology with coefficients in any field ${\mathbb F}$, so that 
$$\chi(E) = 1 - H_{1}(E, {\mathbb F}) + H_{2}(E, {\mathbb F}) = 1,$$
since again $H_{3}(E, {\mathbb F}) = H_{4}(E, {\mathbb F}) = 0$. Now if 
\eqref{e3.32} does not hold, then it follows that $H_{1}(E, {\mathbb F}) = 0$, 
for all ${\mathbb F}$. By the universal coefficient theorem, this implies 
$H_{1}(E, {\mathbb Z}) = H_{2}(E, {\mathbb Z}) = 0$. 

  Thus, $E$ is an integral homology ball, with finite $\pi_{1}$. It 
then follows by the arguments in [1, Lemma 6.3] that $E$ is flat, giving 
a contradiction. This proves \eqref{e3.32}. The proof of the second statement 
is an immediate consequence of the universal coefficient theorem. 
{\endproof}

  Next, observe that there is an injection
\begin{equation} \label{e3.33}
0 \rightarrow  H_{2}(E, {\mathbb F} ) \rightarrow  H_{2}(M, {\mathbb F} ). 
\end{equation}
To see this, the Mayer-Vietoris sequence for (a thickening of) the 
pair $(E, M \setminus E)$ gives
$$H_{2}(S^{3}/\Gamma, {\mathbb F}) \rightarrow  H_{2}(E, {\mathbb F} 
)\oplus H_{2}(M \setminus E, {\mathbb F}) \rightarrow  H_{2}(M, {\mathbb F}).$$
Thus, it suffices to show that $i_{*}: H_{2}(S^{3}/\Gamma, {\mathbb F}) 
\rightarrow  H_{2}(E, {\mathbb F})$ is the zero-map. If $E$ is simply 
connected, this follows immediately from Proposition 3.10, since 
$H_{2}(S^{3}/\Gamma, {\mathbb Z})$ is torsion, while $H_{2}(E, {\mathbb Z})$ 
is torsion-free. In general, let $\widetilde E$ be the universal cover of 
$E$, so that the covering map $\pi: \widetilde E \rightarrow E$ is finite-to-one. 
One has $\partial \widetilde E = S^{3}/\widetilde \Gamma$, for some finite group 
$\widetilde \Gamma$. Now any 2-cycle in $S^{3}/\Gamma$ with coefficients in 
${\mathbb F}$ may be represented by a map, (or more precisely a collection of 
maps), $f: S^{2} \rightarrow S^{3}/\Gamma$. This lifts to a map $\widetilde f: 
S^{2} \rightarrow S^{3}/\widetilde \Gamma = \partial \widetilde E$. As 
noted above, the 2-cycle $\widetilde f$ bounds a 3-chain $F$ in 
$\widetilde E$. Composing with the projection map $\pi$, it follows 
that $f$ bounds a 3-chain in $E$, which proves the claim.

   Now returning to \eqref{e3.17} and orbifold limits, let $\Sigma$ 
be any non-zero 2-cycle in $H_{2}(E, {\mathbb F})$; ($\Sigma$ exists 
by \eqref{e3.32}). By (3.17), there is a 2-cycle $\Sigma_{\infty}$ in 
$H_{2}(\partial M, {\mathbb F})$ homologous to $\Sigma$, so that there 
is a 3-chain $W$ in $M$ such that $\partial W = \Sigma - \Sigma_{\infty}$. 
But $Z = W \cap (S^{3}/\Gamma) \neq \emptyset$ and so $Z$ represents a 
2-cycle in $H_{2}(S^{3}/\Gamma, {\mathbb F})$. By the preceding argument, 
$Z$ bounds a 3-chain $W'$ in $E$. Hence $(W\cap E) + W'$ is a 3-chain in 
$E$, with boundary $\Sigma$, i.e.~$[\Sigma] = 0$ in $H_{2}(E, {\mathbb F})$, 
giving a contradiction. 

 Thus, the hypothesis (3.17) rules out any possible orbifold 
degeneration of the sequence $(M, g_{i})$. Having ruled out the 
possibilities (ii) and (iii), it follows that (i) must hold. The fact 
that $N$ is weakly embedded in $M$ follows immediately from Remark 3.6 
and the definition of smooth convergence preceding (3.8). This 
completes the proof of Theorem 3.7.
{\endproof}

\begin{remark} \label{r 3.10.}{\bf (i).}
  {\rm Note that Theorem 3.7 does not imply that the limit manifold $N$ 
topologically equals $M$. For instance,  it is possible at this stage 
that some of the topology of $M$ escapes to infinity, and is lost in 
the limit $N$. All of the results and discussion in \S 3.2 hold without 
change if the 4-manifold $M$ is replaced by a sequence of 4-manifolds 
$M_{i}$, with a common boundary $\partial M_{i} = \partial M$, so that 
one has $(M_{i}, g_{i})$ in place of $(M, g_{i})$. 

{\bf (ii).}
 We also observe that the proof shows that there is a uniform bound 
$\Lambda  = \Lambda (\Lambda_{0},d,D,T,\{\gamma_{i}\}) < \infty$ such 
that, for $x$ satisfying (3.22),
\begin{equation} \label{e3.34}
|R_{g_{i}}|(x) \leq  \Lambda , 
\end{equation}
i.e. there is a uniform bound on the sectional curvatures of 
$\{g_{i}\}$ in this region. 

{\bf (iii).}
 The existence and behavior of $(M, g_{i})$ at base points $x_{i}$ 
where $dist_{\bar g_{i}}(x_{i}, \partial M) \rightarrow  \infty $ will 
be discussed at the end of \S 4.}
\end{remark}

\section{Compactness II: Boundary Behavior.}
\setcounter{equation}{0}

  Theorem 3.7 essentially corresponds to a uniform interior regularity 
result in that one has uniform control of the AHE metrics $g_{i}$ in 
the interior of the 4-manifold $M$, i.e.~on compact subsets of $M$. In this 
subsection, we extend this to similar control in a neighborhood of definite 
size about the boundary. The Chern-Gauss-Bonnet theorem in (3.12), cf.~also (3.21), 
plays a crucial role in obtaining control in the interior. A significant point in 
this section is to prove that such an apriori $L^{2}$ bound on the curvature is not 
necessary to obtain control near the boundary; thus we show that control of the 
boundary metric itself gives control of the Einstein metric in a neighborhood of 
the boundary of definite size, cf. Corollary 4.10. Much of the proof of this result 
holds in fact in all even dimensions and conjecturally in all dimensions. 

 We begin with the following result, most of which was proved in [6], cf.~also [7]. 
We fix $p > 4$ and denote the $L^{p}$ curvature radius by $\rho$ instead of $\rho^{p}$. 

\begin{proposition} \label{p 4.1.}
  On a fixed $4$-manifold $M$, let $g \in E_{AH}^{m,\alpha}$, $m \geq 4$, with 
boundary metric $\gamma$ satisfying $||\gamma||_{C^{m,\alpha}} \leq K$ in 
a fixed coordinate atlas for $\partial M$. Let $U_{\delta} = \{x\in\bar M: 
dist_{\bar g}(x, \partial M) \leq \delta\}$ where $\bar g$ is the geodesic 
compactification of $(M, g)$ with boundary metric $\gamma$. Suppose the 
$L^{p}$ curvature radius of $\bar g$ satisfies
\begin{equation} \label{e4.1}
\rho (x) \geq  \rho_{0}, 
\end{equation}
for some constant $\rho_{0} > 0$, for all $x \in  U_{\delta_{0}}$, 
some $\delta_{0} > 0$. 

 Then there is a $\delta_{1} = \delta_{1}(\rho_{0}, \delta_{0}, K) > 0$ 
such that the $C^{m-1,\alpha}$ geometry of the geodesic compactification 
$\bar g$ is uniformly controlled in $U_{\delta_{1}}$, in that 
\begin{equation}\label{e4.2}
r_{h}^{m-1,\alpha} \geq r_{0},
\end{equation}
where $r_{0}$ depends only on $\rho_{0}$, $\delta_{0}$, and $K$. The same 
statement holds if $m = \infty$ or $m = \omega$.
\end{proposition}

\noindent
{\bf Proof:}
 By the discussion in \S 2, boundary regularity implies that $(M, g)$ has a 
$C^{m,\alpha}$ conformal compactification $\widetilde g$, with boundary metric 
$\gamma$. The metric $\widetilde g$ may be chosen to be of constant scalar 
curvature $\widetilde s = const$ and $\widetilde g$ is $C^{m,\alpha}$ in local 
harmonic coordinates. Again as noted in \S 2, the geodesic compactification $\bar g$ 
of $g$ with boundary metric $\gamma$ is then at least $C^{m-1,\alpha}$.
Moreover, by [6, Prop.~2.1, Thm.~2.4 and Prop.~2.7], the $C^{m-1,\alpha}$ geometry 
of the metric $\bar g$ is uniformly controlled in $U_{\delta_{1}}$, i.e.~(4.2) holds, 
for $\delta_{1} = \delta_{1}(\rho_{o}, \gamma)$, provided the $C^{m,\alpha}$ geometry 
of the boundary metric $\gamma$ is uniformly bounded, (4.1) holds, and provided 
there is a uniform lower bound on the volume radius of $\bar g$ in $U_{\delta_{1}}$. 

 The control on the boundary metric and (4.1) are given by hypothesis. 
To obtain a lower bound on the volume radius $\nu$, observe that Lemma 3.8, 
cf.~also (3.27), gives a lower bound on the volume ratios of balls, $vol 
B_{q}(s)/s^{4}$ centered at points $q\in\partial M,$ depending only on 
the $C^{0}$ geometry of the boundary metric. In particular, there is a 
constant $\nu_{0} > 0$ such that $vol B_{q}(\rho_{0})/\rho_{0}^{4} 
\geq  \nu_{0}$, $q\in\partial M$, for all $(M, \bar g_{i})$. Within the 
$L^{p}$ curvature radius, i.e. for balls $B_{x}(r) \subset  
B_{q}(\rho_{0})$, standard volume comparison results imply then a 
lower bound $vol B_{x}(r)/r^{4} \geq  \nu_{1}$, where $\nu_{1}$ depends 
only on $\nu_{0}$, cf. [4,\S 3]. This gives a uniform lower bound on 
the volume radius within $U_{\rho_{0}}$, and hence the result follows.

{\endproof}

  Proposition 4.1 also holds locally, in that if (4.1) holds at all points 
in an open set $V$ containing a domain of fixed size in $\partial M$, then 
the conclusion holds in $V\cap U_{\delta_{1}}$. Also, this result holds for 
$m \geq 3$. Of course the arguments concerning the lower bound on the volume 
radius hold in all dimensions. 

\medskip

   Proposition 4.1 should be understood as a natural strengthening of the boundary 
regularity property that AH Einstein metrics with $C^{m,\alpha}$ boundary metric 
have a $C^{m,\alpha}$ conformal compactification. Thus, it gives uniform control 
or stability of the metric $g\in E_{AH}^{m,\alpha}$ near the boundary, depending 
only on the boundary metric and a lower bound on $\rho$ as in (4.1). The main 
result of this section, Theorem 4.4 below, removes the assumption (4.1) and proves 
the estimate (4.2) independent of $\rho_{0}$ and $\delta_{0}$. Of course this could 
not be true locally; Theorem 4.4 is necessarily a global result. 

  The weak hypothesis (4.1) on the curvature radius may obviously be replaced 
by a stronger assumption. Thus, Proposition 4.1 implies, for instance, that in 
the geodesic compactification,  
\begin{equation}\label{e4.3}
r_{h}^{k,\beta} \geq r_{0} \Rightarrow r_{h}^{m-1,\alpha} \geq r_{1},
\end{equation}
for any $2 < k + \beta < m-1 + \alpha$ with $r_{1} = r_{1}(r_{0}, \delta_{0}, K)$. 
Thus, control of the metric in a weak norm, such as $L^{k,p}$ or $C^{k,\beta}$ 
implies control of the metric in a stronger norm, governed by the regularity of 
the metric at the boundary. We formalize this as follows:

\begin{definition}\label{d4.2}
{\rm In any dimension, an AH Einstein metric $g \in E_{AH}$ satisfies the 
{\it strong control property} if (4.3) holds, for some $(k, \beta)$ and 
$(m, \alpha)$ with $k + \beta < m-1 + \alpha$, and $m-1 \geq n$. }
\end{definition}

  The proof of boundary regularity in [6], (cf.~also [7]), uses the fact that AH 
Einstein metrics and their conformal compactifications satisfy the conformally invariant 
Bach equation in dimension 4. In a suitable gauge, this is a non-linear elliptic system 
of PDE and boundary regularity follows essentially from standard boundary regularity for 
such elliptic systems. The strong control property (4.3), (or its original version in 
Proposition 4.1), then follows from the fact that concurrent with boundary regularity, 
one has uniform elliptic estimates for solutions of the Bach equations near and up to 
the boundary, given weak control in such a region. 

\begin{remark}\label{r4.3}
{\rm It will be noted in the proofs of the results below that all of the results of 
this section hold in all dimensions in which the strong control property holds. 
A recent result of Helliwell [25] proves boundary regularity for AH Einstein metrics 
in all even dimensions, by generalizing the proof in dimension 4 in [6], [7]. The idea is 
to use the Fefferman-Graham obstruction tensor in place of the Bach equations in higher 
dimensions; again in a suitable gauge, this gives a non-linear elliptic system. The proof 
given in [25] requires that $g$ be $C^{n, \alpha}$ conformally compact, although we 
suspect that it suffices for $g$ to be $C^{2,\alpha}$ conformally compact. In any case, 
it is straightforward to show from [25] that the strong control property (4.3) holds in 
all even dimensions with $n < k + \beta < m-1 + \alpha$. 

 In odd dimensions, conformal compactifications of AH Einstein metrics are not $C^{n}$ 
up to the boundary, exactly due to the presence of the obstruction tensor ${\mathcal H}$. 
Nevertheless, one can define the spaces $C^{m,\alpha}$, $m \geq n$ in a polyhomogeneous 
sense. We conjecture that with such a suitable modification of the definition of harmonic 
radius, the strong control property holds in all odd dimensions; this is an interesting 
open question. }
\end{remark}

  The main result of this section is the following:

\begin{theorem} \label{t 4.4.}
  Let $(M, g)$ be an AH Einstein metric on a $4$-manifold with a $L^{2,p}$ conformal 
compactification, $p > 4$, and boundary metric $\gamma$. Suppose that $\gamma \in 
C^{m,\alpha}$, $m \geq 4$ and that, in a fixed coordinate system for $\partial M$, 
$||\gamma||_{C^{m,\alpha}} \leq K$, for some $\alpha > 0$. If $\bar g$ is the 
associated $C^{m-1,\alpha}$ geodesic compactification of $g$ determined by $\gamma$, 
then there are constants $\mu_{0} > 0$ and $\rho_{0} > 0$, depending only on $K$, 
$\alpha$, and $p$, such that
\begin{equation} \label{e4.4}
\rho (x) \geq  \rho_{0}, \ \ {\rm for \ all} \ x \ {\rm with} \ \ t(x) \leq  \mu_{0}. 
\end{equation}

  The same result holds in all dimensions in which the strong control property holds, 
with $m \geq n+1$. 
\end{theorem}

   The proof is rather long, and so is broken into several Propositions and Lemmas. 
Before beginning, it is worth pointing out that the curvature of $\bar g$ {\it  at} 
$\partial M$ is uniformly bounded by the $C^{2}$ geometry of the boundary metric 
$\gamma$, cf.~(A.8)-(A.10) in the Appendix. Thus, the geometry of $g$ at $\partial M$ 
is well-controlled to $2^{\rm nd}$ order; one needs to extend this to control of $g$ near 
$\partial M$.

  We will give the proof in all dimensions, assuming the strong control property 
holds. Also, we work with the $L^{p}$ curvature radius $\rho = \rho^{p}$, $p > n$, 
but one may work equally well with the $L^{k-2,q}$ curvature or $L^{k,q}$ harmonic 
radius; the arguments are exactly the same for these radii. 

   Let $\tau_{0}$ be the distance to the cutlocus of the normal exponential map from 
$\partial M$; $\tau_{0}(x) = dist_{\bar g}(x, \bar C)$, where $\bar C$ is the 
cutlocus, as in (3.19). The first step is to prove that the function $\rho$ is 
bounded below by $\tau_{0}$ near $\partial M$. The second, more difficult, step is to 
prove that $\tau_{0}$ itself is uniformly bounded below near $\partial M$.

\medskip

 We begin with the following result.
\begin{proposition} \label{p 4.5.}
  For $(M, g)$ as in Theorem $4.4$, there is a constant $c_{0} > 0$, depending only on 
$K$, $\alpha$ and $p$, such that
\begin{equation} \label{e4.5}
\rho(x) \geq  c_{0}\tau_{0}(x), \ \ {\rm for \ all} \ x \ {\rm with} \ \  t(x) \leq  1. 
\end{equation}
\end{proposition}

\noindent
{\bf Proof:}
 The proof is by contradiction. If (4.5) is false, then there must exist a sequence 
of AH Einstein metrics $g_{i}$ on $M_{i}$, with $C^{n,\alpha}$ (polyhomogeneous) geodesic 
compactifications $\bar g_{i}$, for which the boundary metrics $\gamma_{i}$ satisfy
\begin{equation} \label{e4.6}
||\gamma_{i}||_{C^{n,\alpha}} \leq  K, 
\end{equation}
but for which the $L^{p}$ curvature radius $\rho = \rho_{i}$ of $\bar g_{i}$ 
satisfies 
\begin{equation} \label{e4.7}
\rho(x_{i})/\tau_{0}(x_{i}) \rightarrow  0 \ \ {\rm as} \ \  i \rightarrow  \infty , 
\end{equation}
on some sequence of points $x_{i}\in (\bar M_{i}, \bar g_{i})$. Note that the ratio in 
(4.7) is scale-invariant. Without loss of generality, assume that $x_{i}$ realizes 
the minimal value of $\rho/\tau_{0}$ on $(M_{i}, \bar g_{i})$. Of course $x_{i}$ may 
occur at $\partial M$. 

 Now blow-up or rescale the metrics $\bar g_{i}$ at $x_{i}$ to make $\rho (x_{i}) = 1$, 
i.e. set
\begin{equation} \label{e4.8}
g_{i}'  = \lambda_{i}^{2} \bar g_{i}, 
\end{equation}
where $\lambda_{i} = \rho (x_{i})^{-1} \rightarrow  \infty$. Let $\rho'  = 
\lambda_{i}\rho$ and $\tau_{0}'  = \lambda_{i}\tau_{0}$ be the $L^{p}$ curvature 
radius and distance to the cutlocus, with respect to $g_{i}'$. Then 
\begin{equation} \label{e4.9}
\rho' (x_{i}) = 1, 
\end{equation}
and by the minimality property of $\rho (x_{i})$,
\begin{equation} \label{e4.10}
\rho' (y_{i}) \geq  \rho' (x_{i})\frac{\tau_{0}' (y_{i})}{\tau_{0}(x_{i})} = 
\frac{\tau_{0}' (y_{i})}{\tau_{0}(x_{i})}. 
\end{equation}
for any $y_{i}\in \bar M$. Since (4.7) implies $\tau_{0}'(x_{i}) \rightarrow \infty$, 
it follows that 
\begin{equation} \label{e4.11}
\rho' (y_{i}) \geq  \frac{1}{2}, 
\end{equation}
for all $y_{i}$ with $dist_{g_{i}'}(x_{i}, y_{i}) \leq D$, for any given $D$, provided 
$i$ is sufficiently large. Of course one also has $\tau_{0}' (y_{i}) \rightarrow \infty$ 
for any such $y_{i}$.

 It follows that the metrics $g_{i}' $ have uniformly bounded curvature, on the 
average in $L^{p}$ on all unit balls of $g_{i}'$-bounded distance to $x_{i}$. 
For clarity, we divide the proof into two cases, according to whether 
$dist_{g_{i}'}(x_{i}, \partial M)$ remains bounded or not.

\medskip

{\bf Case I. Bounded distance.}

 Suppose $dist_{g_{i}'}(x_{i}, \partial M) \leq D$, for some $D < \infty$. We 
may then apply Proposition 4.1, or more precisely its local version as following 
the statement of Proposition 4.1, to conclude that a subsequence of the metrics 
$(M_{i}, g_{i}', x_{i})$ converges in the $C^{n,\alpha'}$ (or stronger) topology to a 
limit $(N, g' , x_{\infty})$. (Of course the convergence is modulo diffeomorphisms). 
Since $\tau_{0}' (y_{i}) \rightarrow  \infty$ for any $y_{i}$ within bounded 
$g_{i}'$-distance to $x_{i}$, the limit $(N, g')$ is a complete manifold with 
boundary. Moreover, since the unscaled boundary metrics $\gamma_{i}$ are uniformly 
bounded in $C^{n,\alpha}$, the limit $(\partial N, \gamma') = ({\mathbb R}^{n}, \delta)$, 
where $\delta$ is the flat metric.

  Since the $L^{p}$ curvature radius is continuous in the $C^{n,\alpha'}$ topology, 
it follows from (4.9) that
\begin{equation} \label{e4.12}
\rho' (x_{\infty}) = 1. 
\end{equation}
The rest of the proof in this case is to prove that $(N, g')$ is flat. This clearly 
contradicts (4.12), and so will complete the proof.

 Let
\begin{equation} \label{e4.13}
t' (x) = dist_{g'}(x, \partial N) = \lim_{i\rightarrow\infty}t_{i}' (x_{i}), 
\end{equation}
where $t_{i}'$ is the geodesic defining function for $(\partial M_{i}, g_{i}' )$ and 
$x_{i} \rightarrow x$. Since $\partial M$ is totally geodesic in $(M, g_{i}')$, 
the smooth convergence implies $\partial N = {\mathbb R}^{n}$ is totally geodesic 
in $(N, g')$. Since $\tau_{0}'  = \infty$ on $N$, $t'$ is globally defined and smooth 
on $N$.

 We now use some of the curvature properties of geodesic compactifications given 
in the Appendix. The equations (A.1)-(A.3) for the curvatures $R'$, $Ric'$ and $s'$ 
also hold on $(N, g')$, with $t'$ in place of $t$. We note that the Hessian 
$D^{2}t' = A$, where $A$ is the $2^{\rm nd}$ fundamental form of the level sets of 
$t'$, and $\Delta t' = H$, the mean curvature of the level sets.

 By (A.3), the scalar curvature $s'$ of $(N, g')$ is given by
\begin{equation} \label{e4.14}
s'  = -2n\frac{\Delta t'}{t'}, 
\end{equation}
and by (A.11) satisfies
\begin{equation} \label{e4.15}
\dot s' \geq \frac{1}{2n^{2}}t'(s')^{2}, 
\end{equation}
where $\dot s'  = \partial s' /\partial t'$. Since $s_{\gamma'} = 0$, (A.8) implies 
that $s'  = 0$ on $\partial N$, and hence by (4.15), $s'  \geq 0$ everywhere on $N$. 
Moreover, elementary integration of (4.15) implies that if, along a geodesic 
$\sigma  = \sigma (t')$ normal to $\partial N$, $s'(t_{0}) >  0$, for some $t_{0} > 0$, 
then $(t')^{2} \leq t_{0}^{2} + 4n^{2}/(s' (t_{0}))$, which is impossible, 
since $t'  \rightarrow \infty$ along all geodesics normal to $\partial N$. It follows 
that one must have
\begin{equation} \label{e4.16}
s'  \equiv  0 \ \ {\rm on} \ \ (N, g'). 
\end{equation}
From (4.14) and (4.16), we then have
\begin{equation} \label{e4.17}
\Delta t'  = 0, 
\end{equation}
so that $t'$ is a smooth harmonic function on $N$, with $|\nabla t'| = 1$.

 Now the Riccati equation (A.7) for the $t'$ geodesics on $(N, g')$ gives
\begin{equation} \label{e4.18}
|A|^{2} + Ric(\nabla t',\nabla t') = 0, 
\end{equation}
On the other hand, the formula (A.2) for the Ricci curvature holds, 
so that on $(N, g')$,
\begin{equation} \label{e4.19}
Ric = - (n-1)t^{-1}D^{2}t -  t^{-1}(\Delta t) g' = - (n-1)t^{-1}D^{2}t. 
\end{equation}
Here and below, we drop the prime from the notation. Clearly $D^{2}t(\nabla t,\nabla t) 
\equiv  0$ on $N$, and so by (4.18), $|A|^{2} = 0$, i.e. all the level 
sets are totally geodesic. Since $A = D^{2}t$, (4.19) implies that $Ric \equiv 0$, 
so that $(N, g')$ is Ricci-flat. The vector field $\nabla t$ is thus a parallel 
vector field, and so $(N, g')$ splits as a product along the flow lines of 
$\nabla t$. Since $\partial N = \{t=0\}$ is flat ${\mathbb R}^{n}$, it follows that 
$(N, g')$ is flat, as claimed. 

\medskip

{\bf Case II. Unbounded distance.}

 Suppose $t_{i}'(x_{i}) = dist_{g_{i}'}(x_{i}, \partial M) \rightarrow  \infty$, 
as $i \rightarrow \infty$. Using (4.7) again, exactly the same arguments as in 
(4.14)-(4.19) applied to $(M_{i}, g_{i}', x_{i})$ show that, for $i$ sufficiently 
large, $(M_{i}, g_{i}')$ almost splits in regions of bounded diameter about $x_{i}$, 
in that such regions are topologically products and the flow by the integral curves 
of $t_{i}'$ are almost isometries of $g_{i}'$. 

  More precisely, suppose first that $(M_{i}, g_{i}', x_{i})$ does not collapse near 
$x_{i}$. Then, exactly as in Case I, by the strong control property, 
$(M_{i}, g_{i}', x_{i})$ is close in the $C^{n,\alpha'}$ topology, to a product 
metric in regions of $g_{i}'$-bounded diameter about $x_{i}$. Suppose instead 
$(M_{i}, g_{i}', x_{i})$ collapses near $x_{i}$. The strong control property 
implies the collapse is with bounded curvature. Clearly the collapse, which is caused 
by short geodesic loops, is transverse to the flow lines of $t_{i}'$. By the general 
collapse theory, cf.~[15], the collapse may be unwrapped in local covering spaces. Thus, 
there are covering spaces of balls of fixed but small diameter about $x_{i}$ such that 
the lifted metrics on the covers do not collapse, and so converge in $C^{n,\alpha'}$ to 
a local limit. All of the estimates (4.14)-(4.19) remain valid on such local covers, and 
since these are pointwise estimates transverse to the collapsing directions, it follows 
again that the original metrics $(M_{i}, g_{i}', x_{i})$ are close in $C^{n,\alpha'}$ 
to a product metric in regions of $g_{i}'$-bounded diameter about $x_{i}$. 

  This shows that $\rho_{i}'$ is almost constant in such regions, and, (again by the 
strong control property),
\begin{equation}\label{e4.20}
\frac{\partial \rho_{i}'}{\partial t_{i}'} \sim 0. 
\end{equation}
Note that $\partial \rho / \partial t$ is scale-invariant, (as are the ratios 
$\rho / t$ and $\rho / \tau_{0}$). The argument above was carried out at the base 
points $x_{i}$. However, exactly the same reasoning shows that (4.20) 
holds at all points $p_{i}$ where $(\rho / \tau_{0})(p_{i}) < < 1$. Now integrate 
(4.20) along integral curves $\sigma_{i}$ of $t_{i}'$ starting at points $y_{i}$ 
of bounded $g_{i}'$-distance to $x_{i}$ in the direction toward $\partial M$. 
Since $\rho_{i}'(y_{i}) \sim 1$ and  $t_{i}'(y_{i}) \rightarrow \infty$, this 
implies in particular that the scale-invariant ratio $\rho /t$ on $(M_{i}, g_{i}')$ 
satisfies
\begin{equation}\label{e4.20a}
\frac{\rho_{i}'}{t_{i}'} < < 1,
\end{equation}
along such curves, all the way to $\partial M$, since $\rho/\tau_{0} <<1$ along such 
curves. However, this is clearly impossible, since again the arguments 
in Case I imply that $(M_{i}, g_{i}')$ is almost flat within bounded $g_{i}'$-distance 
to $\partial M$, so that $\rho_{i}' \geq (1-\delta)t_{i}'$ in such regions, where 
$\delta$ is small for $i$ large. In other words, \eqref{e4.20} and \eqref{e4.20a} 
prevent the ratio $\rho / t$, which is arbitrarily small around $x_{i}$, from 
increasing to near 1 near $\partial M$. This contradiction completes the proof.  

{\endproof}

  Note that the opposite inequality to (4.5) also holds, i.e.
\begin{equation}\label{e4.21}
\rho(x) \leq \tau_{0}(x).
\end{equation}
This is essentially a tautology, since the metric $\bar g$ is singular at the 
cutlocus and so $\rho$-balls are not defined past the cutlocus. The proof of 
Proposition 4.5 does not use the strong control property with $m-1 \geq n$. 
Basically, it suffices to have this property in any norm stronger than the 
$L^{2,p}$ norm, so that the $L^{p}$ curvature radius is continuous. However, 
the condition $m-1 > n$ will be crucial in the next result. 

   The next step in the proof of Theorem 4.3 is to prove that $\tau_{0}$ 
is uniformly bounded below near $\partial M$.
\begin{proposition} \label{p 4.6.}
  For $(M, g)$ as in Theorem $4.4$, there is a constant $\mu_{1} > 0$, depending only 
on $K$, $\alpha$ and $p$, such that, for any $x\in\partial M,$
\begin{equation} \label{e4.22}
\tau_{0}(x) \geq  \mu_{1}. 
\end{equation}
\end{proposition}

\noindent
{\bf Proof:}
 The proof of this estimate is much more subtle than that of Proposition 4.5. 
Propositions 4.1 and 4.5 are essentially local results; they hold for AH Einstein 
metrics defined only in a geodesic collar neighborhood of $\partial M$, $t \in 
[0, \varepsilon_{0}]$, for some $\varepsilon_{0} > 0$. However, Proposition 4.6, 
(and hence Theorem 4.4), is global; it requires $g$ to be an AH Einstein metric 
defined on a compact manifold $M$, with conformal boundary $\partial M$ at infinity. 
 
  We first prove some preliminary lemmas, which are basically straightforward consequences 
of the strong control property. The actual proof begins after the proof of Lemma 4.9; 
on a first reading, one might want to start at this point to understand how the Lemmas 
are used. The crux of the proof is the use of the relation (2.22) in Proposition 2.4. 
It is this relation in particular, (as well as Lemma 4.7), which requires $(M, g)$ 
is global. 

\medskip

  To begin, we examine the geometry of the geodesic compactification $\bar g$ 
in regions where $\tau_{0}$ is very small, in fact possibly arbitrarily small. 
Thus suppose $(M_{i}, g_{i})$ is a sequence of AH Einstein metrics with boundary 
metrics $\gamma_{i}$ on $\partial M_{i} = \partial M$ uniformly controlled in 
$C^{m,\alpha}$, $m \geq n+1$. Hence, a subsequence of $\gamma_{i}$ converges in 
$C^{m,\alpha'}$ to a limit $C^{m,\alpha}$ boundary metric $\gamma$. Let 
$x_{i}$ be any sequence of points in $M_{i}$ such that 
\begin{equation}\label{e4.23}
\tau_{0}(x_{i}) \rightarrow 0.
\end{equation}
Rescale the compactified metrics $\bar g_{i}$ to 
\begin{equation}\label{e4.24}
g_{i}' = (\tau_{0}(x_{i}))^{-2}\bar g_{i}, 
\end{equation}
so that 
\begin{equation}\label{e4.25}
\tau_{0}'(x_{i}) = 1.
\end{equation}
By Proposition 4.5 and \eqref{e4.21}, $\rho'(x_{i})$ is uniformly bounded 
above and below. Hence as in the proof of Proposition 4.5, a subsequence of 
$(M, g_{i}', x_{i})$ converges in the $C^{n,\alpha'}$ topology (modulo diffeomorphisms), 
to a maximal limit $(N, g', x_{\infty})$; maximal here means the maximal connected 
domain on which the convergence is $C^{n,\alpha'}$. The manifold $N$ has 
boundary-at-infinity $\partial N$ a domain in ${\mathbb R}^{n}$, with boundary 
metric the flat metric $\gamma' = \delta$, and containing the ball 
$B_{x_{\infty}}(1)$ in ${\mathbb R}^{n}$. However, in this case, $N$ is 
not complete away from its boundary, since $\tau_{0}'(x_{\infty}) \leq 1$.

  Consider now also the sequence of AH Einstein metrics $g_{i}$ on $M_{i}$, with conformal 
compactification $\bar g_{i}$, with base points $x_{i}$ chosen above. Of course 
$dist_{\bar g_{i}}(x_{i}, \partial M) \rightarrow 0$. The metrics $\bar g_{i}$ are 
rescaled up to $g_{i}' = \lambda_{i}^{2}\bar g_{i}$, $\lambda_{i} = (\tau_{0}(x_{i}))^{-1}$, 
and converge, (in a subsequence), in $C^{n,\alpha'}$ to the maximal limit $(N, g', x)$. Now 
the rescaling of $\bar g_{i}$ just corresponds to changing the defining function $t_{i}$ 
to $t_{i}' = \lambda_{i}t_{i}$. Thus it does not change $g_{i}$ itself; in fact 
$g_{i} = (t_{i}')^{-2}g_{i}'$. Since $t_{i}' \rightarrow t'$ and $g_{i}' \rightarrow g'$ 
on $N$, it follows that the pointed sequence of AH Einstein metrics $(M_{i}, g_{i}, x_{i})$ 
converges to a limit AH Einstein metric $({\mathcal N}, g_{\infty}, x_{\infty})$. Again here 
${\mathcal N}$ is the maximal connected domain containing $x_{\infty}$ on which the convergence 
is $C^{n,\alpha'}$; in fact the convergence is now $C^{\infty}$ smooth on compact subsets, 
by regularity properties of Einstein metrics. The limit $({\mathcal N}, g_{\infty})$ is 
of course not conformally {\it compact} in general; the ``compactification'' of 
$({\mathcal N}, g_{\infty})$ by $t$ gives the manifold $(N, g')$, i.e.
\begin{equation}\label{e4.26}
g_{\infty} = t^{-2}g',
\end{equation}
where $N \subset {\mathcal N}$ is the region where $t$ is smooth.

\medskip

  The first result is an analogue of (4.12) in this setting. 
\begin{lemma}\label{l4.7}
Any blow-up limit manifold $(N, g', x_{\infty})$ is not flat.
\end{lemma}

\noindent
{\bf Proof:} If $(N, g')$ is flat, then by \eqref{e4.26}, $({\mathcal N}, g_{\infty})$ is 
hyperbolic, i.e.~of constant curvature -1, and so locally embeds in the hyperbolic 
space ${\mathbb H}^{n+1}(-1)$. By construction, the maximal limit $({\mathcal N}, 
g_{\infty})$ is not complete; it has a ``boundary'' corresponding to the points 
where $\rho_{g_{i}} \rightarrow 0$. However, the pointed sequence $(M_{i}, g_{i}, x_{i})$ 
of complete manifolds has a uniform lower bound on Ricci curvature, since the metrics 
are Einstein, and is non-collapsing at points within bounded distance to $x_{i}$. 
By the Gromov weak compactness theorem [22], a subsequence of $(M_{i}, g_{i}, x_{i})$ 
converges in the pointed Gromov-Hausdorff topology to a limit $(\bar {\mathcal N}, d, x)$, 
where $(\bar{\mathcal N}, d)$ is a {\it complete}, non-compact, length space. In particular, 
$(\bar{\mathcal N}, d)$ is geodesically complete; any two points may be joined by a 
minimizing geodesic and all geodesic balls $B_{x}(r)$ in $\bar{\mathcal N}$ have compact 
closure strictly contained in $\bar{\mathcal N}$. 

   Clearly the smooth domain $({\mathcal N}, g_{\infty})$ embeds in $(\bar{\mathcal N}, 
d)$. By a result of Cheeger-Colding [13], the singular set $\partial {\mathcal N} =  
\bar{\mathcal N} \setminus {\mathcal N}$ of $\bar{\mathcal N}$ is of codimension 
2 in $\bar{\mathcal N}$. In particular the metric boundary $\partial {\mathcal N}$ 
is not a topological boundary, (which would have codimension 1). 

  By the volume comparison theorem (3.31) on $(M_{i}, g_{i})$, the ratio 
$vol B^{i}_{z}(r)/vol B_{-1}(r)$ is monotone non-increasing in $r$, where 
$B^{i}_{z}(r)$ is the geodesic $r$-ball about any $z$ in $(M_{i}, g_{i})$ and 
$vol B_{-1}(r)$ is the volume of the geodesic $r$-ball in ${\mathbb H}^{n+1}(-1)$. 
Further, for any given $i$, $\lim_{r\rightarrow 0}vol B^{i}_{z}(r)/vol B_{-1}(r) = 1$. 
A fundamental result of Colding [14], cf.~also [13], shows that the volume of geodesic 
balls is continuous under Gromov-Hausdorff limits so that in $(\bar{\mathcal N}, d)$, 
$vol B_{z}(r)/vol B_{-1}(r)$ is also monotone non-decreasing, and hence 
$\lim_{r\rightarrow 0}vol B_{z}(r)/vol B_{-1}(r) \leq 1$, with equality for 
$z \in {\mathcal N}$. 

   Now suppose first for simplicity of argument that ${\mathcal N}$ is simply connected.
Since ${\mathcal N}$ is of constant curvature $-1$, the developing map based at 
any point in ${\mathcal N}$ gives an isometric immersion $F: {\mathcal N} 
\rightarrow H^{n+1}(-1)$. Observe that the image of $F$ is of full measure in 
${\mathbb H}^{n+1}(-1)$; this is because $\bar{\mathcal N}$ is complete, so 
that geodesics do not terminate, and $\partial{\mathcal N}$ is of codimension 
2 and so of measure 0. In fact $F$ extends to a continuous map of 
$\bar{\mathcal N}$ onto ${\mathbb H}^{n+1}(-1)$. 

  It follows that $F$ maps $B_{z}(r)$ onto $B_{F(z)}(r) \subset {\mathbb H}^{n+1}(-1)$ 
modulo sets of measure 0. Moreover, $F$ preserves the volumes of these balls 
when counting multiplicities of the image, so that
\begin{equation}\label{e4.27}
\frac{vol B_{z}(r)}{vol B_{-1}(r)} \geq \frac{vol F(B_{z}(r))}{vol B_{-1}(r)} = 1.
\end{equation}
Since by the monotonicity above one has
\begin{equation}\label{e4.28}
\frac{vol B_{z}(r)}{vol B_{-1}(r)} \leq 1,
\end{equation}
for all $r$, it is then immediate that
\begin{equation}\label{e4.29}
\frac{vol B_{z}(r)}{vol B_{-1}(r)} = 1,
\end{equation}
for all $r$. This of course implies that $\bar{\mathcal N}$ is isometric to 
${\mathbb H}^{n+1}(-1)$, by the volume rigidity theorem for Ricci curvature. 
Moreover, by [2, Thm.~3.2], the sequence $(M_{i}, g_{i}, x_{i})$ now converges 
smoothly to $(\bar{\mathcal N}, g_{\infty}, x_{\infty})$ everywhere. In particular, 
$\tau_{0} = \infty$, contradicting the fact that $\tau_{0}(x_{\infty}) = 1$ from 
\eqref{e4.25}. 

   The proof is similar when ${\mathcal N}$ is not simply connected. Thus, 
let $\widetilde{\mathcal N}$ be the universal cover of ${\mathcal N}$. (For illustration, 
it is useful to picture $\widetilde{\mathcal N}$ as a branched cover of $\bar{\mathcal N}$ 
branched over the singular locus $\partial{\mathcal N}$). Then the developing 
map $F$ is an isometric immersion of $\widetilde{\mathcal N}$ into ${\mathbb H}^{n+1}(-1)$. 
Let $D$ be a Dirichlet fundamental domain for the action of $\pi_{1}({\mathcal N})$ 
on $\widetilde{\mathcal N}$, based at a point $z \in {\mathcal N}$. Modulo sets 
of measure 0, $D$ may be identified with ${\mathcal N}$, and $F$ gives an isometric 
immersion $F: D \rightarrow {\mathbb H}^{n+1}(-1)$. Now however $\partial D$ is 
of codimension 1 and $D$, (or $\bar D$), is not geodesically complete; the action of 
$\pi_{1}({\mathcal N})$ identifies subsets of $\partial D$ to obtain the manifold 
${\mathcal N}$ or its completion $\bar {\mathcal N}$. Nevertheless, letting 
$S_{z}(r)$ be the geodesic $r$-sphere about $z$ in $N$, or equivalently in $D$ 
modulo sets of measure 0, and letting $S_{F(z)}(r)$ be the corresponding sphere 
in ${\mathbb H}^{n+1}(-1)$, we claim that 
\begin{equation}\label{e4.29a}
\limsup_{r\rightarrow \infty}\frac{vol S_{z}(r)}{vol_{-1}(S_{F(z)}(r))} \geq 1. 
\end{equation}
To prove \eqref{e4.29a}, recall the construction of the blow-up limit ${\mathcal N}$ and its 
compactification $N$ with conformal infinity given by the domain $\partial N \subset 
{\mathbb R}^{n}$. View the codimension 2 singular set $S = \partial {\mathcal N}$ as 
a subset of the compactification $\bar N$, and suppose first for clarity that $S$ intersects 
the full conformal boundary ${\mathbb R}^{n}$ in a closed set $S_{\infty}$ also of 
codimension 2 in ${\mathbb R}^{n}$. Then $\partial N =  {\mathbb R}^{n} \setminus 
S_{\infty}$, and so in particular $\partial N$ is of full measure in ${\mathbb R}^{n}$. 
As $r \rightarrow \infty$, the geodesic spheres $S_{z}(r)$ in $N$, (or $D$), tend to 
$\partial N$. Then, as in the proof of \eqref{e4.27}, \eqref{e4.29a} follows from the 
fact that $\partial N$ is a set of full measure in the conformal compactification 
$S^{n} = {\mathbb R}^{n} \cup \{\infty\}$ of ${\mathbb H}^{n+1}(-1)$. 

   Now the same argument holds regardless of the exact structure of the singular 
set $S_{\infty} \subset {\mathbb R}^{n}$, since it is still the case that $S_{z}(r) 
\subset N$ is of full measure in $\bar S_{z}(r) \subset \bar N$ and $\lim_{r\rightarrow 
\infty}\bar S_{z}(r) = {\mathbb R}^{n}$. 

  By integration, the estimate \eqref{e4.29a} also holds for balls $B_{z}(r)$ and the 
relation \eqref{e4.29} then follows from \eqref{e4.28} as before. 
{\endproof}

  We point out that the same proof shows that $(N, g', x_{\infty})$ cannot be 
conformally flat. This is because \eqref{e4.26} would then imply that the metric 
$({\mathcal N}, g_{\infty})$ is conformally flat and Einstein, and hence again hyperbolic; 
the proof then proceeds just as before.

  In the following, we write
\begin{equation}\label{e4.30}
\phi \sim 1 ,
\end{equation}
if there is a constant $C < \infty$ such that $C^{-1} \leq \phi \leq C$. The data 
$\phi$ will be determined by $(M, g)$ as in Theorem 4.4, but $C$ is required to be 
independent of $(M, g)$. The next result, which is of independent interest, relates 
the curvature radius $\rho$ with the stress-energy term $\tau_{(n)}$ in (2.12). 

\begin{lemma}\label{l4.8}
For $(M, g)$ as in Theorem $4.4$, and for any $x \in \partial M$ with $\tau_{0}(x)$ 
sufficiently small, one has the estimate
\begin{equation}\label{e4.31}
\rho^{n}(x)\cdot \sup_{y\in B_{x}(\frac{1}{2}\rho(x))}|\tau_{(n)}|(y) \sim 1.
\end{equation}
The same estimate holds for $g_{(n)}$ in place of $\tau_{(n)}$. 
\end{lemma}

\noindent
{\bf Proof:} Observe that the product in \eqref{e4.31} is scale invariant, cf.~also (2.14). 
Suppose first that there exists a sequence of metrics $(M_{i}, g_{i})$ and points $x_{i}
\in \partial M_{i}$, such that 
\begin{equation}\label{e4.32}
\rho^{n}(x_{i})\cdot \sup_{y\in B_{x_{i}}(\frac{1}{2}\rho(x_{i}))}|\tau_{(n,i)}|(y) 
\rightarrow 0.
\end{equation}
We work in the scale $\hat g_{i} = \rho(x_{i})^{-2}\bar g_{i}$ where 
\begin{equation}\label{e4.33}
\hat \rho(x_{i}) = 1,
\end{equation}
so that $\hat \tau_{0}(x_{i}) \sim 1$. The estimate \eqref{e4.32} implies 
$|\hat \tau_{(n,i)}|(y) << 1$, for all $y \in \hat B_{x_{i}}(\frac{1}{2})$. By Proposition 
4.1, (or the strong control property), the rescaled metrics $\hat g_{i}$ converge 
(in a subsequence) in the $C^{n,\alpha'}$ topology to a limit metric $\hat g$ on a 
maximal connected domain $(U, x_{\infty})$ where $\hat \rho$ does not converge to 
0, (i.e.~the region where the curvature of $\hat g_{i}$ does not blow-up in $L^{p}$). 
Of course $B_{x_{\infty}}(1) \subset U$. As noted prior to Lemma 4.7, the blow-up 
limit of the metric $\gamma$ on $\partial M$ is the flat metric on ${\mathbb R}^{n}$. 

  Now the $C^{n,\alpha'}$ convergence implies that
\begin{equation}\label{e4.34}
\tau_{(n)} = 0,
\end{equation}
on $\partial U$. Since the boundary metric is flat, this implies that $g_{(n)} = 0$, 
cf.~[18], [33]. The unique continuation property, Proposition 2.1, then implies that 
$(U, \hat g)$ is flat, which however contradicts Lemma 4.7. 

  The proof of the opposite inequality in \eqref{e4.31} is similar. Thus, suppose there exist 
$(M_{i}, g_{i})$ and points $x_{i}\in \partial M$ such that  
$$\rho^{n}(x_{i})\cdot \sup_{y\in B_{x_{i}}(\frac{1}{2}\rho(x_{i}))}|\tau_{(n,i)}|(y) > > 1.$$
Now work in the scale where $\sup_{y\in B_{x_{i}}(\frac{1}{2}\rho(x_{i}))}|\tau_{(n,i)}|(y) 
= 1$, and at points $y_{i}$ realizing this supremum. In this scale, (using the fact 
that $\rho$ is Lipschitz, with Lipschitz constant 1), $\rho(y_{i}) >> 1$. Hence, via 
\eqref{e4.21}, $\tau_{0}(y_{i}) >> 1$ in this scale. Then the same arguments as in 
(4.14)-(4.19) together with the strong control property imply that the metric is 
$C^{n,\alpha'}$ close to the flat metric, in bounded domains about $y_{i}$. It follows 
that $\tau_{(n)}$ and $g_{(n)}$ are (arbitrarily) small in bounded domains about 
$y_{i}$ for $i$ sufficiently large, which gives a contradiction to the scale 
normalization above. 

{\endproof}

  Note that the strong control property implies further that the product in \eqref{e4.31} 
is bounded away from 0 and $\infty$ in $C^{\alpha'}$. 

  Now consider a conformal rescaling of the metric $\bar g$ in place of the constant 
rescalings used before. Thus, set 
\begin{equation} \label{e4.35}
\widetilde g = \rho(x)^{-2}\bar g.
\end{equation}
Apriori, the function $\rho$ may not be smooth in $x$, and so \eqref{e4.35} should be 
replaced by the expression $\widetilde g = \rho_{s}(x)^{-2}\bar g$, where $\rho_{s}$ 
is a $C^{\infty}$ smoothing of $\rho$ satisfying $\frac{1}{2}\rho \leq \rho_{s} \leq 
\rho$. Thus, in the following, \eqref{e4.35} is understood with respect to $\rho_{s}$ 
in place of $\rho$; however, to keep the notation reasonable, and also because 
this difference is only of minor significance, we continue to work the notation 
$\rho$. 

  If one sets $\hat g = \tau_{0}^{-2}(x)\bar g$, then by \eqref{e4.5} and \eqref{e4.21} 
the two metrics $\widetilde g$ and $\hat g$ are $C$-quasi-isometric, i.e. 
\begin{equation}\label{e4.36}
\hat g \sim \widetilde g,
\end{equation}
in the sense of \eqref{e4.30}. 

  Let $\widetilde \rho$ be the $L^{p}$ curvature radius of 
$(M, \widetilde g)$ and similarly for $\widetilde r_{h}^{m-1,\alpha}$. Also, let 
$\widetilde g_{(n)}$ be the $n^{\rm th}$ term in the Fefferman-Graham expansion 
(2.9)-(2.10) for the geodesic compactification of $(M, g)$ determined by the boundary 
metric $\widetilde \gamma = \widetilde g|_{\partial M}$, and similarly for the 
stress-energy term $\widetilde \tau_{(n)}$.

\begin{lemma}\label{l4.9}
Let $(M, g)$ be as in Theorem $4.4$. Then on $(M, \widetilde g)$, for all $x \in \partial M$, 
one has the estimates
\begin{equation}\label{e4.37}
\widetilde \rho (x) \sim 1, \ \ {\rm and} \ \ 
\widetilde r_{h}^{m-1,\alpha}(x) \sim 1. 
\end{equation}
These estimates also hold for the geodesic compactification determined by $\widetilde \gamma$. 
Further, for $\rho(x)$ sufficiently small,  
\begin{equation}\label{e4.38}
\sup_{y\in B_{x}(\frac{1}{2})}\widetilde \tau_{(n)}(y) \sim 1 .
\end{equation}
\end{lemma}

\noindent
{\bf Proof:} 
We first prove $\widetilde \rho$ is bounded below on $\partial M$. To do this, 
standard formulas for the behavior of the curvature under conformal changes give 
$$\widetilde R = \rho^{-2}[R - g\wedge(D^{2}\log \rho^{-1} - (d\log \rho)^{2} + 
{\tfrac{1}{2}}|d\log \rho|^{2}g)],$$
where $\wedge$ denotes the Kulkarni-Nomizu product, cf.~[10]. Hence, 
\begin{equation}\label{e4.39}
|\widetilde R|_{\widetilde g} \leq c(\rho^{2}|R| + \rho|D^{2}\rho| + |d\rho|^{2}), 
\end{equation}
where the right side is taken with respect to $\bar g$ and we recall the remark 
on smoothing following \eqref{e4.35}. Note also that the right side of 
\eqref{e4.39} is scale invariant. Thus, it suffices to obtain an $L^{p}$ or $L^{\infty}$ bound 
of the terms on the right in \eqref{e4.39}, in the scale where $\rho = 1$, on balls of radius 
$r_{0}$, for a uniform $r_{0} > 0$. Since the function $\rho$ is Lipschitz, with 
Lipschitz constant 1, the last term is bounded in $L^{\infty}$. By the strong 
control property, the first two terms are also bounded in $L^{\infty}$ on the ball of 
radius $\frac{1}{2}$. This gives a uniform upper bound on $|\widetilde R|_{\widetilde g}$ 
on the ball of radius $\frac{1}{2}$ and hence
$$\widetilde \rho \geq \rho_{0},$$
for a fixed $\rho_{0}$, (depending only on the constant $c_{0}$ in (3.4)). Also, since 
$(M, \bar g)$ satisfies the strong control property, in the scale where $\rho(x) = 1$ one 
has a uniform lower bound on $r_{h}^{m-1,\alpha}$, $r_{h}^{m-1,\alpha}(x) \geq r_{0}$, 
which again via \eqref{e4.39} gives a lower bound on $\widetilde r_{h}^{m-1,\alpha}$ on 
$\partial M$. 

  The opposite estimate $\widetilde \rho (x) \leq \rho_{1}\rho (x)$, for a fixed 
$\rho_{1} < \infty$, also holds. For if, on some sequence $(M_{i}, \bar g_{i}, x_{i})$, 
$\widetilde \rho (x_{i}) \rightarrow \infty$ while $\rho (x_{i}) = 1$, then the 
metric $\widetilde g_{i}$ is almost flat in $L^{p}$ in large balls 
$B_{x_{i}}((1-\delta)\widetilde \rho(x_{i}))$, for any fixed $\delta > 0$. It follows 
that any maximal limit $(N, g_{\infty}', x_{\infty})$, as following \eqref{e4.24}, is 
conformally flat. As noted following Lemma 4.7, this contradicts (the proof of) 
Lemma 4.7. Combining the arguments above proves \eqref{e4.37}. 

  Next, the estimate \eqref{e4.37} gives uniform control in $C^{m-1,\alpha}$ of the conformal 
factor relating the geodesic compactification with respect to $\widetilde g$ and the 
rescaling of $\bar g$ by the constant factor $\rho(x)^{-2}$, for any fixed $x$. This 
implies that \eqref{e4.37} also holds for the geodesic compactification with respect to 
$\widetilde g$. In particular, this gives an upper bound on the term 
$\widetilde \tau_{(n)}$ in \eqref{e4.38}. 

  Regarding the lower bound, suppose first $n$ is odd. Then the lower bound on 
$\widetilde \tau_{(n)}$ in \eqref{e4.38} follows immediately from \eqref{e4.31} and the conformal 
transformation rule (2.14). A similar argument holds also for $n$ even. Thus, although 
(2.14) does not hold exactly, it does hold modulo terms which involve lower order derivatives 
of the boundary metric; for $\rho(x)$ sufficiently small, these terms scale at a lower power 
of $\rho(x)^{-1}$ than $g_{(n)}$, so that (2.14) holds to leading order. 

{\endproof}

  We are now in position to begin the proof of Proposition 4.6 per se. Suppose then \eqref{e4.22} 
is false. Then there exist AH Einstein metrics $(M_{i}, g_{i})$ with $C^{m,\alpha}$ 
controlled boundary metrics $\gamma_{i}$ and points $x_{i}$ such that 
$\tau_{0}(x_{i}) \rightarrow 0$, and so 
\begin{equation}\label{e4.40}
\rho_{i}(x_{i}) \rightarrow 0,
\end{equation}
where $\rho_{i}$ is the $L^{p}$ curvature radius of the geodesic compactification 
$\bar g_{i}$ of $g_{i}$ with boundary metric $\gamma_{i}$. From now on, we work in 
the conformally rescaled metric \eqref{e4.35}. The Lipschitz property of 
$\rho$, together with \eqref{e4.40} implies that 
\begin{equation}\label{e4.41}
vol_{\widetilde \gamma_{i}}(\partial M) \rightarrow \infty.
\end{equation}

  The main point now is to use the relation (2.22) in the $\widetilde \gamma_{i}$ metric: 
\begin{equation}\label{e4.42}
\int_{\partial M}\langle {\mathcal L}_{X}\tau_{(n)} + [(1 - {\tfrac{2}{n}})div X] \tau_{(n)}, 
h_{(0)} \rangle dV_{\widetilde \gamma_{i}} = 
\int_{\partial M}\langle \sigma_{(n)} + {\tfrac{1}{2}}tr h_{(0)} \tau_{(n)}, 
\hat{\mathcal L}_{X}\gamma \rangle dV_{\widetilde \gamma_{i}} + b_{(n)}.
\end{equation}
Here we have dropped the $i$ and tilde from the notation; the inner product in 
\eqref{e4.42} as well as $\tau_{(n)}$ and $\sigma_{(n)}$ are taken with respect 
to $\widetilde \gamma_{i}$; also the term $b_{(n)}$ is given by the last term in 
(2.22). Recall from (2.12) that $\tau_{(n)}$ is determined by $g_{(n)}$ and 
the boundary metric. By (2.21) and Lemma 2.2, $\sigma_{(n)}$ is any solution 
of the system
\begin{equation}\label{e4.43}
\delta \sigma_{(n)} = -\delta'(\tau_{(n)}),
\end{equation}
$$tr \sigma_{(n)} = -tr'(\tau_{(n)}) + (a_{(n)})',$$
where $\delta'$, $tr'$ and $(a_{(n)})'$ are the variations of $\delta$, $tr$ and 
$a_{(n)}$ at $(\partial M, \widetilde \gamma_{i})$ in the direction $h_{(0)}$. 
Of course solutions $\sigma_{(n)}$ of \eqref{e4.43} are invariant under the addition 
of transverse-traceless terms, as is \eqref{e4.42}. We point out that the system of 
equations \eqref{e4.43} is conformally invariant, since these equations define formal 
solutions to the conformally compactified linearized Einstein equations. In more detail, 
a conformal change of the boundary metric induces a conformal change of $h_{(0)}$; 
as described in and following (2.14), such conformal changes in turn induce a 
transformation of $\tau_{(n)}$, $a_{(n)}$ and $\sigma_{(n)}$. The space of solutions of 
\eqref{e4.43} is thus transformed into itself under changes in the conformal 
compactification of solutions of the linearized Einstein equations. 

  Consider then the equations \eqref{e4.43} on $(\partial M, \gamma_{i})$, and set 
$\gamma_{i} = \gamma$.  We first claim that the system \eqref{e4.43} is solvable 
for any variation $h_{(0)}$ of $(\partial M, \gamma)$, and hence, by Proposition 2.4, 
\eqref{e4.42} also holds for all such variations $h_{(0)}$. To see this, by 
Lemma 2.2 and Proposition 2.3, we know that each equation in \eqref{e4.43} is 
individually solvable, and thus need to show that the equations are simultaneously 
solvable. As is well-known, any symmetric bilinear form $\sigma$ on 
$(\partial M, \gamma)$ can be written uniquely as $\sigma = \delta^{*}V + f\gamma 
+ k$, where $k$ is transverse-traceless. Setting $-\delta'(\tau_{(n)}) = \phi_{1}$, 
$-tr'(\tau_{(n)}) + a_{(n)}'  = \phi_{2}$, the system \eqref{e4.43} becomes 
\begin{equation}\label{e4.43a}
\delta\delta^{*} V - \frac{1}{n}d\delta V = \phi_{1}+\frac{1}{n}d\phi_{2},
\end{equation}
$$-\delta V + nf = \phi_{2}.$$ 
By Proposition 2.3, the term $\phi_{1} \in Im \delta$, and similarly $d\phi_{2} \in 
Im \delta$, since for any Killing field $X$ on $(\partial M, \gamma)$, 
$\langle d\phi_{2}, X \rangle = \langle \phi_{2}, \delta X \rangle = 0$. Thus, 
the first equation is solvable for $V$, and hence the second equation is solvable 
for $f$, which proves the claim. Note that the first equation in \eqref{e4.43a} is 
an elliptic equation for $V$. 

  Next, we claim there is a solution $\sigma_{(n)} = \delta^{*}V + f\gamma$ of 
\eqref{e4.43} on $(\partial M, \gamma)$, ($\gamma = \gamma_{i}$), such that 
\begin{equation}\label{e4.44}
||\sigma_{(n)}||_{C^{m-n,\alpha}} \leq C ||h_{(0)}||_{C^{m-n,\alpha}}||\tau_{(n)}||_{C^{m-n,\alpha}} ,
\end{equation}
where $C$ depends only on $K$ in Theorem 4.4. This follows directly from standard 
elliptic estimates for $V$ associated to the equation \eqref{e4.43a}, (choosing 
$V$ orthogonal to the kernel of the operator), together with the fact that 
$\gamma_{i}$ are uniformly bounded in $C^{m,\alpha}$. 

  Transforming then to the conformally related metrics $\widetilde \gamma_{i}$ and using 
the fact \eqref{e4.37} that $\tau_{(n)}$ is uniformly controlled in $C^{\alpha}$ then gives 
\begin{equation}\label{e4.45}
||\sigma_{(n)}||_{C^{\alpha}} \leq C||h_{(0)}||_{C^{\alpha}},
\end{equation}
for a fixed constant $C$, independent of $i$, where $\sigma_{(n)}$ and the norms are 
taken with respect to the $\widetilde \gamma_{i}$ metric. 

  Let $Z_{i} = \{x\in \partial M: \rho_{i}(x) \leq \varepsilon_{i}\}$, for some sequence 
$\varepsilon_{i} \rightarrow 0$. For any $x_{i} \in Z_{i}$, by Lemma 4.9 the pointed sequence 
$(\partial M, \widetilde \gamma_{i}, x_{i})$ converges in $C^{n,\alpha'}$, and uniformly on 
compact sets, to a complete conformally flat metric $\widetilde \gamma_{\infty}$ on a domain 
$V \subset \partial M$. The metric $\widetilde \gamma_{\infty}$ is conformally flat because 
one is conformally blowing up the sequence of smoothly controlled metrics $\gamma_{i}$ on 
$\partial M$. Of course, the convergence is modulo diffeomorphisms which blow-up small 
regions in $\partial M$ to unit size. Thus, there exist embeddings $F_{i}: V \rightarrow 
\partial M$ such that $(F_{i})^{*}\widetilde \gamma_{i} \rightarrow \widetilde \gamma_{\infty}$. 
For example, if $x_{i}$ realizes the minimal value of $\rho_{i}$ in \eqref{e4.40}, then 
$\widetilde \gamma_{\infty}$ is a complete conformally flat metric on ${\mathbb R}^{n}$. 
Other choices of base points may lead to complete conformally flat metrics on 
${\mathbb R}\times S^{n-1}$ for example. In general, since the conformal class 
$[\gamma_{i}]$ is controlled in $C^{m,\alpha}$, the limit domain $V$ is contained in 
${\mathbb R}^{n}$, with closure $\bar V = {\mathbb R}^{n}$. Similarly, the forms 
$\tau_{(n,i)}$ and $\sigma_{(n,i)}$ converge in $C^{\alpha'}$ to limit forms 
$\tau_{(n)}$ and $\sigma_{(n)}$ on $V$; recall that the tilde has been dropped from 
the notation. 

  For the same reasons, since the metrics $\widetilde g_{i}$ satisfy \eqref{e4.37}, the pointed 
sequence $(M_{i}, \widetilde g_{i}, x_{i})$ also converges, in $C^{n,\alpha'}$, to a maximal 
connected limit $(N, \widetilde g_{\infty}, x_{\infty})$, with boundary data 
$(V, \widetilde \gamma_{\infty})$. 

  Let $Y$ be any $C^{2}$ smooth vector field on $\partial M$. The vector fields $X_{i} = 
(F_{i})^{*}Y$ are blow-ups of $Y$ based at $x_{i}$ and converge to a conformal Killing 
field $X_{\infty}$ on $(V, \widetilde \gamma_{\infty})$, so that
\begin{equation}\label{e4.46}
\hat {\mathcal L}_{X_{i}}\widetilde \gamma_{i} \rightarrow 
\hat {\mathcal L}_{X_{\infty}}\widetilde \gamma_{\infty} = 0 .
\end{equation}
Similarly, since the metrics $\gamma_{i}$ are bounded in $C^{n,\alpha}$, the terms $a_{(n)}$ 
are also bounded, and hence, wirh respect to the boundary metrics $\widetilde \gamma_{i}$, 
$(a_{(n)})_{i} \rightarrow 0$ pointwise. Combining then \eqref{e4.42} with \eqref{e4.41} 
and \eqref{e4.45}, it thus follows that for $i$ large,
\begin{equation}\label{e4.47}
\sup_{|Y| \leq 1}\sup_{|h_{(0)}| \leq 1}|\int_{\partial M}\langle 
\widetilde{\mathcal L}_{X_{i}}\tau_{(n,i)}, h_{(0)} \rangle dV_{\widetilde \gamma_{i}}| << 
vol_{\widetilde \gamma_{i}}(\partial M),
\end{equation}
where $\widetilde{\mathcal L}_{X}\tau_{(n)} = {\mathcal L}_{X}\tau_{(n)} + 
[(1 - {\tfrac{2}{n}})div X] \tau_{(n)}$, $|h_{(0)}| = ||h_{(0)}||_{C^{\alpha}(\widetilde 
\gamma_{i})}$, and $|Y| = ||Y||_{C^{2}(\gamma_{i})} \sim ||Y||_{C^{2}(\gamma_{\infty})}$. 
Since $h_{(0)}$ is arbitrary, it follows from \eqref{e4.47} that there exist (possibly many) 
regions where $\rho_{i} \rightarrow 0$, for which, on a complete, conformally 
flat limit $(V, \widetilde \gamma_{\infty}, x_{\infty})$ of the boundary, one has
\begin{equation}\label{e4.48}
\widetilde{\mathcal L}_{X_{\infty}}\tau_{(n)} = 0,
\end{equation}
for any conformal Killing field $X_{\infty}$ on $(V, \widetilde \gamma_{\infty})$. 
We claim that the only solution of \eqref{e4.48}, (for all conformal Killing $X_{\infty}$), is 
\begin{equation}\label{e4.49}
\tau_{(n)} = 0.
\end{equation}
To see this, a well-known relation gives ${\mathcal L}_{X}\tau_{(n)} = \nabla_{X}\tau_{(n)} + 
2\delta^{*}X\circ \tau_{(n)}$, so that if $X = X_{\infty}$ is conformal Killing, 
then ${\mathcal L}_{X_{\infty}}\tau_{(n)} = \nabla_{X_{\infty}}\tau_{(n)} + \frac{2}{n}
(div X_{\infty})\tau_{(n)}$. Hence \eqref{e4.48} gives 
\begin{equation}\label{e4.50}
\nabla_{X_{\infty}}\tau_{(n)} = -(div X_{\infty})\tau_{(n)},
\end{equation}
for all conformal Killing fields $X_{\infty}$. Pairing \eqref{e4.50} with $\tau_{(n)}$ 
gives 
\begin{equation}\label{e4.51}
\frac{1}{2}X_{\infty}(|\tau_{(n)}|^{2}) = -(div X_{\infty})|\tau_{(n)}|^{2}.
\end{equation}
Hence, if $\tau_{(n)} \neq 0$, then $X_{\infty}(\log |\tau_{(n)}|) = -div X_{\infty}$, 
so that by \eqref{e4.50}, 
$$\nabla_{X_{\infty}}(|\tau_{(n)}|^{-1}\tau_{(n)}) = 0.$$
Since $(V, \widetilde \gamma_{\infty})$ is conformally flat and $X_{\infty}$ is an 
arbitrary conformal Killing field, this easily implies that $|\tau_{(n)}|$ is constant, 
which from \eqref{e4.51} implies \eqref{e4.49}. 

  Since the convergence of $(M_{i}, \widetilde g_{i}, x_{i})$ to $(N, \widetilde g_{\infty}, 
x_{\infty})$ is $C^{n,\alpha'}$, \eqref{e4.49} contradicts \eqref{e4.37}, which proves 
Proposition 4.6. (Equivalently, \eqref{e4.49} can be shown to imply that the the blow-up 
limit $(N, g', x_{\infty})$ of $(M_{i}, \bar g_{i}, x_{i})$ as preceding Lemma 4.7 is flat, 
contradicting Lemma 4.7). Propositions 4.5 and 4.6 together imply Theorem 4.4, which also 
completes the proof of this result.  

{\endproof}

  Combining Proposition 4.1 and Theorem 4.4 gives the following main result of 
this section. 
\begin{corollary}\label{c4.10}

  On a fixed $4$-manifold $M$, suppose $g \in E_{AH}^{m,\alpha}$, $m \geq 4$, has 
boundary metric $\gamma$ satisfying $||\gamma||_{C^{m,\alpha}} \leq K$ in 
a fixed coordinate atlas for $\partial M$. Then there exists 
$\delta_{1} = \delta_{1}(K) > 0$ such that the $C^{m-1,\alpha}$ geometry 
of the geodesic compactifications $\bar g$ is uniformly controlled in 
$U_{\delta_{1}}$, in that \eqref{e4.2} holds, with $r_{0}$ depending only on $K$. 

  The same result holds in all dimensions in which the strong control property 
holds. 
\end{corollary}

{\endproof}

\section{Properness of the Boundary Map.}
\setcounter{equation}{0}

  In this section, we combine the results of \S 3 and \S 4 to prove Theorem A, 
i.e.~the properness of $\Pi$. Following this, in \S 5.2, we also analyse the 
behavior when $\Pi$ is not proper and prove this is solely due to the formation 
of Einstein cusp metrics, cf.~Theorem 5.4 and the discussion following it. 

\medskip

 {\bf 5.1.} Corollary 4.10 gives uniform $C^{m-1,\alpha}$ control of the geodesic 
compactification of an AH Einstein metric $(M, g)$ in a neighborhood of definite 
size about $\partial M$, in terms of $C^{m,\alpha}$ control of the boundary metric 
$\gamma$. Theorem 3.7 addresses control in the interior, away from the boundary. 
Regarding the hypotheses (3.17)-(3.21) of Theorem 3.7, the assumption (3.19) on a 
lower bound on the inradius of $\bar g$ and the assumption (3.20) on an upper bound 
on the diameter $T$ of $S(t_{1})$ are now immediate consequences of Corollary 4.10. 
Next we show that there is a global $L^{2}$ bound on $W$, i.e.~(3.21) holds. 

\begin{proposition} \label{p 5.1.}
  Let $(M, g)$ be an AH Einstein metric on a $4$-manifold $M$, with boundary metric 
$\gamma$. Then there is a constant $\Lambda$, depending only on the topology of $M$ 
and the $C^{4,\alpha}$ norm of $\gamma$, such that 
\begin{equation} \label{e5.1}
\int_{M}|W_{g}|^{2}dV_{g} \leq  \Lambda . 
\end{equation}
\end{proposition}

\noindent
{\bf Proof:}
 The integral (5.1) is conformally invariant. In the tubular neighborhood $U = 
U_{\mu_{1}}$ about $\partial M$, (cf.~\eqref{e4.22}), we compute the $L^{2}$ integral 
(5.1) with respect to the compactification $\bar g$, while in $M \setminus U$, 
the integral is computed with respect to the Einstein metric $g$. 

  Thus, by Corollary 4.10, the curvature $\bar R$ of $\bar g$ is uniformly 
bounded in $U$. Since $vol_{\gamma}\partial M$ is also uniformly bounded, 
the curvature bound also gives a uniform upper bound on the volume of 
$(U, \bar g)$. Hence, by conformal invariance, the integral over $U$ in 
(5.1) is uniformly bounded. 

  For the integral over $M \setminus U$, we use the Chern-Gauss-Bonnet 
theorem for manifolds with boundary, as in [5]. Let $S(t_{0}) = 
\partial U_{\mu_{0}}$, $\mu_{0} = \mu_{1}/2$, viewed as a boundary 
in $(M, g)$, and let $\Omega  = M \setminus U_{\mu_{0}}$. Since $(\Omega, g)$ 
is Einstein, one has,
$$\frac{1}{8\pi^{2}}\int_{\Omega}|W|^{2}dV 
=\frac{1}{8\pi^{2}}\int_{\Omega}(|R|^{2}- 6) = \chi(\Omega)-\frac{3}{4\pi^{2}}vol\Omega  
+ \int_{S(t_{1})}B(R,A) \leq  \chi (M) + \int_{S(t_{1})}B(R,A);$$
compare also with (3.16). Here $B(R,A)$ is a boundary term, depending 
on the curvature $R$ and second fundamental form $A$ of $S(t_{0})$ in 
$(M, g)$. But again by Corollary 4.10, $R$ and $A$ are uniformly 
controlled on $S(t_{0}) \subset  (M, g)$, as is $vol_{g}S(t_{0})$. 
Hence, the $L^{2}$ norm of $W$ over $\Omega$ is uniformly bounded 
and so (5.1) follows.

{\endproof}

 We now assemble the work above to obtain the following result, which represents 
the major part of Theorem A. Recall from [6] that if $(M, g)$ is an AH Einstein 
metric and $\bar g$ is a geodesic compactification, then its width $Wid_{\bar g}M$ 
is defined by
\begin{equation} \label{e5.2}
Wid_{\bar g}M = \sup\{t(x): x\in M\}. 
\end{equation}
The width depends on a choice of the boundary metric $\gamma$ for the 
conformal infinity. However, if $\gamma$ and $\gamma'$ are 
representatives in $[\gamma]$, then
$$C^{-1}Wid_{\bar g}M \leq  Wid_{\bar g'}M \leq  C Wid_{\bar g}M, $$
where $C$ depends only on the $C^{1}$ norm of $\gamma^{-1}\gamma' $ and 
$(\gamma')^{-1}\gamma$ in a fixed coordinate system on $\partial M.$
\begin{theorem} \label{t 5.2.}
  Let $\{g_{i}\}$ be a sequence of AH Einstein metrics on $M$, with 
boundary metrics $\gamma_{i}$ and suppose $\gamma_{i} \rightarrow  
\gamma $ in the $C^{m,\alpha'}$ topology on $\partial M$, $m \geq 4$. 
Suppose further that 
\begin{equation} \label{e5.3}
H_{2}(\partial M, {\mathbb R} ) \rightarrow  H_{2}(\bar M, {\mathbb R} ) 
\rightarrow  0, 
\end{equation}
and there is a constant $D < \infty$ such that
\begin{equation} \label{e5.4}
Wid_{\bar g_{i}}M \leq  D. 
\end{equation}

 Then a subsequence of $\{g_{i}\}$ converges smoothly and uniformly on 
compact subsets to an AH Einstein metric $g$ on $M$ with boundary 
metric $\gamma .$ The geodesic compactifications $\bar g_{i}$ converge 
in the $C^{m-1,\alpha'}$ topology to the geodesic compactification 
$\bar g$ of g, within a fixed collar neighborhood $U$ of $\partial M.$
\end{theorem}

\noindent
{\bf Proof:} 
Corollary 4.10 proves the last statement. As noted above, together with 
Proposition 5.1, it follows that the all the hypotheses (3.17)-(3.21) of 
Theorem 3.7 are satisfied. Thus if one chooses base points $x_{i}$ satisfying 
(3.22), then a subsequence of $(M, g_{i}, x_{i})$ converges smoothly and uniformly 
on compact sets, to a limit AH Einstein metric $(N, g, x)$.

 For a fixed $d > $ 0 small, let 
\begin{equation} \label{e5.5}
(M_{i})_{d} = \{x\in (M, g_{i}): t_{i}(x) \leq  d\}, 
\end{equation}
and define $N_{d}$ in the same way. Then the smooth convergence of the 
compactifications $\bar g_{i}$ implies that $(M_{i})_{2d}$ is 
diffeomorphic to $N_{2d}$ and each is a collar neighborhood of 
$\partial M.$ Further, the limit metric $g$ on $N_{2d}$ is AH, with 
boundary metric $\gamma$.

 On the other hand, (5.4) implies that the complementary domains 
\begin{equation} \label{e5.6}
(M_{i})^{d} = \{x\in (M, g_{i}): t_{i}(x) \geq  d\} 
\end{equation}
have uniformly bounded diameter with respect to $g_{i}$, and so Theorem 3.7 
implies that $(M_{i})^{d/2}$ is diffeomorphic to the limit domain $N^{d/2}$ 
in $(N, g)$. It follows that $N = M$ and so $g$ is an AH Einstein metric 
on $M$, with boundary metric $\gamma$. 

{\endproof}
 
 Combining the results above leads easily to the proof of Theorem A. 

\medskip
\noindent
{\bf Proof of Theorem A.}

 Let $\gamma_{i}$ be a sequence of boundary metrics in ${\mathcal C}^{o}$, 
with $\gamma_{i} \rightarrow \gamma \in {\mathcal C}^{o}$ in the 
$C^{m,\alpha'}$ topology on $\partial M$, $m \geq 4$, with $\Pi (g_{i}) = 
[\gamma_{i}]$. Since only the conformal classes are uniquely 
determined, one may choose for instance $\gamma_{i}$ to be metrics of 
constant scalar curvature. To prove $\Pi^{o}$ is proper, one needs to 
show that $\{g_{i}\}$ has a convergent subsequence in ${\mathcal E}_{AH}$ 
to a limit metric $g\in{\mathcal E}_{AH}$ with $\Pi [g] = [\gamma]$.

 Suppose first there is a constant $s_{0} > 0$ such that
$$s_{\gamma_{i}} \geq  s_{0}. $$
where $s_{\gamma_{i}}$ is the (intrinsic) scalar curvature of the 
boundary metric $\gamma_{i}.$ It is then proved in [6, Prop.5.1] that
\begin{equation} \label{e5.7}
Wid_{\bar g_{i}}M \leq  \sqrt{3}  \pi / \sqrt{s_{0}}, 
\end{equation}
cf. also (A.12). Hence, in this case Theorem A follows directly from 
Theorem 5.2.

 Next, suppose only $s_{\gamma_{i}} \geq $ 0. If there is some constant 
$D <  \infty $ such that $Wid_{\bar g_{i}}M \leq  D$, then again 
Theorem 5.2 proves the result. Suppose instead
\begin{equation} \label{e5.8}
Wid_{\bar g_{i}}M \rightarrow  \infty .
\end{equation}
In this case, there is a rigidity result associated with the limiting 
case of (5.7) as $s_{0} \rightarrow 0$, proved in [6, Rmk.5.2, Lem.5.5]. 
Namely, $s_{\gamma_{i}} \geq 0$ and (5.8) imply that the sequence 
$(M, g_{i}, x_{i})$, for $x_{i}$ as in (3.22), converges in 
the Gromov-Hausdorff topology to a hyperbolic cusp metric
\begin{equation} \label{e5.9}
g_{C} = dr^{2}+r^{2}g_{F}, 
\end{equation}
where $g_{F}$ is a flat metric on $\partial M$. It follows from the 
proof of Theorem 5.2, as following (5.5), that $\gamma  = g_{F}$, so that, 
by definition, $\gamma \notin {\mathcal C}^{o}$. This implies that 
necessarily $Wid_{\bar g_{i}} \leq  D$, for some $D < \infty$, which 
completes the proof.

{\endproof}

\begin{remark} \label{r 5.3}
{\rm   Theorem A implies that the space ${\mathcal E}_{AH}^{o}$ on a given 
4-manifold $M$ satisfying (5.3) has only finitely many components whose 
conformal infinities intersect any given compact set in ${\mathcal C}^{o}$. 

  A similar result holds if the filling manifold $M$ is also allowed to 
vary. Thus, there are only finitely many 4-manifolds $M_{i}$ having a common 
boundary $\partial M$, which satisfy (5.3) and $\chi(M_{i}) \leq C$, for some 
constant $C < \infty$, for which 
$$\cap_{i} \Pi({\mathcal E}_{AH}^{o}(M_{i})) \neq \emptyset .$$
The proof of this is exactly the same, using Theorem 3.5 and Remark 3.6 with 
varying manifolds $M_{i}$. }
\end{remark}

{\bf \S 5.2.}
 In this section, we characterize the possible degenerations when $\Pi$ 
is not proper. Observe first that Theorem 5.2 implies that if the 
manifold $M$ satisfies (5.3), then the "enhanced" boundary map
\begin{equation} \label{e5.10}
\Psi : E_{AH} \rightarrow  Met(\partial M)\times {\mathbb R} , 
\end{equation}
$$\Psi (g) = (\Pi (g), Wid_{\Pi (g)}M) $$
is proper. Thus, degenerations of AH Einstein metrics with controlled 
conformal infinity can only occur when the width diverges to $\infty$. 
On the other hand, as indicated in the Introduction, in general the 
boundary map $\Pi$ is not proper.

 Define an AH {\it  Einstein metric with cusps}  $(N, g)$ to be a 
complete Einstein metric $g$ on a 4-manifold $N$ which has two types of 
ends, namely AH ends and cusp ends. A cusp end of $(N, g)$ is an end 
$E$ such that $vol_{g}E < \infty$. Thus, $N$ has a compact (possibly 
disconnected) hypersurface $H$, disconnecting $N$ into two non-compact 
connected components $N = N_{1}\cup N_{2}$ where $(N_{1}, g)$ is an AH 
Einstein metric with boundary $H$ and $(N_{2}, g)$ has finite volume, 
so that each end of $N_{2}$ is a cusp end. A natural choice for $H$ is 
the level set $t^{-1}(1),$ where $t$ is a geodesic defining function 
for the conformally compact boundary $\partial_{AH}N$ of $N$. Then 
$N_{1} = \{x\in N: t(x) \leq  1\}$, $N_{2} = \{x\in N: t(x) \geq  1\}$.

 Note that any cusp end $E$ is not conformally compact. As one diverges 
to infinity in $E$, the metric $g$ is collapsing, in that 
$vol_{g}B_{x}(1) \rightarrow 0$ as $x \rightarrow \infty$ in $E$, 
and hence the injectivity radius satisfies $inj_{g}(x) \rightarrow 0$, 
as $x \rightarrow  \infty$ in $E$, see \S 3.1.

\begin{theorem} \label{t 5.4.}
  Let $M$ be a $4$-manifold satisfying \eqref{e1.4}, and let $g_{i}\in E_{AH}$ 
be AH Einstein metrics with boundary metrics $\gamma_{i}\in 
C^{m,\alpha}$, with $\gamma_{i} \rightarrow \gamma$ in the 
$C^{m,\alpha'}$ topology on $\partial M$, $m \geq 4$. Let $t_{i}$ be the 
geodesic defining function associated with $\gamma_{i}$ and choose base 
points $x_{i}\in H_{i}=t_{i}^{-1}(1)$. 

 Then a subsequence of $\{g_{i}\}$ converges, modulo diffeomorphisms in 
${\mathcal D}_{1}$, either to an AH Einstein metric $g$ on $M$, or to an AH 
Einstein metric with cusps $(N, g, x)$, $x = \lim x_{i}$. The 
convergence is smooth and  uniform on compact subsets of $M$, $N$ 
respectively. In both cases, the conformal infinity is given by 
$(\partial M, [\gamma])$. Further, the manifold $N$ weakly embeds in 
$M$, as in \eqref{e3.23}.
\end{theorem}

\noindent
{\bf Proof:} 
For any $D <  \infty$, let $(M_{i})_{D} = \{x\in (M, g_{i}): t_{i}(x) 
\leq  D\}$, as in (5.5). Theorem 5.2 implies that a subsequence of 
$((M_{i})_{D}, g_{i}, x_{i})$ converges smoothly to a limit AH Einstein 
metric $g$ on a domain $N_{D}$, with $N_{D}$ diffeomorphic to 
$(M_{i})_{D}$, and with conformal infinity of $N_{D}$ given by 
$(\partial M, [\gamma])$. 

 If there is a fixed $D < \infty$ such that (5.4) holds, then the 
result follows from the proof of Theorem 5.2 or Theorem A. Thus, we may 
suppose 
\begin{equation} \label{e5.11}
Wid_{\bar g_{i}}M \rightarrow  \infty . 
\end{equation}
In this case, it follows from [6, Lemma 5.4] that there is a constant 
$V_{0} < \infty$, depending only on $\{\gamma_{i}\}$ and the Euler 
characteristic $\chi (M)$, such that
\begin{equation} \label{e5.12}
vol_{g_{i}}(M_{i})^{1} \leq  V_{0}, 
\end{equation}
where $(M_{i})^{1} = \{x\in (M, g_{i}): t_{i}(x) \geq $ 1\} is the 
complementary domain to $(M_{i})_{1}$. (The estimate (5.12) is a 
straightforward consequence of (2.18), or more precisely the bound $V 
\leq \frac{4\pi^{2}}{3}\chi(M)$, given uniform control of the metrics 
$g_{i}$ and $\bar g_{i}$ on $(M_{i})_{1}$).

 By Theorem 3.5 and Remark 3.6, the pointed manifolds $(M, g_{i}, 
x_{i})$, for $x_{i}$ base points in $S_{i}(1) = t_{i}^{-1}(1)$, 
converge in a subsequence and in the Gromov-Hausdorff topology to a 
complete Einstein manifold $(N, g, x)$, $x = \lim x_{i}$. The 
convergence is also in the $C^{\infty}$ topology, uniform on compact 
sets. The domain $N^{1} = \{x\in (N, g): t(x) \geq  1\}$, with $t(x) = 
\lim t_{i}(x)$, is the limit of the domains $(M_{i})^{1}$. 

 The bound (5.12) implies that $N^{1}$ is of finite volume while (5.11) 
implies that $N^{1}$ is non-compact. It follows that the full limit $N 
= N^{1}\cup N_{1}$ with limit metric $g$ is a complete AH Einstein 
manifold with conformal infinity $(\partial M, [\gamma])$ and with a 
non-empty collection of cusp ends. The fact that $N$ weakly embeds in 
$M$ follows exactly as in the proof of Theorem 3.7.

{\endproof} 

 The fact that the AH cusp metric $(N, g)$ has conformal infinity 
$(\partial M, [\gamma])$ and that it weakly embeds in $M$ implies that 
there is a sequence $t_{j} \rightarrow \infty$ such that $N_{t_{j}} 
\subset  N$ embeds in $M$, for $N_{t_{j}}$ as above. Hence, for any $j$ 
large, the manifold $M$ may be decomposed as
\begin{equation} \label{e5.13}
M = N_{t_{j}} \cup  (M\setminus N_{t_{j}}). 
\end{equation}
With respect to a suitable diagonal subsequence $j = j_{i}$, the 
metrics $g_{i}$ on $M$ push the region $M\setminus N_{t_{j_{i}}}$ off 
to infinity as $i \rightarrow \infty$ and $t_{j_{i}} \rightarrow  
\infty$, giving rise to cusp ends in the limit $(N, g)$.

\medskip

 Theorem 5.4 suggests the construction of a natural completion $\bar 
{\mathcal E}_{AH}$ of ${\mathcal E}_{AH}$. Thus, for any $R < \infty$ large, 
let ${\mathcal C} (R)$ denote the space of conformal classes $[\gamma]$ on 
$\partial M$ which contain a representative $\gamma$ satisfying 
$||\gamma||_{C^{4,\alpha}} \leq  R$, with respect to some fixed coordinate 
system for $\partial M$. Let ${\mathcal E}_{AH}(R) = \Pi^{-1}(R)$ and for 
any $g\in{\mathcal E}_{AH}(R)$, choose a base point $x\in t^{-1}(1) \equiv  
H_{g}$, where $t$ is the geodesic defining function associated to the 
boundary metric $\gamma$.

 Then Theorem 5.4 implies that the completion $\bar{\mathcal E}_{AH}(R)$ of 
${\mathcal E}_{AH}(R)$ in the pointed Gromov-Hausdorff topology based at 
points $x\in H_{g}$ is the set of AH Einstein metrics on $M$, together 
with AH Einstein metrics with cusps $(N, g)$. If $\gamma_{i}$ are 
metrics in ${\mathcal C}(R)$, then a subsequence of $\gamma_{i}$ converges 
in $C^{m,\alpha'}$ to a limit metric $\gamma\in{\mathcal C}(R)$. If 
$g_{i}\in{\mathcal E}_{AH}(R)$ satisfy $\Pi (g_{i}) = \gamma_{i}$, then the 
corresponding subsequence of $g_{i}$ converges in the Gromov-Hausdorff 
topology based at $x_{i}$ to $(N, g)\in\bar {\mathcal E}_{AH}(R)$, and the 
conformal infinity of $(N, g)$ is $\gamma  = \lim \gamma_{i}$. The 
convergence is also in the $C^{\infty}$ topology, uniform on compact 
sets. Of course one may have $N = M$, in which case the limit is an AH 
Einstein metric on $M$.

 For $R < R'$, ${\mathcal E}_{AH}(R) \subset {\mathcal E}_{AH}(R')$ and so we 
may form the union $\bar{\mathcal E}_{AH} = \cup_{R}\bar{\mathcal E}_{AH}(R)$ 
with the induced topology. It is clear that in this topology, the 
boundary map $\Pi$ extends to a continuous map
\begin{equation} \label{e5.14}
\bar {\Pi}: \bar {\mathcal E}_{AH} \rightarrow {\mathcal C} . 
\end{equation}
Further, the following corollary is an immediate consequence of Theorem 
5.4.

\begin{corollary} \label{c 5.5.}
  Let $M$ be a $4$-manifold satisfying \eqref{e1.4}. Then the extended map 
$\bar \Pi: \bar {\mathcal E}_{AH} \rightarrow  {\mathcal C}$ is proper.
\end{corollary}
{\endproof}

 It follows in particular that the image 
$$\bar \Pi(\bar {\mathcal E}_{AH}) \subset{\mathcal C}$$
is a closed subset of ${\mathcal C}$. However, it is not known if $\bar {\mathcal E}_{AH}$ 
is a Banach manifold. Even if $\bar {\mathcal E}_{AH}$ is not a Banach manifold, it would 
be interesting to understand the 'size' of the set of AH Einstein metrics with cusps. 
The set of such metrics corresponds of course exactly to 
$$\partial\bar {\mathcal E}_{AH} = \bar {\mathcal E}_{AH}\setminus{\mathcal E}_{AH}.$$ 
In particular, with regard to the work to follow in \S 6-\S 7, one would like to 
know if the image $\bar \Pi(\partial{\mathcal E}_{AH})$ disconnects ${\mathcal C}$ or not, 
cf. also Remark 7.8 below. Again the fact that $\bar \Pi$ is proper implies that 
$\bar \Pi(\partial{\mathcal E}_{AH})$ is a closed subset of ${\mathcal C}$.

\smallskip

 Let
\begin{equation} \label{e5.15}
\hat {\mathcal C} = {\mathcal C} \setminus \bar \Pi(\partial{\mathcal E}_{AH}); 
\end{equation}
this is the space of conformal classes on $\partial M$ which are not 
the boundaries of AH Einstein metrics with cusps associated to the 
manifold $M$. Also let $\hat {\mathcal E}_{AH} = \bar \Pi^{-1}(\hat {\mathcal C})$, 
so that $\hat {\mathcal E}_{AH}$ is the class of AH Einstein metrics 
on $M$ whose conformal infinity is not the conformal infinity of any AH 
Einstein metric with cusps associated to $M$. If one gives $\hat {\mathcal C}$ 
the relative topology, as a subset of ${\mathcal C} ,$ then the 
following result is also an immediate consequence of the results above.
\begin{corollary} \label{c 5.6.}
  Let $M$ be a $4$-manifold satisfying \eqref{e1.4}. Then the map 
\begin{equation} \label{5.16}
\hat \Pi : \hat {\mathcal E}_{AH} \rightarrow  \hat {\mathcal C}
\end{equation}
is proper.
\end{corollary}{\endproof}
\begin{remark} \label{r 5.7.}
  {\rm The first examples where $\Pi$ is not proper, for example $\Pi^{-1}(pt)$ 
is non-compact in ${\mathcal E}_{AH}$, is the sequence of AH Einstein metrics $g_{i}$ 
on ${\mathbb R}^{2}\times T^{n-1}$, converging to the hyperbolic cusp metric on 
${\mathbb R}\times T^{n}$, discussed in Remark 2.6; cf.~[6, Prop.4.4] for the explicit 
construction of $\{g_{i}\}$. 

 However, these metrics lie in distinct components of the moduli space 
${\mathcal E}_{AH}$ on ${\mathbb R}^{2}\times T^{2}$, and so cannot be connected 
by a curve of metrics on ${\mathbb R}^{2}\times T^{2}$. Thus ${\mathcal E}_{AH}$ 
has infinitely many components, and this is the cause of $\Pi$ being non-proper. 
Currently, there are no known examples where $\Pi$ restricted to a component of 
${\mathcal E}_{AH}$ is not proper. On the other hand, if one passes to the quotient of 
${\mathcal E}_{AH}$ by the full diffeomorphism group, then the metrics $g_{i}$ above 
may be joined by a curve. 

  More recently, G. Craig [17] has constructed infinite sequences of AH Einstein 
metrics $(M_{i}, g_{i})$ on non-diffeomorphic manifolds $M_{i}$, which all have 
the same conformal infinity $(\partial M, [\gamma_{0}])$. These limit on a complete 
hyperbolic manifold $(N, g_{-1})$, with conformal infinity given by 
$(\partial M, [\gamma_{0}])$, but with additional cusp, i.e.~parabolic, ends 
analogous to the hyperbolic cusp end above. In fact, the Einstein metrics are 
constructed by Dehn filling each cusp end of $N$. }
\end{remark}

\section{Degree of the Boundary Map.}
\setcounter{equation}{0}
 
 In this section, we prove that the boundary map $\Pi^{o}$ has a 
well-defined ${\mathbb Z}$-valued degree, following Smale 
[34] and White [37], cf.~also [34]. Throughout this section, we work 
componentwise on ${\mathcal E}_{AH}^{o}$ and ${\mathcal C}^{o}$, but we will 
not distinguish components with extra notation. Thus, we assume 
$$\Pi^{o}: {\mathcal E}_{AH}^{o} \rightarrow  {\mathcal C}^{o}$$ 
where ${\mathcal E}_{AH}^{o}$ and ${\mathcal C}^{o}$ are connected, that is connected 
components of the full spaces. The results of this section also hold 
for the restricted boundary map $\hat \Pi: \hat {\mathcal E}_{AH} 
\rightarrow \hat {\mathcal C}$ in (5.16), but since there is no intrinsic 
characterization of $\hat {\mathcal C}$, we work with ${\mathcal C}^{o}$; see 
also Remark 7.8.

 Since $\Pi^{o}$ is a proper Fredholm map of index 0, the Sard-Smale 
theorem [34] implies that the regular values of $\Pi^{o}$ are open and 
dense in ${\mathcal C}^{o}$. For $\gamma $ a regular value, the fiber 
$(\Pi^{o})^{-1}(\gamma)$ thus consists of a finite number of points, 
i.e.~(equivalence classes of) AH Einstein metrics on $M$. By [34]
\begin{equation} \label{e6.1}
deg_{2}\Pi^{o} = \#(\Pi^{o})^{-1}([\gamma ]) \ \ {\rm mod \ 2}, 
\end{equation}
is well-defined, for any regular value $[\gamma ]$ in ${\mathcal C}^{o}$. 
We recall that if $[\gamma] \notin {\rm Im}\Pi^{o}$, then $[\gamma]$ 
is tautologically a regular value of $\Pi^{o}$. 

 Next, we show that $\Pi^{o}$ has a well-defined degree in ${\mathbb Z}$, 
essentially following [37]. Thus, for $[\gamma]\in {\rm Im}\Pi^{o}$, 
let $g$ be any AH Einstein metric on $M$, with $\Pi^{o}(g) = [\gamma]$. 
Consider the linearization of the Einstein equations, i.e.~as in 
(2.4), the elliptic operator
$$L = \tfrac{1}{2}D^{*}D -  R, $$
acting on $L^{2,2}(M, g)$. The operator $L$ is bounded below on 
$L^{2,2}$ and as in the theory of geodesics or minimal surfaces, let
\begin{equation} \label{e6.2}
ind_{g} \in  {\mathbb Z} , 
\end{equation}
be the $L^{2}$ index of the operator $L$ at $(M, g)$, i.e.~the maximal 
dimension of the subspace of $L^{2,2}(M, g)$ on which $L$ is a negative 
definite bilinear form, with respect to the $L^{2}$ inner product. The 
nullity of $(M, g)$ is the dimension of the $L^{2}$ kernel $K$.

 The main result of this section is the following:

\begin{theorem} \label{t 6.1.}
   Let $\gamma $ be a regular value of $\Pi^{o}$ on $(M, g)$ and define
\begin{equation} \label{e6.3}
deg \Pi^{o} = \sum_{g_{i}\in (\Pi^{o})^{-1}([\gamma ])}(-1)^{ind_{g_{i}}}. 
\end{equation}
Then deg $\Pi^{o}$ is well-defined, i.e. independent of the choice of 
$[\gamma]$ among regular values of $\Pi^{o}$. 
\end{theorem}

\noindent
{\bf Proof:}
 In [37], White presents general results guaranteeing the existence of 
a ${\mathbb Z}$-valued degree, and we will show that the current situation 
is covered by these results. Thus, we refer to [37] for some further 
details.

 Let $[\gamma_{1}]$ and $[\gamma_{2}]$ be regular values of $\Pi^{o}$ 
and let $[\bar \sigma(t)], t\in [0,1]$ be an oriented curve in ${\mathcal C}^{o}$ 
joining them. Choose representatives $\gamma_{1} \in  [\gamma_{1}]$ 
and $\gamma_{2} \in [\gamma_{2}]$ and let $\bar \sigma(t)$ be a curve 
in $Met(\partial M)$ joining $\gamma_{1}$ to $\gamma_{2}$. By [34], we 
may assume that $\bar \sigma$ is transverse to $\Pi$, so that the lift 
$\sigma  = \Pi^{-1}(\bar \sigma)$ is a collection of curves in $E_{AH}$, 
with boundary in the fibers over $\gamma_{1}$ and $\gamma_{2}$. Define an 
orientation on $\sigma$ by declaring that $\Pi$ is orientation preserving 
at any regular point of $\sigma$ which has even index, while $\Pi$ is 
orientation reversing at regular points of $\sigma$ of odd index. Thus, 
provided this orientation is well-defined, the map $\Pi|_{\sigma}: \sigma  
\rightarrow  \bar \sigma$ has a well-defined mapping degree, as a map 
of 1-manifolds. By construction, this 1-dimensional degree is given by 
(6.3) at any regular point of $\sigma$ and hence it follows that (6.3) 
is well-defined. Thus, it suffices to prove that the orientation 
constructed above is well-defined.

 If $\Pi_{*}(d\sigma /dt) \neq  0$ for all $t$, so that all points of 
$\sigma$ are regular, then the index of $\sigma (t)$ is constant, and 
so there is nothing more to prove. Suppose instead that 
$\Pi_{*}(\sigma'(t_{0})) = 0$, so that $\sigma (t_{0})$ is a critical 
point of $\Pi$; (without loss of generality, from here on assume $\sigma$ 
is connected). Hence $\sigma'(t_{0}) = \kappa_{0}\in K$, and for $t$ near 
$t_{0}$, $\sigma(t) = \sigma(t_{0}) + (t-t_{0})\kappa_{0} + O((t-t_{0})^{2})$. 
Without loss of generality, (cf.~[37]), one may assume that $K =  \langle 
\kappa_{0} \rangle$ is the span of $\kappa_{0}$. Thus, $\sigma(t)$ may be 
viewed as a graph over $K_{0}$.

  Let $S_{AH}^{m,\alpha}$ be the space of complete metrics $g$ on $M$ which 
have a $C^{m,\alpha}$ conformal compactification and which satisfy the slice 
condition $\beta_{\sigma({t_{0}})}(g) = 0$. It is proved in [7, Thm.4.1] that the 
map
\begin{equation} \label{e6.4}
S_{AH}^{m,\alpha} \rightarrow {\mathbb S}^{m-2,\alpha}, \ \ g \rightarrow 
Ric_{g} + 3g,
\end{equation}
is a submersion at any $g\in E_{AH}^{m,\alpha}$. Since $Ric_{g} + 3g = 0$ on the 
curve $\sigma(t)$, the implicit function theorem implies that there is a 
2-parameter family $\sigma(t,s)\in S_{AH}^{m,\alpha}$, with $\Pi(\sigma(t,s)) = 
\Pi(\sigma(t)) = \gamma(t)$, for $t$ near $t_{0}$ and 
$s$ near $0$, such that
\begin{equation} \label{e6.5}
Ric_{\sigma(t,s)} + 3\sigma(t,s) \in K.
\end{equation}
It is well-known that forms $\kappa \in K$ satisfy $|\kappa| = O(t^{n})$, 
i.e.~$O(t^{3})$ in the case at hand, cf.~[11], [30]. This, together with the 
fact that such $\kappa$ are transverse-traceless, (cf.~(2.6)), implies that the 
metrics $\sigma(t,s)$ are Einstein to order 3 at conformal infinity $\partial M$, 
in that $\frac{d\sigma}{ds} \in \widetilde T_{\sigma(t,s)}E_{AH}$, cf.~(2.20). 

 By Lemma 2.2, the renormalized action $I_{EH}^{ren}$ is then well-defined on 
the family $\sigma(t,s)$. When viewed as a function of $t,s$, one has
\begin{equation} \label{e6.6}
\frac{d}{ds}I_{EH}^{ren}(\sigma(t,s)) = \int_{M}\langle Ric_{\sigma(t,s)} + 
3\sigma(t,s), \frac{d\sigma}{ds} \rangle dV, 
\end{equation}
since $d\sigma/ds$ has 0 boundary values on $\partial M$. (The 
Euler-Lagrange equations for $I_{EH}^{ren}$ are the same as those 
for the usual Einstein-Hilbert action). By (2.6) and (6.5), 
$$\frac{d}{ds}(Ric_{\sigma(t,s)} + 3\sigma(t,s)) = L(\frac{d\sigma(t,s)}{ds}) \in 
K,$$
which implies that $\frac{d}{ds}\sigma(t,s) \in K$. From this, it follows that 
the integral in (6.6) vanishes only on the curve $\sigma(t)$; thus the curve 
$\sigma (t)$ in the plane $P$ spanned by $t,s$ is the set where $d I_{EH}^{ren}/ds 
= 0$. (The functional $I_{EH}^{ren}$ plays the role of $g$ in [37]). It follows that 
$\sigma (t)$ is (locally) the boundary of the open set $\{d I_{EH}^{ren}/ds > 
0\}$ in $P$:
\begin{equation} \label{e6.7}
\sigma (t) = \partial\{d I_{EH}^{ren}/ds >  0\}. 
\end{equation}
If $ind_{\sigma (t)}$ is even, for $t$ near $t_{0}$, $t \neq t_{0}$, 
give $\sigma (t)$ the boundary orientation induced by this open domain, 
while if $ind_{\sigma (t)}$ is odd for $t$ near $t_{0}$, $t \neq t_{0}$, 
give $\sigma (t)$ the reverse orientation. It is now straightforward to 
see that this definition coincides with the orientation defined at the 
beginning of the proof. Thus, the point $t_{0} \rightarrow \sigma (t_{0})$ 
is a critical point for the map $\pi_{1}\circ\sigma$, where 
$\pi_{1}:(t,s) \rightarrow t$ is projection onto the first factor. If this 
critical point is a folding singularity for $\pi_{1}\circ\sigma$, then the 
index of $\sigma (t)$ changes by 1 in passing through $\sigma (t_{0})$ and 
reverses the orientation of $\lambda (t) = \pi_{1}(\sigma (t))$, (exactly as is 
the case with the standard folding singularity $x \rightarrow  x^{2}$). On 
the other hand, if $\pi_{1}\circ\sigma$ does not fold with respect to $\pi_{1}$, 
(so that one has an inflection point), then the index of $\sigma(t)$ does 
not change through $t_{0}$ and $\pi_{1}$ maps $\sigma (t)$ to $\lambda (t)$ 
in an orientation preserving way. We refer to [37] for further details, and 
also to \S 7 for concrete examples of such folding behavior.

{\endproof}

\begin{remark} \label{r 6.2.}
  {\rm As discussed in \S 2, the inclusion ${\mathcal E}_{AH}^{(m',\alpha')} 
\subset {\mathcal E}_{AH}^{(m,\alpha)}$, is dense, for any $(m',\alpha') \geq 
(m,\alpha)$ and these spaces are diffeomorphic. Hence, the degree of 
$\Pi^{o}$ is independent of $(m, \alpha)$ and defined on the spaces 
${\mathcal E}_{AH}^{(m,\alpha)}$ for all $m \geq 4$. Moreover, the degree is 
also defined on the spaces $\bar {\mathcal E}_{AH}$ and $\hat {\mathcal E}_{AH}$ in 
(5.14) and (5.16), and is independent of $(m, \alpha)$. } 
\end{remark}

\section{Computations of the Degree.}
\setcounter{equation}{0}

 We conclude the paper with computations of $deg\Pi^{o}$ for several 
interesting examples of 4-manifolds. In all cases, the evaluation 
of the degree is made possible by symmetry arguments, using Theorem 2.5.

 Recall that if $deg\Pi^{o} \neq 0$, then $\Pi$ is surjective onto 
${\mathcal C}^{o}$, and so any conformal class $[\gamma]$ in ${\mathcal C}^{o}$ 
on $\partial M$ has a filling by an AH Einstein metric $g$ with 
conformal infinity $[\gamma]$. Of course, $deg\Pi^{o} = 0$ does 
not imply that $\Pi$ cannot be surjective.

 We begin with the proof of Theorem B.

\medskip
\noindent
{\bf Proof of Theorem B.}
 
 As seed metric, choose the hyperbolic (Poincar\'e) metric $g_{0}$ on 
$B^{4}$. This has conformal infinity $[\gamma_{0}]$, where $\gamma_{0}$ 
is the round metric on $S^{3} = \partial B^{4}$. Now the boundary 
metric $\gamma_{0}$ admits a large connected group $SO(4)$ of 
isometries. Theorem 2.5 implies that for any AH Einstein metric $g$ 
on $B^{4}$, $SO(4)$ acts effectively by isometries on $(B^{4}, g)$. 
It is then standard that $g$ must be the Poincar\'e metric on $B^{4}$. 

 It follows that up to isometry, $g_{0}$ is the unique metric with 
conformal infinity $[\gamma_{0}]$. (This also follows from a rigidity 
theorem in [9, Thm.5.2]). Further, it is well-known that the 
$L^{2}$ kernel $K$ of $g_{0}$ is trivial, i.e. $K = \{0\}$; in fact, 
this holds for any Einstein metric of negative sectional curvature, 
cf. [10], [11]. Thus, $g_{0}$ is a regular point of $\Pi$ and since 
$\Pi^{-1}[\gamma_{0}] = g_{0}$, the point $[\gamma_{0}]$ is a regular 
value of $\Pi$. The result then follows from (6.3).

{\endproof}

 We point out the following immediate consequence of the proof.

\begin{corollary} \label{c 7.1.}
  Let $M$ be any $4$-manifold satisfying $\pi_{1}(M, \partial M) = 0$, 
and \eqref{e1.4}, with $\partial M = S^{3}$ and $M \neq  B^{4}$. 
Then on any component of ${\mathcal E}_{AH}^{o}$, 
\begin{equation} \label{e7.1}
deg\Pi^{o} = 0, 
\end{equation}
and $\Pi$ is not surjective. In fact the conformal class of the round 
metric $\gamma_{0}$ on $S^{3}$ is not in $Im \Pi$. 
\end{corollary}

\noindent
{\bf Proof:}
 Suppose the class $[\gamma_{0}] \in {\rm Im}\Pi$, so that there is 
an AH Einstein metric $g$ on $M$, with boundary metric $\gamma_{0}$. 
The proof of Theorem B above implies that necessarily $g = g_{0}$, 
where $g_{0}$ is the Poincar\'e metric on $B^{4}$. Hence, $M = B^{4}$, 
a contradiction. 

{\endproof}

 Next, we have:
\begin{proposition} \label{p 7.2.}
  Let $M = {\mathbb R}^{2}\times S^{2},$ so that $\partial M = S^{1}\times 
S^{2}.$ Then 
\begin{equation} \label{e7.2}
deg\Pi^{o} = 0, 
\end{equation}
and $\Pi $ is not surjective.
\end{proposition}

\noindent
{\bf Proof:}
 As seed metric(s) in this case, we take the remarkable 1-parameter 
family of AdS Schwarz-schild metrics, discussed in detail in [23], cf. 
also [38]. Thus, on ${\mathbb R}^{2}\times S^{2},$ consider the metric
\begin{equation} \label{e7.3}
g_{m} = F^{-1}dr^{2} + Fd\theta^{2} + r^{2}g_{S^{2}(1)}, 
\end{equation}
where $F = F(r) = 1+r^{2}-\frac{2m}{r}.$ The mass parameter $m > $ 0 
and $r \in  [r_{+}, \infty ),$ where $r_{+}$ is the largest root of the 
equation $F(r_{+}) =$ 0. The locus $\{r_{+} =$ 0\} is a totally 
geodesic round 2-sphere $S^{2},$ of radius $r_{+}.$ Smoothness of the 
metric at $\{r_{+} =$ 0\} requires that the circle parameter $\theta $ 
run over the interval $[0,\beta ),$ where $\beta $ is given by
\begin{equation} \label{e7.4}
\beta  = \frac{4\pi r_{+}}{1+3r_{+}^{2}}. 
\end{equation}
It is easily seen that as $m$ varies from $0$ to $\infty$, $r_{+}$ 
varies monotonically from $0$ to $\infty$.

 The metrics $g_{m}$ are isometrically distinct, for distinct values of 
$m$, and form a smooth curve in ${\mathcal E}_{AH}$, with conformal 
infinity given by the conformal class of the product metric $\gamma_{m} 
= S^{1}(\beta )\times S^{2}(1)$. Notice however that the length $\beta$ 
has a maximum value as $m$ ranges over $(0, \infty)$, namely
\begin{equation} \label{e7.5}
\beta  \leq  \beta_{max} = 2\pi / \sqrt{3}, 
\end{equation}
achieved at $r_{+} = 1/ \sqrt{3}$, $m = m_{0} = 2/  \sqrt{3}$. As $m 
\rightarrow 0$ or $m \rightarrow \infty$, one has $\beta \rightarrow 0$.

 Thus, the boundary map $\Pi$ on the curve $g_{m}$ is a fold map, 
folding the ray $m\in (0,\infty )$ onto the $\beta$-interval $(0, 
\beta_{max}]$. The map $\Pi$ restricted to the curve $g_{m}$ is a 2-1 
map, except at the point $g_{m_{0}}$. The metric $g_{m_{0}}$ is a critical 
point of $\Pi$, and, (as will be seen below), the tangent vector 
$(dg_{m}/dm)_{m=m_{0}}$ spans the $L^{2}$ kernel $K_{g_{m_{0}}}$.

 Next, we claim that the metrics $g_{m}$ are the only AH Einstein 
metrics on ${\mathbb R}^{2}\times S^{2}$ with conformal infinity given by 
a product $\gamma_{L} = S^{1}(L)\times S^{2}(1)$. The isometry group of 
$\gamma_{L}$ contains $G = SO(2)\times SU(2)$ and by Theorem 2.5, any AH 
Einstein metric $g_{L}$ on $M$ with boundary metric $\gamma_{L}$ has 
$G$ acting effectively by isometries. As in (2.7), we may choose a geodesic 
defining function for $g_{L}$ so that $g_{L}$ has the form
$$g_{L} = ds^{2} + g_{s}, $$
where $g_{s}$ is a curve of metrics on $S^{1}\times S^{2}$ invariant 
under the $G$-action. Thus, the Einstein metric $g_{L}$ 
has cohomogeneity 1. It then follows from the classification given in 
[31] for instance that the metric $g_{L}$ on ${\mathbb R}^{2}\times S^{2}$ 
is isometric to $g_{m},$ for some $m = m(L)$.

  A similar argument shows that the metrics $g_{m}$, for $m \neq m_{0}$ are 
regular points of $\Pi$. For suppose there exists $\kappa \in K = K_{g_{m}}$. 
Then $\kappa$ is the tangent vector to a curve $\sigma(s) \in E_{AH}$ with 
$\sigma(0) = g_{m}$ and $\Pi_{*}(\kappa) = 0$, so that the boundary metric 
of $\sigma(s)$ is fixed to first order in $s$. The proof of Theorem 2.5 in 
[8] then shows that $\kappa$ inherits the symmetries of $g_{m}$, i.e.~$\kappa$ 
is $G$-invariant. Hence, $\kappa$ is tangent to the curve $g_{m}$, 
so that, without loss of generality, $\sigma(s) = g_{m(s)}$. Since 
$\Pi_{*}(\frac{d}{dm}g_{m}) \neq 0$, when $m \neq m_{0}$, this proves the claim. 

 We may thus compute $deg\Pi^{o}$ by evaluating the formula (6.3) on a 
pair of distinct metrics $g_{m_{1}}$ and $g_{m_{2}}$ with $\Pi 
(g_{m_{1}}) = \Pi (g_{m_{2}})$. From [23, \S 3], one has
\begin{equation} \label{e7.6}
ind_{g_{m_{1}}} = +1, \ \ ind_{g_{m_{2}}} = 0, 
\end{equation}
for $m_{1} < m_{0} < m_{2}$, which proves (7.2). Alternately, since 
$\Pi$ is a 2-1 fold map on $g_{m}$, the proof of Theorem 7.1 shows 
directly that (7.2) holds.

 Further, the symmetry argument above implies that the metrics 
$S^{1}(L)\times S^{2}(1)$ are not in Im$\Pi$, whenever
\begin{equation} \label{e7.7}
L >  \beta_{max}, 
\end{equation}
and hence $\Pi$ is not surjective.
{\endproof}

 This result should be compared with the following:
\begin{proposition} \label{p 7.3.}
  Let $M = S^{1}\times {\mathbb R}^{3},$ so that $\partial M = S^{1}\times S^{2}$. 
Then 
\begin{equation} \label{e7.8}
deg\Pi^{o} = 1, 
\end{equation}
and $\Pi^{o}$ is surjective.
\end{proposition}

\noindent
{\bf Proof:}
 As seed metrics in this situation, we take a family of hyperbolic 
metrics, namely the metrics ${\mathbb H}^{4}(- 1)/{\mathbb Z}$, where the 
${\mathbb Z}$-quotient is obtained by a hyperbolic or loxodromic 
translation of length $L$ along a geodesic in ${\mathbb H}^{4}(- 1)$. The 
conformal infinity is the product metric $S^{1}(L)\times S^{2}(1)$ in 
the case of a hyperbolic translation, and the bent product metric 
$S^{1}(L)\times _{\alpha} S^{2}(1)$ on the same space when the 
translation is loxodromic; the angle $\alpha$ between the factors 
$S^{1}(L)$ and $S^{2}(1)$ corresponds to the loxodromic rotation.

 As in Proposition 7.2, any AH Einstein metric on $S^{1}\times {\mathbb R}^{3}$ 
with boundary metric $S^{1}(L)\times S^{2}(1)$ has $SO(2)\times SO(3)$ 
acting effectively by isometries. Again, the classification in 
[31] implies that, on this manifold, the only such metrics are 
hyperbolic. Hence the result follows as in the proof of Theorem B.
{\endproof}

 Next, we turn to non-trivial disc bundles over $S^{2}$. For the disc 
bundle of degree 1 over $S^{2}$, i.e. $M = {\mathbb C}{\mathbb P}^{2}\setminus B^{4}$, 
with $\partial M = S^{3}$, Corollary 7.1 implies that $\Pi$ is not surjective 
so that
\begin{equation} \label{e7.9}
deg \Pi^{o} = 0. 
\end{equation}

\begin{remark} \label{r 7.4.}
  {\rm Actually, to justify this statement one needs to remove the 
hypothesis (1.4), since in this situation $H_{2}(\partial M)$ does not 
surject onto $H_{2}(\bar M)$. Recall that (1.4) was only used to rule out 
orbifold degenerations in the proof of Theorem 3.7. However, it is 
sometimes possible to rule out orbifold degenerations in the case of 
sufficiently low Euler characteristic directly, without the use of (1.4). 
This situation was treated in [1] in the case of orbifold degenerations on 
compact manifolds, and the argument for AH metrics is similar; thus we refer 
to [1], [3] for some further details. 

 Let $M$ be the disc bundle of degree $k$ over $S^{2}$, so that $\chi (M) = 2$, 
and let $\{g_{i}\}$ be a sequence of AH Einstein metrics on $M$ which converge 
to an AH Einstein orbifold $(X, g)$, with boundary metrics $\gamma_{i}$ 
converging in $C^{m,\alpha}$ to the boundary metric $\gamma$ for $(X, g)$. 
By (2.18), one has 
$$\frac{1}{8\pi^{2}}\int_{M}|W_{g_{i}}|^{2} = 2 -  
\frac{3}{4\pi^{2}}V(g_{i}). $$
The results of Theorem 3.7, Corollary 4.10 and Proposition 5.1 show that 
$V(g_{i}) \rightarrow  V(g)$, where $V(g)$ is the renormalized volume 
of the Einstein orbifold $(X, g)$. Now as discussed following Proposition 3.10, 
orbifold singularities arise by crushing essential 2-cycles in $M$ to points. 
Since $M$ has only one essential 2-cycle, (up to multiplicity), given by the 
zero-section, there can only be one singular point on $X$. For the same 
reason, there is only one Ricci-flat ALE space $(E, g_{\infty})$ associated 
to the singularity. It is then easy to see, cf.~[1, \S 6], [3, \S 3.3] and 
references therein, that 
$$\lim_{i\rightarrow\infty}\int_{M}|W_{g_{i}}|^{2} = 
\int_{X}|W_{g}|^{2} +\int_{E}|W_{g_{\infty}}|^{2}. $$
The formula (2.18) for Einstein orbifolds $(X, g)$ 
gives
$$\frac{1}{8\pi^{2}}\int_{X}|W_{g}|^{2} = \chi (X_{0}) + 
\frac{1}{|\Gamma_{0}|} -  \frac{3}{4\pi^{2}}V(g), $$
where $X_{0}$ is the regular set of $X$ and the orbifold singularity of $X$ 
is $C(S^{3}/\Gamma_{0})$, (cf. [1, (6.2)]). Also 
$$\int_{E}|W_{g_{\infty}}|^{2} = \chi(E) - \frac{1}{|\Gamma|},$$
where $E$ is asymptotic to $C(S^{3}/\Gamma)$ at infinity, (cf. [1, (6.3)]). 
Combining these equations gives $\Gamma_{0} = \Gamma = {\mathbb Z}_{k}$, and hence
$\chi(X_{0}) = 0$, so that $E$ is diffeomorphic to $M$. In particular, this 
rules out orbifold degenerations in the case $k = 1$, where $M = {\mathbb C}{\mathbb P}^{2} 
\setminus B^{4}$, since $\Gamma$ must be non-trivial, $\Gamma \neq \{e\}$. This gives 
(7.9). 

  For general $k$, one thus has
\begin{equation} \label{e7.10}
\int_{E}|W_{g_{\infty}}|^{2} = 2 - \frac{1}{k}.
\end{equation}
On the other hand, the signature formula for ALE metrics, cf. [27] for instance, gives
\begin{equation} \label{e7.11}
\tau(E) = \frac{1}{12\pi^{2}}\int_{E}|W_{+}|^{2} - |W_{-}|^{2} + \eta(S^{3}/{\mathbb Z}_{k}),
\end{equation}
and 
$$\eta(S^{3}/{\mathbb Z}_{k}) = \frac{(k-1)(k-2)}{3k}.$$
Since $\tau(E) = 1$, simple arithmetic as in the Hitchin-Thorpe inequality gives the estimate
\begin{equation} \label{e7.12}
4k-2 \geq |3k - (k-1)(k-2)|.
\end{equation}
The estimate (7.12) implies $k \leq 9$, so that orbifold degenerations are not possible 
if $k \geq 10$. 

  We conjecture that there are in fact no orbifold degenerations for $k \geq 3$. However, 
for $k = 2$, it will be seen below that there are orbifold degenerations. In fact, in 
this case, (7.12) shows that $(E, g_{\infty})$ must be self-dual, and hence by Kronheimer's 
classification [29], $(E, g_{\infty})$ is the Eguchi-Hanson metric. }
\end{remark}

\begin{example} \label{e 7.5}
 {\rm It is interesting to compare the result (7.9) with an explicit 
family of AH Einstein metrics on $M = {\mathbb C}{\mathbb P}^{2}\setminus 
B^{4}$, namely the AdS Taub-Bolt family, cf. [24], [31]. This is a 
1-parameter family of metrics given by
\begin{equation} \label{e7.13}
g_{s} = E_{s}\{(r^{2}-1)F^{-1}(r)dr^{2} + (r^{2}-1)^{-1}F(r)(d\tau 
+ \cos \theta d\phi )^{2} + (r^{2}-1)g_{S^{2}(\frac{1}{2})}\},  
\end{equation}
where the parameters $r$, $s$ satisfy $r \geq s$ and $s > 2$, the 
constant $E_{s}$ is given by
\begin{equation} \label{e7.14}
E_{s} = \frac{2}{3}\frac{s-2}{s^{2}-1}, 
\end{equation}
and the function $F(r) = F_{s}(r)$ is
\begin{equation} \label{e7.15}
F_{s}(r) = 
Er^{4}+(4-6E)r^{2}+\{-Es^{3}+(6E-4)s+\frac{1}{s}(3E-4)\}r+(4-3E).  
\end{equation}
The parameter $\tau$ runs over the vertical $S^{1}$, and $\tau\in 
[0,\beta)$, with
\begin{equation} \label{e7.16}
\beta  = 2\pi . 
\end{equation}
The bolt $S^{2}$ $\{r = s\}$ is a round, totally geodesic 2-sphere, of area 
$A_{s} = \frac{2}{3}\pi (s-2)$. These metrics are AH, with conformal 
infinity given by a Berger (or squashed) $S^{3},$ with base $S^{2}(\frac{1}{2})$, 
and Hopf fiber $S^{1} = S^{1}(L)$ of length $L = 2\pi \sqrt{E}$. 

 In analogy to (7.4)-(7.5), notice that $E \rightarrow 0$ as $s 
\rightarrow 2$ or $s \rightarrow \infty$, and has a maximal value 
$E_{max} = (2 - \sqrt {3}) / 3$ at $s = s_{0} = 2 + \sqrt {3}$. 
Note in particular that since 
$$E_{max} < 1,$$
the round metric $S^{3}(1)$ is not in Im $\Pi (g_{s})$ for any $s$.

 We see that on the curve $g_{s},$ $s \in (2, \infty)$, the boundary 
map $\Pi$ has exactly the same 2-1 fold behavior as for the AdS 
Schwarzschild metric. The use of Theorem 2.5 as in Proposition 7.2 
implies that $\Pi$ is not surjective, and in fact the Berger spheres 
with Hopf fiber length $L$, for $L >  2\pi \sqrt{E_{max}}$ are not in Im$\Pi$.}
\end{example}

 In contrast, one has the following behavior on more twisted disc 
bundles over $S^{2}$. 

\begin{proposition} \label{p 7.6.}
  Let $M = M_{k}$ be the disc bundle of degree $k$ over $S^{2}$, 
$k \geq 10$, so that $\partial M = S^{3}/{\mathbb Z}_{k}$. Then
\begin{equation} \label{e7.17}
deg\Pi^{o} = 1, 
\end{equation}
and $\Pi^{o}$ is surjective.
\end{proposition}

\noindent
{\bf Proof:}
 As seed metrics, choose again the AdS-Taub bolt metrics on $M_{k}$, 
cf.~[24], [31]. These have exactly the same form as 
(7.13)-(7.16), except that the parameter $s$ satisfies $s > 1$, 
$E_{s}$ in (7.14) is replaced by $E_{s,k}$ of the form
\begin{equation} \label{e7.18}
E_{s,k} = \frac{2ks-4}{3(s^{2}-1)}, 
\end{equation}
and the period $\beta$ for $\tau$ is given by $\beta = 2\pi /k$. For 
each $k$, these metrics are AH, and the conformal infinity on 
$S^{3}/{\mathbb Z}_{k}$ is given by the Berger metric, with, as before, 
Hopf circle fibers of length $L = 2\pi\sqrt{E_{s,k}}/k$. When $E = 1$, 
conformal infinity is given by the round metric on $S^{3}/{\mathbb Z}_{k}$. 

 Note however from (7.18) that now the function $E_{s,k}$ is a monotone 
decreasing function of $s$, as $s$ increases from 1 to $\infty .$ 
Hence, $\Pi $ is 1-1 on this curve, and in particular, the metric 
$g_{s_{0}}$, $s_{0} = \frac{1}{3}(k+\sqrt{k^{2}-3})$ has conformal infinity 
the constant curvature metric on $S^{3}/{\mathbb Z}_{k}$. The symmetry 
argument using Theorem 2.5 as before, together with [31], implies these 
metrics are the unique metrics with these boundary values. 
The metrics $g_{s}$ are regular points for $\Pi$, and so (7.17) follows, 
since by Remark 7.4, orbifold singularities cannot occur when $k \geq 10$.
{\endproof}

\begin{remark} \label{r7.7} 
{\rm We expect Proposition 7.6 holds for all $k \geq 3$; the proof above holds 
for such $k$, provided there are no orbifold degenerations, cf. Remark 7.4. 

   However, Proposition 7.6 does not hold when $k = 2$. In this case, (7.18) becomes
\begin{equation} \label{e7.19}
E_{s,2} = \frac{4}{3}\frac{1}{s+1}.
\end{equation}
This is of course monotone decreasing in $s \in (1, \infty)$, but it has a 
finite value at $s = 1$ with $E_{1,2} = \frac{2}{3}$. As $s \rightarrow 1$, 
the area of the bolt $S^{2}$ at $\{r = s\}$ tends to 0, and vanishes at $\{r = 1\}$ 
when $s = 1$. Thus, the Taub-Bolt metric is an orbifold singular metric on 
$C({\mathbb R}{\mathbb P}^{3})$ when $s = 1$. Note this shows that the $k = 2$ discussion 
in Remark 7.4 is sharp, in that as $s \rightarrow 1$ on the Taub-Bolt curve, 
the Ricci-flat ALE space associated to the orbifold singularity is necessarily 
the Eguchi-Hanson metric. In particular, only the values $E \in (0,\frac{2}{3})$ are 
achieved on the Taub-Bolt curve $g_{s}$. I am grateful to Michael Singer for 
pointing out this behavior of the Taub-NUT curve. 

  Hence the round metric $\gamma_{0}$ on ${\mathbb R}{\mathbb P}^{3}$ is not 
in $Im \Pi$ on the curve $g_{s}$. Again, Theorem 2.5 and the classification 
in [31] show that $\gamma_{0} \notin \Pi({\mathcal E})$, where ${\mathcal E}$ is the 
moduli space of AH Einstein metrics on the disc bundle of degree 2 over $S^{2}$. }
\end{remark}

\begin{remark} \label{r 7.8.}
  {\rm We close the paper with some observations on whether the full 
boundary map 
$$\Pi : {\mathcal E}_{AH} \rightarrow  {\mathcal C} $$ 
might be surjective, or almost surjective, at least when 
$deg \Pi^{o} \neq 0$.

 First, recall from Corollaries 5.5 and 5.6 that both the extended map 
$\bar \Pi: \bar {\mathcal E}_{AH} \rightarrow {\mathcal C}$ and the 
restricted map $\hat \Pi : \hat {\mathcal E}_{AH} \rightarrow  \hat {\mathcal C}$ 
are proper. In this regard, it would be very interesting to know if the 
set of boundary values of AH Einstein metrics with cusps disconnects 
${\mathcal C} $ or not, that is whether
$$\bar \Pi(\partial\bar {\mathcal E}_{AH}) \subset  {\mathcal C}  $$
disconnects ${\mathcal C} ,$ or whether $\hat {\mathcal C} = {\mathcal C} \setminus 
\bar \Pi(\partial\bar {\mathcal E}_{AH})$ is path connected, (for instance 
if $\bar \Pi(\partial\bar {\mathcal E}_{AH})$ is of codimension at least 
2). If $\hat {\mathcal C}$ is path connected, then the degree theory 
arguments of \S 6 and \S 7 hold without any change and give a 
well-defined degree $deg\Pi$ on each component of ${\mathcal E}_{AH}$. In 
particular, if this holds and $deg\Pi \neq 0$, then $\Pi$ is almost 
surjective, in that $\Pi $ surjects onto $\hat {\mathcal C}$. 

 On the other hand, if $\bar \Pi(\partial\bar {\mathcal E}_{AH})$ 
disconnects ${\mathcal C}$, then $\bar \Pi(\partial\bar {\mathcal E}_{AH})$ 
represents a ``wall'', past which it may not be possible to fill in 
boundary metrics with AH Einstein metrics. This would be the case for 
instance if $\bar {\mathcal E}_{AH}$ is a Banach manifold with boundary 
$\partial\bar {\mathcal E}_{AH}$ and $\bar \Pi$ maps $\partial\bar {\mathcal E}_{AH}$ 
onto a set of codimension 1 in ${\mathcal C}$. 

 The same issue arises if the condition (1.4) does not hold, so that orbifold 
singular metrics might arise, as discussed in Remark 7.7.

 Finally, it would also be interesting to know if there are topological 
obstructions to the possible formation of cusps, as is the case with 
the formation of orbifold singularities as in (1.4). Thus, with respect 
to the decomposition (5.13), one would like to know for instance if 
some of the homology of one of the two factors injects into that of the 
union $M$.}
\end{remark}

\section*{Appendix}
\setcounter{equation}{0}
\begin{appendix}
\setcounter{section}{1}

 In this appendix, we collect several formulas for the curvature of 
conformal compactifications of AH Einstein metrics. 

 Let $\bar g = \rho^{2}\cdot g$, where $\rho$ is a defining function 
with respect to $\bar g$. The curvatures of the metrics $g$ and $\bar g$ are 
related by the following formulas:
\begin{equation} \label{eA.1}
\bar K_{ab} = \frac{K_{ab} + |\bar \nabla \rho|^{2}}{\rho^{2}} -  
\frac{1}{\rho}\{\bar D^{2}\rho (\bar e_{a},\bar e_{a})+\bar D^{2}\rho 
(\bar e_{b},\bar e_{b})\}. 
\end{equation}
\begin{equation} \label{eA.2}
\bar Ric = - (n-1)\frac{\bar D^{2}\rho}{\rho} + (n\rho^{-2}(|\bar \nabla 
\rho|^{2}- 1) -  \frac{\bar \Delta \rho}{\rho})\bar g, 
\end{equation}
\begin{equation} \label{eA.3}
\bar s = - 2n\frac{\bar \Delta \rho}{\rho} + n(n+1)\rho^{-2}(|\bar \nabla 
\rho|^{2} -  1). 
\end{equation}
The equation (A.2) is equivalent to the Einstein equation (1.2). Observe 
that (A.2) implies that if the compactification $\bar g$ is $C^{2}$, then
\begin{equation} \label{eA.4}
|\bar \nabla \rho| \rightarrow  1 \ \ {\rm at} \ \ \partial M, 
\end{equation}
and so by (A.1), $|K_{ab}+1| = O(\rho^{2}).$ For $r = -\log \rho$, 
a simple calculation gives
\begin{equation} \label{eA.5}
|\bar \nabla \rho| = |\nabla r|, 
\end{equation}
where the norm and gradient on the left are with respect to $\bar g$, and 
on the right are with respect to $g$.

 A defining function $\rho  = t$ is a geodesic defining function if
\begin{equation} \label{eA.6}
|\bar \nabla t| = 1, 
\end{equation}
in a collar neighborhood $U$ of $\partial M$. Clearly, the formulas 
(A.1)-(A.3) simplify considerably in this situation. The function $t$ 
is the distance function to $\partial M$ on $(M, \bar g)$, and 
similarly by (A.5), the function $r$ is a (signed) distance function on 
$(M, g)$. The integral curves of $\bar \nabla t$ and $\nabla r$ are 
geodesics in $(M, \bar g)$ and $(M, g)$ respectively.

 The $2^{\rm nd}$ fundamental form $\bar A$ of the level sets $S(t)$ of 
$t$ in $(M, \bar g)$ is given by $\bar A = D^{2}t$, with $\bar H = \bar 
\Delta t$ giving the mean curvature of $S(t)$. Along $t$-geodesics of 
$(M, \bar g),$ one has the standard Riccati equation
\begin{equation} \label{eA.7}
\bar H'  + |\bar A|^{2} + \bar Ric(\bar \nabla t, \bar \nabla t) = 0, 
\end{equation}
where $\bar H'  = d\bar H/dt$.

 The following formulas relating the curvatures of $(M, \bar g)$ at 
$\partial M$ to the intrinsic curvatures of $(\partial M, \gamma)$ may 
be found in [6, \S 1].

 Let $\bar g$ be a $C^{2}$ geodesic compactification of $(M, g)$, with 
$C^{2}$ boundary metric $\gamma$. Then at $\partial M$,
\begin{equation} \label{eA.8}
\bar s = 2n\bar Ric(N, N) = \frac{n}{n-1}s_{\gamma}, 
\end{equation}
where $N = \bar \nabla t$ is the unit normal to $\partial M$ with respect to 
$\bar g$. If $X$ is tangent to $\partial M$, then
\begin{equation} \label{eA.9}
\bar Ric(N, X) = 0, 
\end{equation}
while if $T$ denotes the projection onto $T(\partial M)$, then
\begin{equation} \label{eA.10}
(\bar Ric)^{T}= \frac{1}{n-2}((n-1)Ric_{\gamma} -  \frac{s_{\gamma}}{2(n-1)}\gamma) . 
\end{equation}
In particular, the full curvature of ambient metric $\bar g$ at 
$\partial M$ is determined by the curvature of the boundary metric 
$\gamma$. 

 Finally, we have the following formula for $\bar s'  = d\bar s/dt$:
\begin{equation} \label{eA.11}
\bar s'  = 2nt^{-1}|\bar D^{2}t|^{2} \geq  \frac{1}{2n^{2}}t\bar s^{2} \geq 0. 
\end{equation}
Hence, if $\bar s (0) > 0$, then 
\begin{equation} \label{eA.12}
t^{2} < 4n^{2}/ \bar s (0). 
\end{equation}

\end{appendix}

\bibliographystyle{plain}

\bigskip

\begin{center}
March, 2008
\end{center}
\medskip
\noindent
\address{Department of Mathematics\\
S.U.N.Y. at Stony Brook\\
Stony Brook, N.Y.11794-3651\\}
\email{anderson@math.sunysb.edu}

\end{document}